\documentclass[11pt]{article}
\usepackage[letterpaper,margin=1in]{geometry}
\usepackage[T1]{fontenc}
\usepackage{lmodern,microtype}
\usepackage{amsmath,amssymb,amsthm,mathtools,bm}
\usepackage{booktabs,array,longtable,graphicx,xcolor,enumitem}
\usepackage[colorlinks=true,linkcolor=blue!45!black,citecolor=blue!45!black,urlcolor=blue!45!black]{hyperref}
\usepackage[nameinlink,noabbrev]{cleveref}
\allowdisplaybreaks[2]
\setlist{itemsep=3pt,topsep=5pt}
\newtheorem{theorem}{Theorem}[section]
\newtheorem{lemma}[theorem]{Lemma}
\newtheorem{proposition}[theorem]{Proposition}
\newtheorem{corollary}[theorem]{Corollary}
\theoremstyle{definition}

\theoremstyle{remark}
\newtheorem{remark}[theorem]{Remark}
\newcommand{\E}{\mathbb E}
\newcommand{\Pp}{\mathbb P}
\newcommand{\F}{\mathcal F}
\newcommand{\R}{\mathbb R}
\newcommand{\as}{\mathrm{a.s.}}
\newcommand{\ind}{\mathbf 1}
\newcommand{\ee}{\mathrm e}
\newcommand{\Ps}{\mathcal P_s}
\newcommand{\Pm}{\mathcal P_{p,M}}

\DeclareMathOperator{\sign}{sign}
\numberwithin{equation}{section}
\title{Gaussian-efficient nonasymptotic confidence sets\\ for unbounded means}
\author{Aaditya Ramdas, Diego Martinez-Taboada\\
\texttt{aramdas@stanford.edu},  \texttt{diegomar@andrew.cmu.edu}}
\date{\today}
\begin{document}
\maketitle
\vspace{-1.5em}
\begin{abstract}
Consider the classical statistical problem of sub-Gaussian mean estimation under martingale dependence:
given $n$ observations with conditional mean $\mu$ and a
known conditional sub-Gaussian proxy $s>0$, we seek a nonasymptotically valid
confidence interval for $\mu$. The usual Chernoff interval has half-width
$s\sqrt{2\log(2/\alpha)/n}$; it neither adapts to the actual
standard deviation $\sigma\le s$ nor attains the Gaussian quantile
$z_{1-\alpha/2}$. We construct an analytically closed-form nonasymptotically-valid interval under martingale dependence whose
endpoints, at every iid law with $\sigma > 0$, agree almost surely to first
order with $\bar X_n\pm\sigma z_{1-\alpha/2}/\sqrt n$. The two new building blocks
are a sub-Gaussian exponential supermartingale with an \emph{empirical} quadratic variance and a
cubic correction term, and a sub-Gaussian continuation supermartingale, the former delivering variance adaptivity (``Empirical-Bernstein'' for unbounded data), the latter also providing the Gaussian quantile (``Empirical-Berry-Esseen'' for unbounded data).  
The same
constructions yield a closed-form time-uniform confidence sequence with asymptotically unimprovable half-width
$\sigma\sqrt{2\log\log n/n}$. We also give closed-form multivariate ellipsoids
with the true-covariance chi-square limit, and extend the scalar Gaussian and
LIL constructions to handle heavy-tailed settings where one only knows a bound on any absolute central moment of order
greater than two. Experiments show
substantial width reductions when the true variance is smaller than the implied variance bound. The surprising (to us) takeaway message is: given any more than two moments, one can construct fully data-dependent nonasymptotic confidence sets for the mean that are asymptotically efficient. We anticipate significant benefits to downstream applications in multi-armed bandits with unbounded data, and other settings where decisions were previously guided by non-adaptive (worst-case) concentration inequalities.
\end{abstract}

\section{Introduction}\label{sec:introduction}

Sub-Gaussian and related concentration bounds are classical, extensively
developed subjects. The exponential-moment method is associated with
Chernoff~\cite{concChernoff1952}; Bennett~\cite{concBennett1962} and
Hoeffding~\cite{concHoeffding1963} provide foundational inequalities for
independent sums, while Azuma~\cite{concAzuma1967} and
Freedman~\cite{concFreedman1975} develop martingale counterparts with bounded
increments and predictable-variance control. Broad treatments of exponential
concentration, sub-Gaussian variables, and their statistical applications are
available in the books of Ledoux~\cite{concLedoux2001}, Boucheron, Lugosi,
and Massart~\cite{concBLM2013}, Vershynin~\cite{concVershynin2018}, and
Wainwright~\cite{concWainwright2019}; de la Pe\~na, Lai, and
Shao~\cite{concDLP2009} give a systematic treatment of self-normalization,
including deviation inequalities and iterated-logarithm laws.

However, the confidence intervals resulting from inverting Chernoff bounds are asymptotically inefficient. 
For concreteness, take the sub-Gaussian mean estimation problem: given a sequence of conditionally $s$-sub-Gaussian data with mean $\mu$, it is well known that for every fixed $\lambda\in\R$, the natural exponential process
\[
 M_n^{\mathrm{MGF},\lambda}
 =\exp\left\{\lambda\sum_{i=1}^n(X_i-\mu)-\frac12ns^2\lambda^2\right\}
\]
is a nonnegative supermartingale, formed as the product of $n$ conditionally-valid e-values. Optimizing $\lambda$ and union bounding over positive and negative $\lambda$ gives the confidence interval $\bar X_n\pm s\sqrt{2\log(2/\alpha)/n}$. Its limitation
is visible in the exponent: the variance charge is the supplied proxy $s^2$,
irrespective of the observed variation. When $s$ is much larger than the actual
standard deviation, collecting more observations does not remove this
first-order inefficiency. Even if $s=\sigma$, the exponential-tail constant
$\sqrt{2\log(2/\alpha)}$ is larger than the Gaussian quantile
$z_{1-\alpha/2}$ at a fixed confidence level. So Chernoff bounds are loose in both these senses ($s,\alpha$).

Variance adaptation is also well established for \emph{bounded} random variables. Empirical-Bernstein bounds,
including those of Maurer and Pontil~\cite{concMaurerPontil2009}, replace a
worst-case range contribution by observed variation. Modern time-uniform
empirical-Bernstein methods~\cite{howard2021} and betting intervals for
bounded means~\cite{concWaudbySmithRamdas2024} sharpen this principle further.
The objectives here are
to extend these ideas to the stated \emph{unbounded} conditional-mean models,
analytically explicit reported intervals, and the exact true-variance
Gaussian quantile at fixed confidence or the sharp true-variance LIL
half-width. 

Our aim is to retain the finite-sample guarantee and interpretability of this
argument while replacing its deterministic variance charge by observed
variation. We choose the headline procedures in the main paper to have
\emph{analytically closed-form intervals or ellipsoids}, to be asymptotically
efficient at the relevant Gaussian or LIL benchmark, and to use simple,
computationally inexpensive statistics. Further, while we know that some of our one-round factors are inadmissible (Appendix~\ref{app:admissible}), our focus is on achieving asymptotic sharpness, and our methods were chosen to maximize intuition relative to the standard Chernoff constructions.

\subsection{Model, validity, and efficiency}
Let $X_i$ be $\F_i$-measurable, and suppose that an unknown constant $\mu$ and a
known $s>0$ satisfy, almost surely,
\begin{equation}\label{eq:model}
 \E[X_i\mid\F_{i-1}]=\mu,\qquad
 \E[\exp\{\theta(X_i-\mu)\}\mid\F_{i-1}]
 \le\exp\{\theta^2s^2/2\}\quad\text{for every }\theta\in\R.
\end{equation}
The observations need not be independent, identically distributed, symmetric,
or bounded. The proxy need not equal the true standard deviation. A confidence
interval $I_n=[L_n,U_n]$ is valid at a planned time $n$ if
$\Pp(\mu\in I_n)\ge1-\alpha$ under every such process. A confidence sequence
$(C_n)$ satisfies the stronger guarantee
$\Pp(\mu\in C_n\text{ for every }n\ge1)\ge1-\alpha$.

For efficiency only, specialize to prefixes of one fixed iid law with variance
$\sigma^2>0$, keeping $\alpha$ and the supplied envelope fixed. Our fixed-time
endpoint target is
\begin{equation}\label{eq:intro-GE}
 \sqrt n(\bar X_n-L_n)\to\sigma z_{1-\alpha/2},\qquad
 \sqrt n(U_n-\bar X_n)\to\sigma z_{1-\alpha/2}
 \quad\as.
\end{equation}
Our CS half-width target is
\begin{equation}\label{eq:intro-LIL}
 r_n\sim\sigma\sqrt{2\log\log n/n}\quad\as.
\end{equation}
These are different calibrations: a fixed-time Gaussian limit is not a
time-uniform claim. The vector target, for fixed dimension and positive definite
covariance $\Sigma$, is the ellipsoid centered at $\bar X_n$ with squared
Mahalanobis radius $\chi^2_{d,1-\alpha}/n$. The theorems require no known
positive lower bound on $\sigma^2$ or on the eigenvalues of $\Sigma$.
The constant conditional mean in \eqref{eq:model} is part of the model;
coverage of an average of varying conditional means is not asserted here.

\subsection{Prior work on Gaussian-efficient nonasymptotic intervals}
\label{subsec:dgg}

Gaussian-efficient finite-sample inference for unbounded observations predates
the present constructions. Recently, Derumigny, Girard, and
Guyonvarch~\cite{derumigny2025} consider the setting of iid
observations with unknown $\sigma^2>0$ and a known standardized kurtosis bound
$\E(X-\mu)^4/\sigma^4\le K$, their Propositions~3.4--3.5 and equation~(14)
provide a nonasymptotically valid interval whose endpoints agree to first order
with $\bar X_n\pm\sigma z_{1-\alpha/2}/\sqrt n$, using Berry-Esseen bounds.


Their assumptions differ from ours: a standardized kurtosis bound supplies no
bound on central moments, while a known loose sub-Gaussian proxy supplies no
uniform standardized kurtosis bound. Moreover, their fixed-time iid guarantee
does not cover arbitrary conditional-mean martingales or continuous
monitoring. These distinctions motivate our constructions. Despite there being no direct way to compare our methods since neither is more general than the other, initial experiments (that provide the exact kurtosis to theirs and fourth moment to ours) suggest that the methods are incomparable in general.

In contrast, there is significant work on the Gaussian efficiency in the bounded setting~\cite{romano2000,austern2022,ming2026gaffke,mtr2026}. 
Out of these, only the last work has martingale validity and iid efficiency, like the current work. In particular, we apply their betting strategy to the unbounded setting using particular nonnegative supermartingales that we introduce next.

\subsection{A brief summary of our new supermartingales}
Observed-quadratic self-normalization is a classical route to variance
adaptation. The results collected by Howard et al.~\cite{howard2020} include
exponential processes with realized quadratic variation and additional
predictable quadratic charges. Examining the proof of their result, also in Fan, Grama, and Liu~\cite{fan2015}, actually yields a better exponential supermartingale:
\[
 \exp\left\{\lambda n(\bar X_n-\mu)
 -\frac{\lambda^2}{2}\sum_{i=1}^n(X_i-\mu)^2
 -\gamma s^3n|\lambda|^3\right\},\qquad
 \gamma=\sqrt3e^{-3/2}.
\]
The main difference from the aforementioned works is that we retain the cubic term in $\lambda$, while the earlier ones converted it to a quadratic term in $\lambda$, which inflates the first-order variance coefficient.

While this version has the right local variance term, the quadratic term depends on $m$ and this does not yield a closed-form interval for $\mu$. Some deterministic inequalities yield a different supermartingale:
\begin{equation}\label{eq:intro-Qprocess}
 \exp\left\{\lambda n(\bar X_n-\mu)
 -\frac{\lambda^2}{2}(Q_n+s^2H_n)
 -s^3(n-H_n)|\lambda|^3\right\},\quad
 Q_n=\sum_{i=1}^n(X_i-\bar X_n)^2,\quad H_n=\sum_{i=1}^n\frac1i.
\end{equation}
For $|\lambda|\le1/(5s)$ this is a test supermartingale. Its inversion is
linear in the candidate mean. Its effective variance satisfies
$(Q_n+s^2H_n)/n\to\sigma^2$, while its cubic charge is negligible at the
relevant scales. These constructions are one core idea of the
paper.  An
asymptotically fine grid of bets gives \eqref{eq:intro-LIL} directly. 

But the above supermartingale does not give the
exact Gaussian quantile in \eqref{eq:intro-GE}. For that goal we use a very different robust \emph{Gaussian continuation} supermartingale
\[
 W_i(m)=\Phi\!\left(\frac{x_i(m)-b}{\sqrt{v_i(m)}}\right)
        +18\gamma s^3(N-i)\tau^{-3/2},
\]
where $(v_i(m),x_i(m))$ are the guarded states defined in
Section~\ref{subsec:ci-heat}, and $A,b,N,\tau$ are specified there.
The supermartingale property holds at $m=\mu$.
This develops the Gaussian continuation principle underlying the recent
bounded-data betting strategy of Martinez-Taboada and Ramdas~\cite{mtr2026} --- there it provided bets in a standard bounded martingale construction, here it directly provides a new nonnegative supermartingale. 
The vector construction uses a matrix heat potential
and a covariance-adaptive cubic mixture. Under a $p$th moment bound, a
power remainder replaces the cubic one and a coarse localization step retains
closed-form inversion.

\subsection{Paper Outline}
Section~2 proves \eqref{eq:intro-GE} for a scalar confidence interval under model
\eqref{eq:model}, with finite-sample coverage also valid conditionally on its
pilot. Section~3 gives a confidence sequence satisfying \eqref{eq:intro-LIL}.
Section~4 extends to multivariate data, and achieves the multivariate true-covariance chi-square limit.
Section~5 replaces sub-Gaussianity by a known conditional absolute central
$p$th moment for any $p>2$, and proves both scalar limits without a finite MGF (and thus this technique applies whenever a known two-sided MGF bound holds on a neighborhood of zero).
Every inference procedure in these four sections reports a closed-form interval
or ellipsoid --- simplicity and interpretability were preferred, even though some of the underlying one-round factors are inadmissible.
Appendix~A gives
explicit admissible improvements that dominate the studied e-values, and recalls a complete admissible-class characterization for moment-bounded laws. 
The natural MGF-based sub-Gaussian e-value is itself admissible, but we provide many others. 
These admissible improvements may not result in analytically closed-form intervals, which is
why the main constructions retain simpler exponential forms. Appendix~B uses
Gaffke's merger to give a different Gaussian-efficient interval for independent
observations, with numerical rather than analytic endpoints.

All method sections end with their own experiments. For the scalar CI with
$s=1$, at $n=1000$ the mean widths are about $0.77$ times Chernoff when
$\sigma=1/2$, and about $0.63$ times Chernoff when $\sigma=1/4$ across the
three reported confidence levels. They are slightly wider when
$\sigma=s$, which appears to be a small price relative to the large gains. For the scalar CS at $n=10000$, the width ratios to the standard
fixed-prior normal-mixture CS range from $0.48$ to $0.53$ when $\sigma=s/2$.
The moment experiments include Student distributions with no finite MGF.



\section{Variance-adaptive sub-Gaussian confidence intervals}\label{sec:subgaussian-ci}

We first derive the new exponential processes, explaining carefully the change in their quadratic terms from the standard construction.  Throughout, a test
supermartingale is nonnegative and starts at one. As in the standard Chernoff proof, all we do is obtain a particular e-value at time $n$, so Markov's inequality followed by rearrangement yields a confidence interval (usually one-sided, and one repeats the argument by negating the tuning parameter and taking a union bound). We later point out how to modify the argument for more general moment generating function (MGF) bounds.
We then give the explicit betting strategy needed for the exact Gaussian quantile.


\subsection{New exponential supermartingales}

Throughout, we denote
\[
 \xi_i:=X_i-\mu, 
 \quad S_n:=\sum_{i\le n}\xi_i,
 \quad T_n(\mu):=\sum_{i\le n}(X_i-\mu)^2, 
 \quad Q_n:=\sum_{i\le n}(X_i-\bar X_n)^2.
\]
The natural MGF process is thus
\[
 \exp\{\lambda S_n-ns^2\lambda^2/2\}.
\]
Its deterministic quadratic term cannot distinguish a law of variance $s^2$
from a law of much smaller variance. Our empirical constructions replace
that quadratic term, at least partly, by the realized squared residuals. In particular, the generic
second-moment methods collected in Howard et al.~\cite{howard2020},
Lemma~3(f,g) and Corollary~8, show how to do this with an additional quadratic
charge. Their ``cubic self-normalization'' result, Lemma~3(i) and Corollary~9, is particularly relevant.
Its proof, equation~(6.38), uses the sharper scalar inequality
\[
 \mathbb E_{i-1}\exp\{\lambda\xi_i-\lambda^2\xi_i^2/2\}
 \le \exp\{\lambda^3\mathbb E_{i-1}(\xi_i)_-^3/3\},\qquad\lambda\ge0.
\]
That proof credits Fan, Grama, and Liu~\cite{fan2015}, Corollary~2.2. The
sub-gamma formulation displayed in~\cite{howard2020} subsequently absorbs the third moments into the variance term, and thereby upper-bounds this cubic charge by a quadratic-order one; retaining the cubic order avoids the resulting first-order proxy inflation.

We keep the cubic order term, and we use sub-Gaussianity only
to bound its expectation in Lemma~\ref{lemma:cubic_inequality}. Subsequently, in Proposition~\ref{proposition:candidate_centered_supermartingale}, we build an exponential supermartingale whose quadratic term scales only with the variance, and the dependence on $s$ appears only in the cubic.


\begin{lemma}[Cubic inequality] \label{lemma:cubic_inequality}
For every real $u$,
\(
 e^{u-u^2/2}\le1+u+\frac{(u_-)^3}{3}.
\)
If $U$ is centered and $t$-sub-Gaussian, then
\[
 \mathbb E U^2\le t^2,\qquad
 \mathbb E U_-\le t/2,\qquad
 \mathbb E U_-^3\le3\gamma t^3,\qquad
 \gamma=\sqrt3e^{-3/2}<2/5.
\]
The same assertions hold conditional on the past if required.
\end{lemma}
\begin{proof}
For $u\ge0$, $\log(1+u)-u+u^2/2$ has derivative $u^2/(1+u)\ge0$.
For $u=-x\le0$, the function
$g(x)=1-x+x^3/3-e^{-x-x^2/2}$ satisfies $g(0)=g'(0)=0$ and
\[
 g''(x)=x\{2-(2+x)e^{-x-x^2/2}\}\ge0.
\]
The variance bound follows from the second-order MGF comparison, and centering
gives $\mathbb E U_-=\mathbb E|U|/2\le t/2$. For every $r>0$,
\[
 U_-^3\le(3/(er))^3e^{-rU}.
\]
Take expectations and choose $r=\sqrt3/t$; when $t=0$ the claim is immediate.
All steps are unchanged under conditioning.
\end{proof}

\begin{proposition}[Cubic self-normalized sub-Gaussian supermartingale] \label{proposition:candidate_centered_supermartingale}
For every fixed $\lambda\in\mathbb R$,
\[
 M_n^{\mathrm c,\lambda}(\mu)=
 \exp\{\lambda S_n -\lambda^2 T_n(\mu)/2
                   -\gamma s^3n|\lambda|^3\}
\]
is a test supermartingale (for the class of conditionally $s$-sub-Gaussian distributions with mean $\mu$). 
\end{proposition}
\begin{proof}
Conditionally, $\lambda\xi_i$ is centered and $s|\lambda|$-sub-Gaussian, so
\[
 \mathbb E_{i-1}e^{\lambda\xi_i-\lambda^2\xi_i^2/2}
 \le1+\gamma s^3|\lambda|^3\le e^{\gamma s^3|\lambda|^3}.
\]
The repaired factors therefore have conditional expectation at most one.
\end{proof}

With $M_n^{\mathrm c,\lambda}(m)$ defined by replacing $\mu$ by $m$,
this construction has the right variance coefficient, but its
inversion $\{m : M_n^{\mathrm c,\lambda}(m) < 1/\alpha\}$ is a valid $1-\alpha$ confidence set for $\mu$ but it need not be an interval. Indeed, $T_n(m) \gg |\bar X_n-m|$ for $|m| \gg 0$, so sufficiently
remote candidates $m$ can re-enter the acceptance set. 
It turns out that we can remove that candidate-dependent quadratic as follows. First, Lemma~\ref{lemma:cubic_inequality_uniform} extends Lemma~\ref{lemma:cubic_inequality} to shifted quadratic terms. Second, Theorem~\ref{thm:qprocess} uses such a cubic inequality to define a nonnegative supermartingale with similar properties to Proposition~\ref{proposition:candidate_centered_supermartingale}, but whose quadratic term does not depend on the specific candidate $m$.

\begin{lemma}[Shifted cubic inequality] \label{lemma:cubic_inequality_uniform}
Let $U$ be centered and $t$-sub-Gaussian, $0\le t\le1/5$. For every $a\in\mathbb R$,
\[
 \mathbb E e^{U-(U+a)^2/2}\le1+t^3\le e^{t^3}.
\]
Both $a$ and $t$ may be predictable in the conditional version.
\end{lemma}
\begin{proof}
For $a\ge0$, apply the preceding cubic inequality to $U+a$ to obtain
\[
 \mathbb E e^{U-(U+a)^2/2}
 \le e^{-a}\{1+a+\mathbb E(U+a)_-^3/3\}
 \le1+\gamma t^3\le1+t^3.
\]
For $a=-b<0$, dropping $-U^2/2$ and using the MGF bound gives
\[
 \mathbb E e^{U-(U-b)^2/2}
 \le\exp\{[t^2(1+b)^2-b^2]/2\}\le1
\]
if $b\ge t/(1-t)$. In the remaining case $0<b<t/(1-t)\le1/4$, use
$(U-b)_-\le U_-+b$ and the moment bounds to get
\[
 \mathbb E e^{U-(U-b)^2/2}
 \le e^b\left(1-b+\gamma t^3+bt^2+b^2t/2+b^3/3\right).
\]
Since $e^b\le e^{1/4}<4/3$, $\gamma<2/5$, and
$e^b(1-b)\le1-b^2/2$,
\begin{align*}
 \mathbb E e^{U-(U-b)^2/2}
 &\le1+\frac8{15}t^3+\frac43bt^2-\frac{23}{90}b^2\\
 &\le1+\frac8{15}t^3+\frac{40}{23}t^4
 \le1+\frac{304}{345}t^3\le1+t^3.
\end{align*}
Here the coefficient of $b^2$ is bounded by
$-1/2+(4/3)(1/10+1/12)=-23/90$, and the second line maximizes the
quadratic in $b$. For $t=0$ the claim follows directly from $U=0$ almost surely.
\end{proof}

\begin{theorem}[Another cubic self-normalized sub-Gaussian supermartingale]\label{thm:qprocess}
Let $H_n=\sum_{i=1}^ni^{-1}$. For every fixed $|\lambda|\le1/(5s)$,
\[
 M_n^\lambda=
 \exp\left\{\lambda S_n
       -\frac{\lambda^2}{2}(Q_n+s^2H_n)
       -(n-H_n)s^3|\lambda|^3\right\}
\]
is a test supermartingale, as is the sharper variant that replaces $s^3|\lambda|^3$ by $\log(1+s^3|\lambda|^3)$. 
\end{theorem}
\begin{proof}
For $i\ge2$, put $r_i=(i-1)/i$, $c_i=\bar X_{i-1}$,
$A_i=e^{\lambda\xi_i-\lambda^2(X_i-c_i)^2/2}$ and
$B_i=e^{\lambda\xi_i}$. Conditional H\"older interpolation gives
\[
 \mathbb E_{i-1}(A_i^{r_i}B_i^{1-r_i})
 \le (\mathbb E_{i-1}A_i)^{r_i}(\mathbb E_{i-1}B_i)^{1-r_i}
 \le\exp\{r_i s^3|\lambda|^3+(1-r_i)s^2\lambda^2/2\}.
\]
The first expectation is controlled by the shift-uniform lemma with
$U=\lambda\xi_i$ and $a=\lambda(\mu-c_i)$; the second is the natural MGF
bound. For $i=1$, use only the MGF bound, corresponding to $r_1=0$.
For the logarithmic version, retain
$\mathbb E_{i-1}A_i\le1+s^3|\lambda|^3$ instead of upper-bounding it by an
exponential; H\"older's power $r_i$ yields the charge
$r_i\log(1+s^3|\lambda|^3)$. Finally, Welford's identity
\[
 Q_n=\sum_{i=2}^n r_i(X_i-\bar X_{i-1})^2
\]
follows by expanding $Q_i-Q_{i-1}$, while
$\sum r_i=n-H_n$ and $\sum(1-r_i)=H_n$. Multiply the repaired factors.
\end{proof}

The quadratic variation is now independent of $m$. Set $L=\log(2/\alpha)$. 
Applying Markov's inequality to the nonnegative supermartingales (defined on $\bar X_n$ and $-\bar X_n$) at time $n$ for a positive $\lambda$ immediately gives the interval
\[
 \bar X_n\pm\left\{
   \frac{L}{n\lambda}+\frac{\lambda(Q_n+s^2H_n)}{2n}
   +s^3\left(1-\frac{H_n}{n}\right)\lambda^2\right\}.
\]

Our confidence sets hold under martingale dependence. This implies that these confidence sets can be applied to i.i.d. sub-Gaussian random variables, where the variance $\sigma^2 = \mathbb E\xi_i^2 = \mathbb E X_i^2 - \mu^2$ is fixed. However, the more general case $n^{-1}\sum_{i\le n}\mathbb E_{i-1}\xi_i^2\to\sigma^2>0$ is of especial interest.  

\begin{lemma}[Consistency of quadratic variation]
If $n^{-1}\sum_{i\le n}\mathbb E_{i-1}\xi_i^2\to\sigma^2>0$ almost surely,
then $Q_n/n\to\sigma^2$ almost surely.
\end{lemma}
\begin{proof}
Sub-Gaussian tails give $\mathbb E_{i-1}\xi_i^4\le16s^4$ and
$\mathbb E_{i-1}\xi_i^2\le s^2$. Hence both
$\sum_i\xi_i/i$ and
$\sum_i(\xi_i^2-\mathbb E_{i-1}\xi_i^2)/i$ are $L^2$-bounded martingales
and converge almost surely. Summation by parts (Kronecker's lemma) yields
$S_n/n\to0$ and
$n^{-1}\sum_i(\xi_i^2-\mathbb E_{i-1}\xi_i^2)\to0$. Use
$Q_n=\sum_i\xi_i^2-S_n^2/n$.
\end{proof}

\subsection{A sharp Empirical-Bernstein bound for sub-Gaussian data}
\label{subsec:ci-practical}

The first application of Theorem~\ref{thm:qprocess} is immediate: choose a
bet $\lambda$, apply Markov's inequality to its two signs, and solve a linear inequality
for the mean. 
For a planned sample size $n\ge4$, use a pilot of size
$k=\lfloor\sqrt n\rfloor$ and score the remaining $N=n-k$ observations,
$Z_i=X_{k+i}$. Write
\[
 \widehat v=\frac1k\sum_{i=1}^k(X_i-\bar X_k)^2,\qquad
 z=\frac1N\sum_{i=1}^N Z_i,\qquad
 Q=\sum_{i=1}^N(Z_i-z)^2,\qquad H=\sum_{i=1}^N\frac1i.
\]
The pilot chooses the bet $\lambda$; the remaining data determines the center and
observed variance. 

At error level $\delta\in(0,1)$, set $L=\log(2/\delta)$. Every
pilot-based $0<\lambda\le1/(5s)$ gives the interval
\begin{equation}\label{eq:ci-Qradius}
 \boxed{\ I_Q(\delta)=z\pm r_Q(\delta),\qquad
 r_Q(\delta)=\underbrace{\frac{L}{N\lambda}}_{\text{evidence threshold}}
 +\underbrace{\frac{\lambda(Q+s^2H)}{2N}}_{\text{observed variance}}
 +\underbrace{s^3\left(1-\frac HN\right)\lambda^2}_{\text{cubic repair}}.\ }
\end{equation}
Here $z\pm r$ denotes $[z-r,z+r]$. A simple default is
\begin{equation}\label{eq:ci-simplebet}
 \lambda=\min\left\{\frac1{5s},\
       \sqrt{\frac{2L}{N(\widehat v+s^2/k)}}\right\}.
\end{equation}
This balances the first two terms of the radius after substituting the
pilot variance for $Q/N$. The regularizer $s^2/k$ keeps the formula defined
when the pilot scatter is zero, and the cap enforces the theorem's allowed
bet size. These are fixed choices, rather than parameters to optimize.

\begin{proposition}[Coverage and variance adaptation]\label{prop:ci-Qsimple}
For any allowed pilot-measurable bet,
\begin{equation}\label{eq:ci-Qcoverage}
 \Pp\{\mu\in I_Q(\delta)\mid\F_k\}\ge1-\delta.
\end{equation}
At a fixed iid law with variance $\sigma^2>0$, the simple choice
\eqref{eq:ci-simplebet}, at fixed $\delta$, satisfies almost surely
\begin{equation}\label{eq:ci-Qlimit}
 \sqrt n\,r_Q(\delta)\longrightarrow
       \sigma\sqrt{2\log(2/\delta)},\qquad
 \sqrt n(z-\bar X_n)\longrightarrow0.
\end{equation}
\end{proposition}
\begin{proof}
Condition on the pilot and apply Theorem~\ref{thm:qprocess} to the scoring
sample. Markov's inequality bounds each signed terminal exponent at level
$L$ with error $\delta/2$. Inverting and taking a union bound gives
\eqref{eq:ci-Qcoverage}. For the limit, $\widehat v\to\sigma^2$,
$Q/N\to\sigma^2$, and $N/n\to1$ almost surely, so eventually
$\lambda\sim\sqrt{2L/(N\sigma^2)}$. The harmonic and cubic corrections
are $o(n^{-1/2})$. Finally,
$\sqrt n(z-\bar X_n)=(k/\sqrt n)(z-\bar X_k)\to0$ by the strong law.
\end{proof}

Thus the supplied proxy $s$ disappears from the leading width, just as the
empirical quadratic in Section~2.1 suggests. But to move from the Chernoff limit involving $\sqrt{\log (2/\delta)}$ to the Gaussian quantile, we need a new nonnegative supermartingale, introduced soon. But we first point out that the same construction above yields an empirical Bernstein inequality much beyond the sub-Gaussian class: whenever a known two-sided MGF bound holds on a neighborhood of zero.

\subsection{A sharp Empirical-Bernstein inequality under a local MGF bound}
\label{subsec:subpsi-Q}

Let $(X_i)_{i\ge1}$ be a real-valued process adapted to a filtration
$(\mathcal F_i)_{i\ge0}$, with a common conditional mean
$\mathbb E[X_i\mid\mathcal F_{i-1}]=\mu$. Write
$\xi_i=X_i-\mu$ and $\mathbb E_{i-1}=\mathbb E[\,\cdot\mid\mathcal F_{i-1}]$.
Suppose a known deterministic function $K=v\psi$, with known $v>0$, satisfies
\begin{equation}\label{eq:subpsi-Q-assumption}
 \mathbb E_{i-1}\exp(u\xi_i)\le \exp\{K(u)\},
 \qquad u\in[-R,R],\quad i\ge1,
\end{equation}
for some known $R>0$, where $K$ is finite on $[-R,R]$.
Define
\begin{equation}\label{eq:subpsi-Q-constants}
 A_R=\exp\{K(R)\}+\exp\{K(-R)\}-2,\qquad
 M_3=\frac{\sqrt3\,A_R}{R^3}.
\end{equation}
If $A_R>0$, set $S^2=2A_R/R^2$; if $A_R=0$, choose any $S>0$.
Put
\begin{align}
 \lambda_0&=\min\left\{\frac1{5S},\frac{R}{1+RS}\right\},
       \label{eq:subpsi-Q-domain}\\
 \rho(\lambda)&=\log\left(1+\frac{4M_3}{9}|\lambda|^3
                          +\frac{40S^4}{23}\lambda^4\right).
       \label{eq:subpsi-Q-repair}
\end{align}

\begin{theorem}[A cubic-repaired sub-$\psi$ sample-scatter supermartingale]
\label{thm:subpsi-Q}
Under \eqref{eq:subpsi-Q-assumption}, 
for every fixed $\lambda\in[-\lambda_0,\lambda_0]$, the process
$M_0^{\lambda}=1$ and
\begin{equation}\label{eq:subpsi-Q-process}
 M_n^\lambda=\exp\left\{
 \lambda\sum_{i=1}^n(X_i-\mu)-\frac{\lambda^2}{2}Q_n
 -H_nK(\lambda)-(n-H_n)\rho(\lambda)\right\},\qquad n\ge1,
\end{equation}
is a nonnegative supermartingale. In particular, with
\begin{equation}\label{eq:subpsi-Q-pureconstant}
 C=\frac{4M_3}{9}+\frac{8S^3}{23},
\end{equation}
the process obtained by replacing $\rho(\lambda)$ in
\eqref{eq:subpsi-Q-process} by $C|\lambda|^3$ is also a nonnegative
supermartingale.
In either version, $K(\lambda)$ may be replaced by the sharper envelope
\begin{equation}\label{eq:subpsi-Q-localenvelope}
 K_R(\lambda)=\min\{K(\lambda),S^2\lambda^2/2\}.
\end{equation}
\end{theorem}

\begin{remark}[Constants and the role of two-sided MGF control]
The constants in \eqref{eq:subpsi-Q-constants} use only two known MGF
envelope values. If a sharper local quadratic bound is available, one
may instead use any $S>0$ satisfying
$K(u)\le S^2u^2/2$ for $|u|\le R$, together with any known
conditional bound $\mathbb E_{i-1}|\xi_i|^3\le M_3$.
The rest of the theorem and proof are unchanged.
\end{remark}

\begin{remark}[Observed variance and local quadratic behavior]
Using $K_R$ ensures $0\le K_R(\lambda)\le S^2\lambda^2/2$ even
if the original envelope is loose near zero.
For example, under iid sampling with variance $\sigma^2>0$, for each
fixed real $z$ and $\lambda_n=z/\sqrt n$, almost surely
\[
 \frac{\lambda_n^2Q_n}{2}\longrightarrow\frac{z^2\sigma^2}{2},
 \qquad H_nK_R(\lambda_n)\longrightarrow0,\qquad
 (n-H_n)\rho(\lambda_n)\longrightarrow0.
\]
This identifies the true-variance local quadratic term. 
\end{remark}

The resulting empirical Bernstein confidence interval adapts to the unknown
variance. For a planned $n\ge4$ and fixed $\delta\in(0,1)$, use the pilot
and scoring notation of the previous subsection, set $L=\log(2/\delta)$, and take
\[
 \lambda_n=\min\left\{\lambda_0,
       \sqrt{\frac{2L}{N(\widehat v+S^2/k)}}\right\},\qquad
 r_{\psi,n}(\delta)=\frac{L}{N\lambda_n}
       +\frac{\lambda_n(Q+S^2H)}{2N}
       +\left(1-\frac HN\right)\frac{\rho(\lambda_n)}{\lambda_n}.
\]
Using $K_R(\pm\lambda_n)\le S^2\lambda_n^2/2$, Theorem~\ref{thm:subpsi-Q}
gives conditional coverage at least $1-\delta$ for $z\pm r_{\psi,n}(\delta)$
given the pilot. At a fixed iid law with $\sigma^2>0$, its half-width satisfies
\[
 r_{\psi,n}(\delta)=\sigma\sqrt{\frac{2\log(2/\delta)}n}
       +o(n^{-1/2})\qquad\text{almost surely}.
\]

\subsection{Gaussian-efficiency using a robust Gaussian continuation supermartingale}
\label{subsec:ci-heat}

Variance adaptation and Gaussian calibration are separate tasks. For every
positive bet, arithmetic--geometric mean gives
\[
 r_Q(\delta)\ge\frac{\sqrt{2\log(2/\delta)(Q+s^2H)}}{N}.
\]
Consequently, choosing a better fixed bet cannot replace the Chernoff
constant by $z_{1-\delta/2}$. We need a test designed around a Gaussian
terminal event, following the design principle of~\cite{mtr2026}.

The design comes from the Doob martingale of a Gaussian rejection event.
For a terminal cutoff $b$, define the continuation probability
\begin{equation}\label{eq:ci-heat-main-potential}
 F(v,x)=\Phi\!\left(\frac{x-b}{\sqrt v}\right).
\end{equation}
This is the probability that a Gaussian continuation, starting at sum $x$
with remaining variance $v$, finishes above $b$. To see its role, suppose
temporarily that the scored observations are independent
$N(\mu,\sigma^2)$ variables with known $\mu$ and $\sigma^2>0$. For $i<N$,
\begin{equation}\label{eq:ci-heat-gaussian-doob}
 \begin{aligned}
 x_i&=\sum_{j=1}^i(Z_j-\mu),\\
 F((N-i)\sigma^2,x_i)
 &=\E\!\left[\ind\!\left\{\sum_{j=1}^N(Z_j-\mu)>b\right\}
       \,\middle|\,Z_1,\ldots,Z_i\right].
 \end{aligned}
\end{equation}
At $i=N$, define the value to be the rejection indicator itself. These
conditional probabilities form its Doob martingale, a bounded nonnegative
martingale. It starts at the Gaussian rejection probability and ends at
zero or one according to whether the terminal test rejects. Tracking this
event is what targets the Gaussian cutoff. In this setting, applying Markov's inequality and solving for $\mu$ yields exactly the classic Gaussian confidence interval for the mean (technically one-sided, but you can repeat with $-x_i$ and take a union bound for the two-sided case). 

Our aim is to robustify this martingale to handle sub-Gaussian observations while
also learning their variance. In the Gaussian experiment the remaining
variance is known exactly as $(N-i)\sigma^2$. In our construction the pilot
sets an initial variance budget, and each centered observation $u$ spends
its observed square: update $(v,x)$ to $(v-u^2,x+u)$. This empirical update
needs a drift correction, even for Gaussian observations. The heat identity
$F_v=F_{xx}/2$ gives
\[
 F(v-u^2,x+u)-F(v,x)
 =F_x(v,x)u+
   \underbrace{\left(\tfrac12F_{xx}-F_v\right)}_{0}u^2
   +\text{remainder}.
\]
Conditional centering removes the linear term. The quadratic cancellation
is the same principle as in Section~2.1. The potential can still have a
small positive drift from the remainder, so we add a reserve that pays for
it at each observation.

Here is the resulting nonnegative supermartingale. For a candidate mean
$m$, fix one sign and start at $(v_0,x_0)=(A,0)$, with $A\ge2\tau>0$, for some tuning parameters $A,\tau$ that we specify later.
On scoring observation $i$, propose
\[
 (v_i,x_i)=(v_{i-1}-u_i^2,x_{i-1}+u_i),\qquad
 u_i=\pm(Z_i-m).
\]
Accept the update if its variance is at least $\tau$. Otherwise keep the
previous state and freeze it for all remaining observations. Also freeze
after an accepted update leaves variance below $2\tau$. These guards keep
the Gaussian continuation smooth even for an unusually large observation.
The states depend on $m$ and the fixed sign chosen for $u_i$; we suppress that dependence
in the notation. For each choice of sign, the process is
\begin{equation}\label{eq:ci-heat-main-supermartingale}
 \boxed{\ W_i(m)=
 \underbrace{\Phi\!\left(\frac{x_i-b}{\sqrt{v_i}}\right)}_{
       \text{Gaussian continuation probability}}
 +\underbrace{\frac{18\gamma s^3(N-i)}{\tau^{3/2}}}_{
       \text{remaining drift reserve}},\qquad 0\le i\le N.\ }
\end{equation}
\begin{theorem}[Gaussian continuation supermartingale]
\label{thm:ci-heat-supermartingale}
Condition on the pilot and choose $A\ge2\tau>0$ and $b\in\R$ using only
the pilot. Under the conditional sub-Gaussian model, for each fixed sign,
the process \eqref{eq:ci-heat-main-supermartingale} at $m=\mu$ is a
nonnegative supermartingale with respect to $(\F_{k+i})_{i=0}^N$:
\[
 \E[W_i(\mu)\mid\F_{k+i-1}]\le W_{i-1}(\mu),\qquad 1\le i\le N.
\]
Dividing by $W_0>0$ gives a test supermartingale starting at one.
\end{theorem}
The proof is in Appendix~\ref{app:scalar-heat}. The reserve decreases at
each step to pay the potential's possible upward drift. Thus $W_i$ is the
Gaussian continuation
probability plus the budget for its remaining approximation error.
This robustification gives validity throughout the conditional sub-Gaussian
class, while estimating variance through the pilot and observed scatter.

Here is the complete sub-Gaussian recipe at error level
$\varepsilon\in(0,1)$. All slack is tied to the single prescribed sequence
\[
 \rho=\frac1{8\log(n+e)}.
\]
The three quantities that enter the construction are
\begin{align}
 R&=\frac{s}{\sqrt{2\varepsilon N\rho}},
 &\tau&=s^2\left\{\frac{36\gamma N}
                 {\varepsilon\rho(1-2\rho)}\right\}^{2/3},
 \label{eq:ci-heat-simple-reserves}\\
 A&=N(\widehat v+8s^2\rho)+NR^2+2\tau,
 &&\gamma=\sqrt3e^{-3/2}.
 \label{eq:ci-heat-simple-budget}
\end{align}
Each has a specific role. The coarse radius $R$ localizes the candidate mean.
The reserve $\tau$ keeps the continuation smooth. The initial variance budget
$A$ pays for the pilot variance, its vanishing slack, localization, and
smoothing. None is a freely tuned hyperparameter.

Choose the terminal cutoff
\[
 b=\sqrt A\,z_{1-\varepsilon(1-2\rho)^2/2}.
\]
This fixes the initial Gaussian tail probability. With the reserve above,
\[
 W_0=\frac{\varepsilon(1-2\rho)(1-\rho)}2,\qquad
 W_N(m)=F(v_N,x_N).
\]
Conditional Markov's inequality therefore gives, for each sign,
\begin{equation}\label{eq:ci-heat-main-tail}
 \Pp\{F(v_N,x_N)\ge1-\rho\mid\F_k\}
 \le\frac{\varepsilon(1-2\rho)}2
 \qquad\text{at }m=\mu.
\end{equation}
The remaining error $2\varepsilon\rho$ pays for the coarse radius $R$.

If the observed variance fits the pilot budget,
$Q/N\le\widehat v+8s^2\rho$, set
\begin{equation}\label{eq:ci-heat-simple-radius}
 \boxed{\ r_H(\varepsilon)=\min\left\{R,\
 \frac{\sqrt A\,z_{1-\varepsilon(1-2\rho)^2/2}
       +\sqrt{A-Q}\,z_{1-\rho}}{N}\right\}.\ }
\end{equation}
Otherwise set $r_H(\varepsilon)=R$, and report the interval $I_H(\varepsilon)=z\pm r_H(\varepsilon)$. The check ensures $A-Q>0$ whenever
the square root is used. The first term in the numerator is the Gaussian
terminal cutoff; the second pays for the variance left in the continuation.

\begin{theorem}[Gaussian calibration]\label{thm:ci-heat-main}
Under the conditional sub-Gaussian model, for every pilot-measurable
$\varepsilon\in(0,1)$, $I_H(\varepsilon)=z\pm r_H(\varepsilon)$ has conditional coverage at least
$1-\varepsilon$ given the pilot. At a fixed iid law with $\sigma^2>0$,
if $\varepsilon_n\to\varepsilon_\infty\in(0,1)$, then almost surely
\begin{equation}\label{eq:ci-heat-main-limit}
 \sqrt n\,r_H(\varepsilon_n)\longrightarrow
           \sigma z_{1-\varepsilon_\infty/2},\qquad
 \sqrt n(z-\bar X_n)\longrightarrow0.
\end{equation}
\end{theorem}
\begin{proof}[Why the calibration works]
Conditional Chebyshev gives
$\Pp\{|z-\mu|>R\mid\F_k\}\le2\varepsilon\rho$.
On the observable variance check, every candidate with $|z-m|\le R$ obeys
\[
 \sum_{i=1}^N(Z_i-m)^2=Q+N(z-m)^2\le A-2\tau.
\]
Its process never freezes, so
$x_N=\pm N(z-m)$ and $v_N=A-Q-N(z-m)^2$.
If such a candidate lies beyond the Gaussian radius in
\eqref{eq:ci-heat-simple-radius}, one of its signed terminal potentials
exceeds $1-\rho$. The two bounds in \eqref{eq:ci-heat-main-tail}, together
with localization, give total error
$\varepsilon(1-2\rho)+2\varepsilon\rho=\varepsilon$.
If the check fails, the coarse radius applies directly.

For efficiency, the pilot and scoring variances converge at polynomial
almost-sure rates. Therefore the variance check eventually passes and
\[
 A/N\to\sigma^2,\qquad (A-Q)/N\sim8s^2\rho,\qquad
 \sqrt{\rho}\,z_{1-\rho}\to0.
\]
The leftover-variance term vanishes after multiplication by $\sqrt N$,
while the cutoff tends to $\sigma z_{1-\varepsilon_\infty/2}$.
The coarse radius has $\sqrt N R\to\infty$, so it does not change this
limit. The pilot changes the center by $o(n^{-1/2})$ as before.
\end{proof}

The general moment formulation and the two technical lemmas needed for
this argument appear together in Appendix~\ref{app:scalar-heat}. The
Gaussian calibration uses only a moment of order greater than two; the
exponential interval above additionally uses the sub-Gaussian assumption.

\subsection{Relation to martingale Berry--Esseen bounds}
\label{rem:heat-vs-martingale-BE}
The Gaussian continuation argument belongs to the classical proof tradition
of martingale normal approximation. Bolthausen~\cite{Bolthausen1982}
established foundational Berry--Esseen bounds for martingale difference
arrays. A particularly explicit comparison is Fan~\cite{Fan2019},
Theorems~2.1--2.2 and Section~3.1, equations~(33)--(40): his proof introduces the Gaussian
potential $\Phi((b-x)/\sqrt v)$, telescopes its increments, and uses
$F_{xx}=2F_v$ to cancel the leading quadratic terms. Our potential
$F(v,x)=\Phi((x-b)/\sqrt v)$ uses the complementary tail and the same
identity. Conditional centering removes its linear term, and a conditional
moment of order $q\in(2,3]$ controls the remainder. The heat-kernel
cancellation itself is therefore classical.

Our construction makes the normal-approximation error explicit at
each step. If $\mathbb E_{i-1}|Z_i-\mu|^q\le M$, the guarded update with
variance floor $\tau$ has conditional upward drift at most
$d=3M\tau^{-q/2}$. Thus, for correctly centered data and for each fixed sign,
\[
 W_i=F(v_i,x_i)+(N-i)d
\]
is a nonnegative supermartingale. Conditional on the pilot, with initial
variance budget $A$ and $x_0=0$, its terminal implication is
\[
 \mathbb P\{F(v_N,x_N)\ge1-\eta\mid\text{pilot}\}
 \le\frac{\Phi(-b/\sqrt A)+Nd}{1-\eta},\qquad 0<\eta<1.
\]
In this precise sense, the Gaussian branch is a Berry--Esseen-type
correction of a Gaussian test, organized as a test supermartingale after
normalization by $W_0>0$.

The inferential statement is nevertheless different from directly
applying a martingale Berry--Esseen theorem. Such theorems typically bound
the uniform Kolmogorov distance between a normalized terminal martingale
and a Gaussian law, with assumptions or error terms controlling its
\emph{predictable} quadratic variation. Here the variance budget is spent on the
realized squares $(Z_i-m)^2$,  an observable variance audit,
and a coarse fallback yielding finite-sample coverage without assuming
variance stabilization. We control a guarded, smoothed rejection event;
we do not prove a uniform normal-approximation bound for a Studentized
statistic under arbitrary adaptation. Gaussian efficiency is a separate
conclusion under iid sampling. 

Accordingly, our robust Gaussian continuation supermartingale combines some classical martingale
Berry--Esseen machinery with a particular observed-variance implementation
and explicit confidence-interval inversion. 

\subsection{Combining practical width with Gaussian efficiency}
\label{subsec:ci-hybrid}

The exponential radius is useful at moderate sample sizes; the Gaussian
radius supplies the exact limiting quantile. Both are centered at $z$, so
combining them requires only an error split and a minimum. Set
\begin{equation}\label{eq:ci-budgets}
 \alpha_Q=\frac{10\alpha}{10+\log(n+e)},\qquad
 \alpha_H=\alpha-\alpha_Q,
\end{equation}
and report
\begin{equation}\label{eq:ci-final}
 \boxed{\ I_n=z\pm\min\{r_Q(\alpha_Q),r_H(\alpha_H)\}.\ }
\end{equation}
Use \eqref{eq:ci-simplebet} for the exponential bet, or the optional
refinement in Appendix~\ref{app:scalar-tuning}. For $n<4$, use the ordinary
Chernoff interval. The particular constant ten only controls the
finite-sample error allocation; any fixed positive constant gives the same
conclusions below.

\begin{theorem}[Validity and endpoint Gaussian efficiency]
\label{thm:ci-final}
For every adapted sequence in the sub-Gaussian model,
$\Pp(\mu\in I_n)\ge1-\alpha$. At every fixed iid law with variance
$\sigma^2>0$, both endpoints agree almost surely to first order with
$\bar X_n\pm\sigma z_{1-\alpha/2}/\sqrt n$.
\end{theorem}
\begin{proof}
The intersection has error at most $\alpha_Q+\alpha_H=\alpha$ by a union
bound. For efficiency, $\alpha_H\to\alpha$, so
Theorem~\ref{thm:ci-heat-main} gives the desired limit for the Gaussian
branch. Since $\alpha_Q\to0$, the exponential branch satisfies, for every
allowed positive bet,
\[
 \sqrt N\,r_Q(\alpha_Q)
 \ge\sqrt{2\log(2/\alpha_Q)\,Q/N}\longrightarrow\infty
 \quad\as.
\]
Thus the minimum eventually selects the Gaussian branch. The pilot center
differs from $\bar X_n$ by $o(n^{-1/2})$ almost surely.
\end{proof}

This combination preserves the two guarantees, but it makes no claim that
the conservative Gaussian calibration is useful at moderate $n$. The
experiments below distinguish the source of the observed gains from the
source of the eventual Gaussian constant.

\subsection{Empirical comparison with the Chernoff interval}
\label{subsec:ci-experiments}

We evaluated \eqref{eq:ci-final} with the optional finite-sample tuning
\eqref{eq:ci-pilottuning}--\eqref{eq:ci-newton} in
Appendix~\ref{app:scalar-tuning}, including the pilot and the error allocation
\eqref{eq:ci-budgets}. The displayed numerical results concern this refined
bet, rather than the simpler default \eqref{eq:ci-simplebet}. Each of six laws and each
$n\in\{100,300,1000,3000,10000\}$ used $5000$ independent samples; the three
confidence levels were evaluated on the same samples. The random seed was
$2026102801$. We supplied the common known proxy $s=1$ in every case.
All means were zero. The laws were three Gaussians with standard deviations
$0.25,0.5,1$; a symmetric two-point law on $\{-0.5,0.5\}$; a uniform law on
$[-\sqrt{3}/2,\sqrt{3}/2]$; and $2(B-0.1)$ for $B\sim\mathrm{Bernoulli}(0.1)$.
The bounded laws have range at most two and are therefore $1$-sub-Gaussian.
We compared width against the standard Chernoff interval using all $n$
observations and its full requested error probability $\alpha$.

\begin{table}[t]
\centering
\caption{Mean width divided by full-sample Chernoff width. Values below one
favor the proposed closed-form interval. Each entry averages $5000$ samples.}
\label{tab:ci-fresh}
\small
\begin{tabular}{lrrrrrr}
\toprule
&\multicolumn{3}{c}{$n=1000$}&\multicolumn{3}{c}{$n=10000$}\\
\cmidrule(lr){2-4}\cmidrule(lr){5-7}
Law&$90\%$&$95\%$&$99\%$&$90\%$&$95\%$&$99\%$\\
\midrule
$N(0,0.25^2)$&0.630&0.637&0.654&0.468&0.469&0.476\\
$N(0,0.5^2)$&0.767&0.768&0.776&0.650&0.646&0.644\\
$N(0,1)$&1.189&1.179&1.169&1.140&1.123&1.102\\
Two-point, $\sigma=0.5$&0.765&0.766&0.775&0.650&0.647&0.645\\
Uniform, $\sigma=0.5$&0.765&0.767&0.775&0.650&0.647&0.644\\
$2(B-0.1)$, $\sigma=0.6$&0.845&0.844&0.849&0.738&0.731&0.725\\
\bottomrule
\end{tabular}
\end{table}

At $n=1000$, the interval reduced mean widths by about $35$--$37\%$ when
$\sigma=0.25$ and by about $22$--$24\%$ when $\sigma=0.5$. The gains grew
at $n=10000$ and persisted for the asymmetric law. These gains require some
slack between the known proxy and the true standard deviation: the tight-proxy
Gaussian retained a $10$--$19\%$ width overhead in the displayed range. At
$n=100$, width ratios across the complete design ranged from $0.969$ to
$1.408$, so this construction does not uniformly improve on Chernoff at small
sample sizes. The largest Monte Carlo standard error of a mean ratio across
all designs was below $0.00085$.

Observed coverage was conservative. Across all laws and sample sizes, the
smallest coverage estimates were $0.9946$, $0.9972$, and $0.9998$ at nominal
$90\%$, $95\%$, and $99\%$ confidence, respectively. These observations are
diagnostics, not substitutes for the coverage proof. The heat branch was
selected in none of the $450000$ evaluated intervals. Thus the reported
moderate-sample gains come from the $Q_N$ exponential construction, while the
heat branch supplies the exact Gaussian limiting constant. Improving its
finite-sample remainder calibration remains an identifiable practical goal.

\section{Variance-adaptive sub-Gaussian LIL confidence sequences}
\label{sec:subgaussian-cs}

A confidence interval pays for one planned look at the data. A confidence
sequence pays for every look, and its sharp benchmark is correspondingly
different: under iid sampling with variance $\sigma^2>0$, the law of the
iterated logarithm (LIL) has scale
$\sigma\sqrt{2\log\log n/n}$. The earlier exponential supermartingale provides an
especially direct route to that scale. 

\subsection{One exponential process, many simultaneous bets}
\label{subsec:cs-grid}
Recall $S_n(m)=n(\bar X_n-m)$, $Q_n=\sum_{i\le n}(X_i-\bar X_n)^2$,
and $H_n=\sum_{i\le n}1/i$ and that
Theorem~\ref{thm:qprocess} gives the test supermartingale
\begin{equation}\label{eq:cs-q-process}
 M_n^\lambda=
 \exp\!\left\{\lambda S_n(\mu)
  -\frac{\lambda^2}{2}(Q_n+s^2H_n)
  -(n-H_n)\log(1+s^3|\lambda|^3)\right\},
 \qquad |\lambda|\le b:=\frac1{5s},
\end{equation}
For a fixed positive bet $\lambda$ and threshold $L$, define the explicit
radius
\begin{equation}\label{eq:cs-fixed-radius}
 B_n(\lambda,L)=\frac{L}{n\lambda}
 +\frac{\lambda(Q_n+s^2H_n)}{2n}
 +\frac{n-H_n}{n\lambda}\log(1+s^3\lambda^3).
\end{equation}
If the two processes with bets $\lambda$ and $-\lambda$ stay below
$\exp(L)$, then $|\bar X_n-\mu|\le B_n(\lambda,L)$ simultaneously at all
times. Ville's inequality makes this a simultaneous guarantee with error
at most $2e^{-L}$. Thus the remaining issue is to choose a useful bet without
paying a first-order penalty for that choice.

Following the calibrated-grid approach to confidence sequences~\cite{howard2021},
use the countable grid
\begin{align}
 m_j&=\lceil\log(j+e)\rceil,
 &w_j&=(\log2)\left\{\frac1{\log(j+2)}-\frac1{\log(j+3)}\right\} \asymp \frac1{j(\log j)^2},
 \label{eq:cs-grid-weights}\\
 \lambda_{j,k}&=b\,2^{-j-k/m_j},
 &\pi_{j,k}&=\frac{w_j}{m_j},
 \quad j\ge0,\quad 0\le k<m_j.
 \label{eq:cs-grid-bets}
\end{align}
The $w_j$ telescope to one, so $\sum_{j,k}\pi_{j,k}=1$. The dyadic grid
becomes progressively finer as the bets approach zero. This refinement matters:
a grid with a fixed multiplicative spacing can leave a fixed inflation in the
LIL constant.
For a current-data selection of one grid point, set
\begin{align}
 V_n&=Q_n+s^2H_n,
 &g_n&=\log\log(e^e+n)+\log\frac{20}{\alpha},\label{eq:cs-selector1}\\
 \lambda_n^*&=\min\left\{b,
       \sqrt{\frac{2g_n}{V_n}},
       \left(\frac{g_n}{2s^3n}\right)^{1/3}\right\},
 &q_n&=\log_2\frac b{\lambda_n^*},\label{eq:cs-selector2}\\
 k_n&=\lfloor m_{\lfloor q_n\rfloor}(q_n-\lfloor q_n\rfloor)\rfloor.
 \label{eq:cs-selector3}
\end{align}
Write $\lambda_n=\lambda_{\lfloor q_n\rfloor,k_n}$ and
\begin{equation}\label{eq:cs-grid-threshold}
 L_n=\log\frac{10}{\alpha\pi_{\lfloor q_n\rfloor,k_n}}.
\end{equation}
All quantities are defined at $n=1$, since $V_n>0$. The cube-root cap protects
against a tiny observed variance at early times; it is eventually inactive at
every positive-variance iid law. No search over candidates or over the grid is
performed.

\subsection{A closed-form sub-Gaussian confidence sequence with a practical coarse grid}
\label{subsec:cs-closedform}
The asymptotically fine grid is not always the best finite-sample choice.
We therefore reserve equal mass for four coarse bets and the fine grid:
\begin{equation}\label{eq:cs-main-radius}
 r_n^{\mathrm{CS}}
 =\min\left\{
    B_n(\lambda_n,L_n),\quad
    \min_{0\le\ell\le3}
     B_n\left(b2^{-\ell},\log\frac{10}{\alpha}\right)
        \right\},
 \qquad
 C_n=\left[\bar X_n-r_n^{\mathrm{CS}},\bar X_n+r_n^{\mathrm{CS}}\right].
\end{equation}
This is the minimum of five displayed scalar expressions, all sharing the
same center. It requires only the running mean, scatter, and harmonic number.
In particular, it is analytically closed-form and takes a constant number of
real-arithmetic operations and constant storage per observation.

\begin{theorem}[Finite-sample coverage]\label{thm:cs-coverage}
Under the conditional sub-Gaussian model, the intervals
\eqref{eq:cs-main-radius} satisfy
\[
 \mathbb P\{\mu\in C_n\text{ for every }n\ge1\}\ge1-\alpha.
\]
No independence or variance convergence is required for this coverage claim.
\end{theorem}
\begin{proof}
For each coarse bet, apply Ville's inequality separately to its two signs
with error $\alpha/10$ per sign. For each fine-grid bet $(j,k)$, apply it
with error $\alpha\pi_{j,k}/10$ per sign. The union bound costs
$4\alpha/5+\alpha/5=\alpha$. On the complementary event, at every time
and for every bet in both collections,
\[
 n\lambda|\bar X_n-\mu|
 \le L+\frac{\lambda^2V_n}{2}
      +(n-H_n)\log(1+s^3\lambda^3).
\]
Dividing by $n\lambda$ proves every radius entering the minimum. Selection
of $(\lfloor q_n\rfloor,k_n)$ using the current observations is allowed because this one
event already covers every grid point at every time. This is not an
uncalibrated substitution of a data-dependent bet into one fixed-bet process.
\end{proof}

The four coarse bets are a finite-sample convenience, not an asymptotic
requirement. One may omit them and give the fine grid all of $\alpha$, replacing
$10$ by $2$ in \eqref{eq:cs-grid-threshold}. More generally, any predetermined
finite list of bets and positive error allocations can be added in this way.
The version in \eqref{eq:cs-main-radius} is our fixed default and is the one
evaluated below.

\subsection{Exact adaptation to the true variance at the LIL scale}
\label{subsec:cs-efficiency}
Suppose, in addition to conditional sub-Gaussianity, that
\begin{equation}\label{eq:cs-variance-limit}
 \frac1n\sum_{i=1}^n
 \mathbb E[(X_i-\mu)^2\mid\mathcal F_{i-1}]
 \longrightarrow\sigma^2\in(0,s^2]
 \quad\text{almost surely}.
\end{equation}
In particular, this holds for iid observations with variance $\sigma^2>0$.

\begin{theorem}[Sharp LIL half-width]\label{thm:cs-lil}
Under \eqref{eq:cs-variance-limit}, for each fixed $\alpha\in(0,1)$,
\begin{equation}\label{eq:cs-lil-limit}
 \frac{r_n^{\mathrm{CS}}}
      {\sigma\sqrt{2\log\log n/n}}
 \longrightarrow1\qquad\text{almost surely}.
\end{equation}
\end{theorem}
\begin{proof}
Write $\xi_i=X_i-\mu$ and
$v_i=\mathbb E[\xi_i^2\mid\mathcal F_{i-1}]$. Conditional sub-Gaussianity
implies $\mathbb E\xi_i^2\le s^2$ and $\mathbb E\xi_i^4\le16s^4$. Therefore
the martingales $\sum_{i\le n}\xi_i/i$ and
$\sum_{i\le n}(\xi_i^2-v_i)/i$ are bounded in $L^2$ and converge almost
surely. Kronecker's lemma gives
\[
 \bar X_n\to\mu,\qquad
 \frac1n\sum_{i\le n}(\xi_i^2-v_i)\to0
 \quad\text{almost surely}.
\]
Using \eqref{eq:cs-variance-limit} and
$Q_n/n=n^{-1}\sum_{i\le n}\xi_i^2-(\bar X_n-\mu)^2$ yields
$Q_n/n\to\sigma^2$. Since $H_n=O(\log n)$, we also have
$V_n/n\to\sigma^2$.

On this probability-one event, put $h_n=\log\log n$. Then
$g_n\sim h_n$, and eventually the square-root term is the smallest in
\eqref{eq:cs-selector2}: it tends to zero, and its ratio to the cube-root
term is $O((h_n/n)^{1/6})\to0$. Hence
\[
 \lambda_n^*\sim\sqrt{\frac{2h_n}{n\sigma^2}},\qquad
 \lfloor q_n\rfloor\sim\frac{\log n}{2\log2},\qquad m_{\lfloor q_n\rfloor}\to\infty.
\]
The upward rounding obeys
$1\le\lambda_n/\lambda_n^*<2^{1/m_{\lfloor q_n\rfloor}}$, so the ratio tends to one.
The mean value theorem gives
\[
 w_j\sim\frac{\log2}{j(\log j)^2},\qquad m_j\sim\log j,
\]
and consequently
\[
 L_n=\log\log n+3\log\log\log n+O(1)\sim h_n.
\]
The first two terms of \eqref{eq:cs-fixed-radius} are therefore each
asymptotic to $\frac12\sigma\sqrt{2h_n/n}$. The remaining term is at most
$s^3\lambda_n^2=O(h_n/n)=o(\sqrt{h_n/n})$. Thus the fine-grid radius has
the limit in \eqref{eq:cs-lil-limit}.

For each fixed coarse bet $\lambda>0$, its radius tends to the positive
constant
$\lambda\sigma^2/2+\log(1+s^3\lambda^3)/\lambda$. The fine-grid radius
tends to zero, so the minimum in \eqref{eq:cs-main-radius} eventually selects
it. This proves the claimed limit for the practical default.
\end{proof}

At an iid nondegenerate law, the classical LIL states
\[
 \limsup_{n\to\infty}
 \frac{|\bar X_n-\mu|}{\sigma\sqrt{2\log\log n/n}}=1
 \quad\text{almost surely}.
\]
A sample-mean-centered confidence sequence whose half-width converges to a
strictly smaller multiple would therefore miss $\mu$ infinitely often almost
surely. This is the sense in which the leading constant is optimal. 

\subsection{Empirical results}
\label{subsec:cs-experiments}
We simulated $2{,}000$ independent paths in each of four designs:
$N(0,\sigma^2)$ and the symmetric two-point law on $\{-\sigma,\sigma\}$,
with $\sigma\in\{1/2,1\}$ and the reported proxy fixed at $s=1$. Each path
was followed at every integer time through $10{,}000$. The exact method tested
was \eqref{eq:cs-main-radius}, including the four coarse bets and logarithmic
repair, with no tuning by the realized data beyond the displayed grid
selector. The baseline was the standard normal-mixture confidence sequence
\begin{equation}\label{eq:cs-normalmixture-baseline}
 r_n^{\mathrm{NM}}=\frac{s}{n}
 \sqrt{(n+1)\{\log(n+1)+2\log(1/\alpha)\}},
\end{equation}
which mixes the natural MGF processes against a centered normal prior of
precision $s^2$.

\begin{table}[ht]
\centering
\caption{Mean half-width of the proposed confidence sequence divided by the
normal-mixture half-width \eqref{eq:cs-normalmixture-baseline}, for Gaussian
observations. Values below one favor the proposed method. The same ratios
for the symmetric two-point laws differ by at most $0.001$ after rounding.}
\label{tab:cs-widths}
\begin{tabular}{ccrrr}
\toprule
$\sigma/s$ & $1-\alpha$ & $n=100$ & $n=1{,}000$ & $n=10{,}000$\\
\midrule
$1/2$ & $90\%$ & $0.976$ & $0.641$ & $0.484$\\
$1/2$ & $95\%$ & $1.016$ & $0.667$ & $0.497$\\
$1/2$ & $99\%$ & $1.105$ & $0.723$ & $0.526$\\
$1$ & $90\%$ & $1.220$ & $0.990$ & $0.849$\\
$1$ & $95\%$ & $1.244$ & $0.996$ & $0.880$\\
$1$ & $99\%$ & $1.305$ & $1.018$ & $0.950$\\
\bottomrule
\end{tabular}
\end{table}

At $n=1{,}000$, variance adaptation reduces mean widths by $28$--$36\%$
when $\sigma=s/2$; at $n=10{,}000$ the reductions are $47$--$52\%$.
When the proxy is tight, the proposed widths are comparable by $n=1{,}000$
and smaller by $5$--$15\%$ at $n=10{,}000$. At $n=100$ the proposed
construction is generally wider, particularly when $\sigma=s$. The largest
Monte Carlo standard error of a Gaussian width ratio was below $0.0011$.

For Gaussian observations with $\sigma=s$, the numbers of paths on which
the proposed confidence sequence ever missed the mean were $43$, $15$, and
$1$ out of $2{,}000$ at confidence levels $90\%$, $95\%$, and $99\%$;
the corresponding normal-mixture counts were $129$, $61$, and $12$.
For the two-point law with $\sigma=s$, these counts were $50$, $24$, and
$4$, versus $100$, $54$, and $14$ for normal mixtures. When $\sigma=s/2$,
the proposed method had at most three misses per $2{,}000$ paths across the
tested levels, and the normal mixture had none. These finite-horizon
simulations are checks of implementation and finite-sample behavior, not
proofs of the simultaneous guarantee or of the LIL limit.

The comparison fixes the normal-mixture prior in
\eqref{eq:cs-normalmixture-baseline}. A different prior scale can improve
its widths near a chosen horizon, so the table is not a claim of uniform
finite-sample dominance over all normal-mixture choices. The theoretical
advantage is adaptation to $\sigma$ and the exact $\log\log n$ leading rate;
the coarse bets improve its onset without changing either guarantee.

\section{Covariance-adaptive multivariate sub-Gaussian means}\label{sec:multi}

Let $X_i\in\mathbb R^d$ have the common conditional mean $\mu$ and satisfy
\begin{equation}\label{eq:multi-model}
 \E[\exp\{u^\top(X_i-\mu)\}\mid\mathcal F_{i-1}]
 \le\exp\{s^2\|u\|^2/2\},\qquad u\in\mathbb R^d.
\end{equation}
The dimension is fixed in the efficiency statements. The target is the usual
true-covariance Gaussian ellipsoid: at each fixed iid law with covariance
$\Sigma\succ0$, the reported ellipsoid $\mathcal E_n$ should satisfy
\begin{equation}\label{eq:multi-goal}
 \sqrt n(\mathcal E_n-\bar X_n)
 \longrightarrow
 \{\Sigma^{1/2}z:\|z\|^2\le\chi^2_{d,1-\alpha}\}
 \quad\text{almost surely in Hausdorff distance}.
\end{equation}
The covariance condition is needed only for this efficiency conclusion.

We give a construction with two explicit components. A Gaussian mixture of the
directional cubic supermartingales gives useful covariance adaptation at
moderate sample sizes. A vector version of the additive heat construction
supplies the exact chi-square limit. Both components return explicit ellipsoids using only
matrix operations and simple formulas.


One benchmark from prior work is the usual normal-mixture ball, centered at the sample
mean, whose prior scale is optimized to the number of observations. For $N$ observations and error level $b$, it has Euclidean radius
\begin{equation}\label{eq:multi-normal-radius}
 w(N,d,b,s)=\frac{s}{\sqrt N}\min_{\gamma>0}
 \sqrt{(1+\gamma)\{2\log(1/b)+d\log(1+1/\gamma)\}}.
\end{equation}
The minimum has an analytic formula, given in
Appendix~\ref{app:multi-normal}.
Thus the benchmark already receives its best isotropic prior scale,
while its radius still uses the supplied proxy $s$ in every direction.
The ball has coverage at least $1-b$ under \eqref{eq:multi-model},
including conditionally after a pilot~\cite{abbasi2011}. 

It is possible to take the exponential supermartingale constructions from earlier, apply them after projecting the data onto a unit vector $\theta$, and take a (say Gaussian) mixture over $\theta$, in order to derive an empirical-Bernstein multivariate confidence ellipsoid. The details are in Appendix~\ref{app:multi-mixture}; we move on to directly develop the chi-squared-efficient construction.

\subsection{A matrix heat ellipsoid with an explicit chi-squared cutoff}
\label{subsec:multi-heat}

To obtain the exact chi-square constant, use the Gaussian terminal-event
strategy of Section~\ref{subsec:ci-heat}. The pilot covariance $B$ sets the
coordinates $B^{-1/2}(z-m)$: a ball in these coordinates becomes an
ellipsoid in the original coordinates. In the ideal Gaussian experiment
with known covariance, the conditional probability of rejecting outside
a terminal chi-square ball is a nonnegative Doob martingale. We mimic it
using the estimated coordinates and an observed covariance clock.

At error level $b$, write $$c_d(t)=\sqrt{\chi^2_{d,1-t}}$$ 
and use the same
slack $\rho=1/[8\log(n+e)]$ as before. The explicit calibration in
Appendix~\ref{app:multi-heat} chooses an initial covariance budget $AI_d$,
a smoothing floor $\tau I_d$, and a coarse localization radius $W_G$,
all fixed given the pilot. The continuation probability is
\begin{equation}\label{eq:multi-heat-potential}
 F(V,x)=\Pp\!\left\{\|x+V^{1/2}Z\|>
             \sqrt A\,c_d\bigl(b(1-2\rho)^2\bigr)\right\},
 \qquad Z\sim N_d(0,I_d).
\end{equation}
Starting at $(V_0,x_0)=(AI_d,0)$, each whitened observation
$u_j=B^{-1/2}(X_{k+j}-m)$ adds $u_j$ to the sum and subtracts
$u_ju_j^\top$ from the remaining covariance. Accept an update only if
the proposed covariance is at least $\tau I_d$; otherwise retain the
previous state and freeze. Also freeze after an accepted update leaves
its smallest eigenvalue below $2\tau$. With these guards, the central
process is
\begin{equation}\label{eq:multi-heat-process}
 \boxed{\ W_j(m)=F(V_j,x_j)
      +b\rho(1-2\rho)\left(1-\frac jN\right),\qquad 0\le j\le N.\ }
\end{equation}
The second term is the unused drift reserve: it decreases linearly to
zero, including after freezing. The matrix heat identity cancels the
quadratic drift, and the calibration makes this reserve pay the remainder.
Consequently $W_j(\mu)$ is a nonnegative supermartingale. This correction
is needed for the empirical covariance clock even with Gaussian data.

To report an ellipsoid, compute the observed scatter in the pilot's
coordinates, $Q_B=B^{-1/2}QB^{-1/2}$. The check
$\lambda_{\max}(Q_B)\le N(1+8\rho)$ ensures that covariance spending fits
the budget in every direction for all localized candidates. When it
passes, report
\begin{equation}\label{eq:multi-heat-radius}
 \begin{aligned}
 r_G&=\frac{\sqrt A\,c_d\bigl(b(1-2\rho)^2\bigr)
      +\sqrt{A-\lambda_{\min}(Q_B)}\,c_d(\rho)}N,\\
 \mathcal G_b&=\{m:(m-z)^\top B^{-1}(m-z)\le r_G^2\}.
 \end{aligned}
\end{equation}
The first radius term is the Gaussian cutoff; the second pays for the
remaining covariance. If the check fails, return the radius-$W_G$ ball
centered at $z$. The computation uses matrix operations and chi-square
quantiles, without evaluating $F$. Coverage requires no accuracy guarantee
for the pilot covariance.

\begin{theorem}[Validity and efficiency of the matrix heat ellipsoid]
\label{thm:multi-heat}
With the calibration in Appendix~\ref{app:multi-heat}, conditional on the
pilot, $W_j(\mu)$ is a nonnegative supermartingale and $\mathcal G_b$ has
miscoverage at most $b$ under \eqref{eq:multi-model}. At every fixed iid
law with $\Sigma\succ0$, if $b\to\alpha\in(0,1)$, the check eventually
passes almost surely and
\[
 B\to\Sigma,\qquad \sqrt N r_G\to c_d(\alpha),
\]
thus satisfying~\eqref{eq:multi-goal}.
\end{theorem}
The pilot supplies the limiting covariance shape, and the continuation
supplies the chi-square size. The leftover-covariance term vanishes at
the leading scale. Calibration details and the proof are in
Appendix~\ref{app:multi-heat}.

\subsection{Adaptive linear bandits: efficient confidence sets and sequences}
\label{subsec:linear-bandits}

The same ideas extend to adaptive regression, after replacing centered
observations by regression scores. This step matters: the conditional mean
of $x_iY_i$ changes with the chosen action, whereas the score
$x_i(Y_i-x_i^\top\theta_*)$ is conditionally centered. We therefore apply
the empirical-variance constructions directly to these scores.

Let $x_i\in\R^d$ be chosen using $\F_{i-1}$, with $\|x_i\|\le L$ for a
known $L>0$, and suppose
\begin{equation}\label{eq:lb-model}
 Y_i=x_i^\top\theta_*+\varepsilon_i,\qquad
 \E[e^{u\varepsilon_i}\mid\F_{i-1}]\le e^{u^2s^2/2},\quad u\in\R.
\end{equation}
The unknown parameter needs no known norm bound. For any finite predictable
forecast $\widetilde\theta_{i-1}$, define the score and its observed
variance by
\begin{equation}\label{eq:lb-observed-scores}
 \begin{aligned}
 S_n(\theta)&=\sum_{i\le n}x_i(Y_i-x_i^\top\theta),\\
 R_n&=\sum_{i\le n}x_ix_i^\top
                    (Y_i-x_i^\top\widetilde\theta_{i-1})^2.
 \end{aligned}
\end{equation}
The forecast is chosen from past observations, so the variance is computable
and shared across candidate parameters. Its accuracy is unnecessary for coverage. The
forecast-uniform inequality from Section~2.1 gives the central process
\begin{equation}\label{eq:lb-forecast-process}
 M_n^u(\theta)=\exp\left\{u^\top S_n(\theta)
       -\tfrac12u^\top R_nu-s^3\sum_{i\le n}|u^\top x_i|^3\right\},
 \qquad \|u\|\le(5sL)^{-1}.
\end{equation}
For each fixed $u$, this is a test supermartingale at $\theta=\theta_*$.
Its terms have the same roles as before: directional evidence, observed
score variation, and a cubic repair. The clock $R_n$ is shared across
candidates, which allows explicit ellipsoidal inversion.

For the efficiency statements, let
$D_n=\sum_{i\le n}x_ix_i^\top$ describe the design and
$K_n=\sum_{i\le n}\E[\varepsilon_i^2\mid\F_{i-1}]x_ix_i^\top$
describe the variance of the score. Write $\widehat\theta_n$ for ordinary
least squares on all $n$ observations whenever $D_n$ is nonsingular.
The familiar sandwich covariance combines these two matrices: score
uncertainty is transformed into parameter uncertainty by the inverse
design.

\subsubsection{Gaussian efficiency at a planned time}
At a planned sample size, use the Gaussian terminal-event calibration of
Section~\ref{subsec:multi-heat} on whitened regression scores. A pilot
estimates their covariance; the matrix heat potential and its declining
reserve give the validity guarantee. Mapping the resulting score ball
through the inverse design gives a parameter ellipsoid. A coarse
score-mixture set supplies the fallback when the observed design or
variance budget is insufficient. The complete recipe is in
Appendix~\ref{app:lb-fixed}.

\begin{theorem}[Gaussian-efficient linear-bandit confidence ellipsoid]
\label{thm:lb-fixed}
Under \eqref{eq:lb-model}, the explicit set $\mathcal L_{n,b}$ in
Appendix~\ref{app:lb-fixed} satisfies
$\Pp(\theta_*\in\mathcal L_{n,b})\ge1-b$ at every planned time,
for arbitrary predictable bounded actions. For efficiency, suppose
$D_n/n\to D_\infty$ and $K_n/n\to K_\infty$ almost surely at a polynomial
rate, with deterministic positive-definite limits, as specified in
\eqref{eq:lb-stability}. Then, for fixed $b\in(0,1)$,
\begin{equation}\label{eq:lb-gaussian-limit}
 \sqrt n(\mathcal L_{n,b}-\widehat\theta_n)
 \longrightarrow
 \{z:z^\top(D_\infty^{-1}K_\infty D_\infty^{-1})^{-1}z
                   \le\chi^2_{d,1-b}\}
\end{equation}
almost surely in Hausdorff distance.
\end{theorem}
This is the usual Gaussian regression ellipsoid with the true sandwich
covariance. With homoskedastic noise, $K_\infty=\sigma^2D_\infty$, so its
scale uses the true noise variance $\sigma^2$. Stability is needed for this
efficiency claim; coverage itself permits arbitrary predictable designs.
The pilot is asymptotically negligible, so the limiting center is
full-sample least squares. The proof is in Appendix~\ref{app:lb-fixed}.

\subsubsection{Continuous monitoring at the LIL scale}
For continuous monitoring, apply Ville's inequality to simultaneous
directional bets from \eqref{eq:lb-forecast-process}, following the grid
principle of Section~3. The simultaneous guarantee permits the reported
ellipsoid to use the currently observed covariance. The matrix grids,
bet selection, and explicit inversion are given in
Appendix~\ref{app:lb-cs}; no grid is enumerated by the algorithm.

\begin{theorem}[A true-variance linear-bandit confidence sequence]
\label{thm:lb-cs}
Under \eqref{eq:lb-model}, the sets $\mathcal C_n^{\rm LIL}$ in
Appendix~\ref{app:lb-cs} satisfy
$\Pp(\theta_*\in\mathcal C_n^{\rm LIL}\text{ for every }n)\ge1-\alpha$
for arbitrary predictable bounded actions and finite predictable forecasts.
For efficiency, suppose $D_n/n\to D_\infty\succ0$ and
$K_n/n\to K_\infty\succ0$ almost surely, with deterministic limits, and
$\widetilde\theta_n\to\theta_*$ almost surely. For every fixed nonzero
$a\in\R^d$,
\begin{equation}\label{eq:lb-cs-efficient}
 \sup_{\theta\in\mathcal C_n^{\rm LIL}}
       |a^\top(\theta-\widehat\theta_n)|
 \sim\sqrt{\frac{2\log\log n}{n}
                         a^\top D_\infty^{-1}K_\infty D_\infty^{-1}a}
 \quad\as.
\end{equation}
\end{theorem}
Thus the confidence sequence has the same covariance geometry, with the
LIL cost for monitoring at every time. Its efficiency needs convergence
and a consistent forecast, without the polynomial rate used by the
planned-time theorem. Predictable ridge least-squares forecasts provide
one consistent choice under the stated design convergence. The proof is
in Appendix~\ref{app:lb-cs}.

For practical sample sizes, supplement the sharp construction with the
normal-mixture score set and a cheap forecast expert, allocating error
across the candidates and reporting the ellipsoid of smallest volume.
The experiments use these practical variants; their precise allocations
and the planned-time hybrid are in Appendix~\ref{app:lb-practical}.
At singular designs the planned-time fallback may be a cylinder, and the
sequential sets report $\R^d$. These results provide confidence
certificates and sharp limiting shapes; no improved bandit regret theorem
is claimed.

\section{A known moment of any order above two suffices}
\label{sec:moments}


The variance-adaptive principle does not require an MGF. Replace
sub-Gaussianity by a known conditional central-moment bound:
\begin{equation}\label{eq:mp-model}
 \E[X_i\mid\F_{i-1}]=\mu,\qquad
 \E[|X_i-\mu|^p\mid\F_{i-1}]\le M_p,\qquad p>2.
\end{equation}
Both $p$ and $M_p$ are known. If $M_p=0$, all observations equal $\mu$
almost surely, and the singleton $\{X_1\}$ solves both inference tasks.
Otherwise set $q=\min(p,3)$ and
$\kappa=3^{(1-q)/2}M_p^{q/p}$. The same observed-quadratic bet now needs a
power repair:
\begin{equation}\label{eq:mp-sm}
 M_t^\lambda=\exp\!\left\{
 \lambda\sum_{i=1}^t(X_i-\mu)
 -\frac{\lambda^2}{2}\sum_{i=1}^t(X_i-\mu)^2
 -\kappa t|\lambda|^q\right\}.
\end{equation}
which is a test supermartingale for every fixed real
$\lambda$; predictable bets are also allowed. The conditional moment
bound pays the remainder left after the linear and quadratic terms cancel.
Appendix~\ref{app:mp-process} gives the scalar inequality and proof, and Appendix~\ref{app:mp-obstruction} points out why we in general cannot construct a different supermartingale which does not have $\mu$ in the quadratic term, as we had done in the sub-Gaussian (and sub-$\psi$) cases. Nevertheless, one can still derive intervals from the above supermartingale, just not via a direct inversion; one has to perform an additional localization step, as we explain below but leave the details to the appendices.

The important fact is $q>2$: as the bets shrink, $|\lambda|^q$ is of
smaller order than $\lambda^2$. At a positive-variance iid law, this
makes the repair negligible at the Gaussian and LIL betting scales. The
known moment bound controls validity without fixing the leading width to
its possibly loose envelope.

\subsection{An analytically closed-form Gaussian-efficient interval}
\label{sec:mp-ci}

At a planned sample size, use a small pilot to choose the bet and score
the remaining observations, as in Section~\ref{subsec:ci-practical}.
There are again two calibrations. The exponential process
\eqref{eq:mp-sm} gives the practical interval, and the Gaussian
continuation gives the exact limiting quantile. The continuation argument
needs only a moment of order above two, so the general result in
Appendix~\ref{app:scalar-heat} applies directly.

One difference is essential: the quadratic in \eqref{eq:mp-sm} still
depends on the candidate mean. A coarse, moment-valid interval first
localizes that candidate; an explicit quadratic formula then gives the
practical radius. Intersect it with the Gaussian-calibrated interval,
allocating the total error between them. Both intervals have the same
scoring-sample center, so their intersection simply takes the smaller
radius. All allocations, reserves, and the quadratic formula are in
Appendix~\ref{app:mp-ci}; no candidate search is required.

\begin{theorem}[Moment-bound Gaussian efficiency]\label{thm:mp-ci}
Under \eqref{eq:mp-model}, the interval $I_n^p=[L_n,U_n]$ constructed in
Appendix~\ref{app:mp-ci} has coverage at least $1-\alpha$ for every $n$.
For $n\ge4$, coverage also holds conditionally on the pilot. At every
fixed iid law with variance $\sigma^2>0$, its endpoints satisfy
\[
 \sqrt n(\bar X_n-L_n)\longrightarrow\sigma z_{1-\alpha/2},\qquad
 \sqrt n(U_n-\bar X_n)\longrightarrow\sigma z_{1-\alpha/2}
 \quad\as.
\]
\end{theorem}
Under iid sampling, the pilot and scoring variances learn $\sigma^2$, while the extra moment
pays for the approximation error. The Gaussian branch can be conservative
at moderate sample sizes; the practical branch supplies the gains observed
in the experiments. The proof is in Appendix~\ref{app:mp-ci}.

\subsection{A closed-form confidence sequence with the sharp LIL rate}
\label{sec:mp-cs}

For continuous monitoring, apply Ville's inequality simultaneously to a
calibrated grid of bets from \eqref{eq:mp-sm}, following Section~3.
A coarse moment-valid confidence sequence localizes the true mean at
every time. On this one common event, the same quadratic inversion gives
valid radii for every bet, so the reported interval may select a useful
bet from the current data. A slowly refining grid attains the exact LIL
constant; four fixed bets improve its finite-sample behavior. The grid,
selection rule, and explicit radii are in Appendix~\ref{app:mp-cs}.

\begin{theorem}[Moment-bound LIL confidence sequence]\label{thm:mp-cs}
Under \eqref{eq:mp-model}, the intervals
$C_n^p=[\bar X_n-r_n^p,\bar X_n+r_n^p]$ in
Appendix~\ref{app:mp-cs} satisfy
$\Pp(\mu\in C_n^p\text{ for every }n)\ge1-\alpha$.
At every fixed iid law with variance $\sigma^2>0$,
\begin{equation}\label{eq:mp-lil}
 \frac{r_n^p}{\sigma\sqrt{2\log\log n/n}}\longrightarrow1
 \quad\as.
\end{equation}
\end{theorem}
The construction uses constant storage and a fixed number of elementary
operations per observation. The proof is in Appendix~\ref{app:mp-cs}.

\section{Conclusion}\label{sec:conclusion}

The natural sub-Gaussian MGF process pays for the supplied sub-Gaussian parameter at every step, and leaves itself no room to adapt to more favorable data.
Observed-quadratic self-normalization changes that cost to the realized
variation. Retaining a remainder of order strictly larger than two lets this
change survive at the asymptotic Gaussian and LIL scales. For sub-Gaussian (or sub-$\psi$)
data, we presented an especially
transparent empirical Chernoff construction.

A calibrated grid
of exponential bets yields the true-variance LIL constant through a
direct closed-form inversion. The exact fixed-confidence Gaussian quantile
requires a different terminal-event supermartingale; a bounded Gaussian heat potential
provides one with explicit error accounting. The same principle gives
true-covariance chi-square ellipsoids in fixed dimension. A known central
moment of any order above two suffices for the scalar Gaussian and LIL limits,
even without any finite MGF. Finite-sample validity holds in the stated
conditional-mean models; iid assumptions enter only the efficiency arguments.

The main procedures deliberately prioritize explicit endpoints, interpretable
variance terms, and inexpensive computation. Admissible one-round factors
can improve the underlying simple factors, but their power gains
come at the price of possibly giving up simplicity, interpretability and the interval geometry.

In experiments, the scalar sub-Gaussian CI improves substantially on Chernoff when the known
proxy is loose, while it retains a slight overhead for a tight proxy. The
confidence sequence has a later but clear advantage over the reported
normal-mixture baseline. 
The numerical Gaffke
construction in Appendix~B can be sharper at moderate samples under
independence, at the cost of numerical endpoint computation and the stronger
sampling assumption of independence.

It is abundantly clear to us that this paper will have many  downstream applications in adaptive decision making algorithms (like multi-armed bandits and other sequential experiments) with unbounded data, that are often guided by non-adaptive Chernoff or heavy-tailed concentration inequalities. 
Of course, we anticipate that future work may design slightly more practical variants of our constructions for these goals.
We view our primary contribution as conceptual: simply pointing to the very possibility that one can derive sharp, fully variance-adaptive Chernoff bounds, that behave like asymptotically efficient Gaussian bounds up to lower order corrections. We think that the form of our new cubic-corrected and robust Gaussian potential nonnegative supermartingales are likely of independent and lasting interest to the field of concentration (yielding both confidence intervals and sequences), even if the precise correction terms may vary in the future.

\paragraph{AI Statement.} Both authors interacted extensively with GPT to create this document. The most important building blocks of this paper are the authors' prior works on martingale concentration, variance adaptivity and Gaussian efficiency in bounded settings, such as \cite{howard2021}, \cite{concWaudbySmithRamdas2024} and \cite{mtr2026}. The very possibility of achieving sharp empirical Bernstein bounds or Gaussian efficient intervals for unbounded data is not a frequently asked question in the literature. While the key  constructions were provided by the authors, GPT was essential in identifying the exact presented form of the supermartingales, many of which benefited from fine-tuning of constants and tedious tracking of lower-order error terms for proving Gaussian efficiency. The authors iterated on each section separately, and within each section strove to find conceptually simple constructions that worked. Many interesting ideas and constructions ultimately had to be omitted for brevity and clarity. After the content of each section was in place, the authors iterated on the writing, first with GPT and then by hand, until the entire narrative was sufficiently clean to the authors' tastes.

\appendix

\section{Admissible sub-Gaussian and moment-class e-values}\label{app:admissible}

The main paper keeps exponential forms because their quadratic clocks and
outer inversions are easy to interpret and compute. Some of the charged
one-round factors used to derive them are inadmissible: a different valid
factor can dominate them at every observation. This appendix constructs such
improvements and proves admissibility for the specified null classes.
These are one-round claims; they do not establish admissibility or
inadmissibility of an entire adaptive betting policy, confidence sequence,
or interval. In particular, we do not infer the admissibility status of the
forecast-based $Q_n$ process from that of a candidate-centered factor.

For a candidate mean, write $Y=X-m$ and consider
\begin{align}
 \Ps&=\{P:\E_PY=0,\ \E_Pe^{tY}\le e^{s^2t^2/2}\text{ for all }t\in\R\},
 \qquad s>0,\label{adm:eq:sgclass}\\
 \Pm&=\{P:\E_PY=0,\ \E_P|Y|^p\le M\},
 \qquad p>1,\ M>0.\label{adm:eq:momentclass}
\end{align}
An e-value is a Borel measurable $e:\R\to[0,\infty]$ with
$\sup_{P\in\mathcal P}\E_Pe\le1$. It is \emph{pointwise admissible} if no
other valid $g$ satisfies $g\ge e$ everywhere with a strict gain at even one
point. All competing functions may be merely measurable. The moment-class
characterization below is a special case of the general finite-constraint
theory of Larsson, Ramdas, and Ruf~\cite{lrr2026}; we give a direct proof.

\subsection{Two elementary admissibility principles}\label{adm:sec:principles}

\begin{lemma}[Exact and approximate witness laws]\label{adm:lem:witness}
Let $e$ be an e-value for $\mathcal P$. Suppose that for every $y\in\R$ there
exist laws $P_h\in\mathcal P$ such that $P_h\{y\}>0$ and
\begin{equation}\label{adm:eq:approxwitness}
 \frac{1-\E_{P_h}e(Y)}{P_h\{y\}}\longrightarrow0.
\end{equation}
Then $e$ is pointwise admissible. In particular, it suffices that for each $y$
there is a law $P_y$ with $P_y\{y\}>0$ and $\E_{P_y}e=1$.
\end{lemma}
\begin{proof}
If $g\ge e$ and $g(y)\ge e(y)+d$ for some $d>0$, then
\[
 \E_{P_h}g\ge\E_{P_h}e+dP_h\{y\}>1
\]
for sufficiently small $h$, contradicting validity. A possibly infinite strict
gain contains a finite positive gain and is covered by the same argument.
\end{proof}

\begin{lemma}[Affine completion]\label{adm:lem:affine}
Suppose all laws in $\mathcal P$ have mean zero. Let $E$ be a finite-valued
admissible e-value,
$c>0$, and $\lambda\in\R$. If
\begin{equation}\label{adm:eq:affine}
 F(y)=\frac{1+\lambda y+cE(y)}{1+c}\ge0
 \qquad\text{for all }y,
\end{equation}
then $F$ is admissible.
\end{lemma}
\begin{proof}
Validity follows by taking expectations and using $\E Y=0$. If a valid $G$
dominates $F$, define
\[
 \widetilde E=E+\frac{1+c}{c}(G-F)
             =\frac{(1+c)G-1-\lambda y}{c}.
\]
The first expression shows $\widetilde E\ge E\ge0$. The second gives
$\E_P\widetilde E\le1$ under every null law. Thus a strict improvement of $F$
would strictly improve $E$, which is impossible.
\end{proof}

\begin{corollary}[Mixing with one]\label{adm:cor:mixone}
If $E$ is finite-valued and admissible, then $(1-t)+tE$ is admissible for every $0\le t\le1$,
provided the constant one is itself admissible when $t=0$.
\end{corollary}
\begin{proof}
For $t>0$, a dominator $G$ would give the dominator
$E+[G-((1-t)+tE)]/t$ of $E$. Its null expectation is at most one.
\end{proof}

\subsection{A self-contained sub-Gaussian admissibility certificate}

The usual observation that an e-value has Gaussian expectation one immediately
rules out continuous strict improvements. It does not by itself rule out an
arbitrary measurable improvement at a single point, since a Gaussian assigns
zero mass to that point. The following argument resolves this issue.

\begin{theorem}[Gaussian saturation certifies admissibility]\label{adm:thm:sgsat}
Let $e:\R\to[0,\infty)$ be continuous and valid for $\Ps$. If
\[
 \E e(Z)=1,\qquad Z\sim N(0,s^2),
\]
then $e$ is pointwise admissible, even against arbitrary measurable dominators.
\end{theorem}
\begin{proof}
Fix $y\ne0$ and set $r=|y|$. For $h>0$, let
\[
 A_h=\{r\le|Z|\le r+h\},\qquad
 Y_h=\begin{cases}\sign(Z)r,&Z\in A_h,\\Z,&Z\notin A_h.\end{cases}
\]
The variable $Y_h$ is symmetric, centered, and satisfies $|Y_h|\le|Z|$. Hence,
for every real $t$,
\[
 \E\ee^{tY_h}=\E\cosh(tY_h)
 \le\E\cosh(tZ)=\ee^{s^2t^2/2}.
\]
Its law belongs to $\Ps$, and $\Pp(Y_h=y)=\Pp(A_h)/2>0$.
Put
\[
 \omega_h=\max_{\varepsilon\in\{-1,1\}}
   \sup_{r\le u\le r+h}|e(\varepsilon u)-e(\varepsilon r)|.
\]
Continuity implies $\omega_h\to0$, and
\[
 \left|1-\E e(Y_h)\right|
 =\left|\E[e(Z)-e(Y_h)]\right|
 \le\omega_h\Pp(A_h).
\]
Thus the ratio in \eqref{adm:eq:approxwitness} is at most $2\omega_h$. For $y=0$,
map $Z$ to zero on $\{|Z|\le h\}$ instead; the same proof uses
$\omega_h=\sup_{|z|\le h}|e(z)-e(0)|$ and has no factor two. Apply
Lemma~\ref{adm:lem:witness}.
\end{proof}

\begin{corollary}[MGF factors and a quadratic example]\label{adm:cor:mgf}
For every $\theta\in\R$,
\begin{equation}\label{adm:eq:mgf}
 E_\theta(y)=\exp\{\theta y-s^2\theta^2/2\}
\end{equation}
is admissible for $\Ps$. So is, for every $u\in\R$,
\begin{equation}\label{adm:eq:quadratic}
 H_u(y)=\frac{(y+u)^2}{s^2+u^2}.
\end{equation}
\end{corollary}
\begin{proof}
Both are continuous and have expectation one under $N(0,s^2)$. MGF validity is
part of the null definition. For the quadratic, sub-Gaussianity implies
$\E Y^2\le s^2$: indeed
$\E\cosh(tY)\le\ee^{s^2t^2/2}$ and
$\cosh(tY)\ge1+t^2Y^2/2$; let $t\to0$. Centering now proves validity of
\eqref{adm:eq:quadratic}. Apply Theorem~\ref{adm:thm:sgsat}.
\end{proof}

\begin{remark}
Every point is an atom of some law in $\Ps$, as the proof shows. Thus the only
set negligible under every null law is the empty set. Pointwise admissibility
and admissibility modulo common null sets coincide here. The same will hold for
the moment class. No reference measure assumption is needed.
\end{remark}

\subsection{A class of admissible sub-Gaussian e-values}\label{adm:sec:insurance}

A linear bet $1+\lambda y$ is not nonnegative on the real line. Its negative
tail can be insured by an exponential without sacrificing admissibility.

\begin{theorem}[An admissible insured linear bet]\label{adm:thm:insurance}
For $\lambda\ne0$ and $d>0$, define
\begin{align}
 c_{\lambda,d}
 &=d\exp\left\{-1-\frac1d+\frac{s^2\lambda^2}{2d^2}\right\},\label{adm:eq:insurancefee}\\
 F_{\lambda,d}(y)
 &=\frac{1+\lambda y+d\exp\{-1-(1+\lambda y)/d\}}
          {1+c_{\lambda,d}}.\label{adm:eq:insurance}
\end{align}
Then $F_{\lambda,d}$ is a pointwise admissible e-value for $\Ps$.
\end{theorem}
\begin{proof}
For every real $z$,
\[
 z+d\ee^{-1-z/d}=d\{w-1+\ee^{-w}\}\ge0,
 \qquad w=1+z/d,
\]
by $\ee^{-w}\ge1-w$. Therefore the numerator is nonnegative, vanishing only
when $1+\lambda y=-d$. Also
\[
 d\ee^{-1-(1+\lambda y)/d}
 =c_{\lambda,d}E_{-\lambda/d}(y).
\]
Apply Corollary~\ref{adm:cor:mgf} and Lemma~\ref{adm:lem:affine}. Equivalently, the
factor is continuous, valid, and Gaussian-saturated, so
Theorem~\ref{adm:thm:sgsat} applies directly.
\end{proof}

\begin{proposition}[Optimizing the insurance charge]\label{adm:prop:optinsurance}
Writing $t=s|\lambda|$, the charge in \eqref{adm:eq:insurancefee} is uniquely
minimized over $d>0$ at
\begin{equation}\label{adm:eq:optd}
 d_*(t)=\frac{\sqrt{1+4t^2}-1}{2},\qquad
 c_*(t)=d_*(t)\exp\left\{-\frac12-\frac1{2d_*(t)}\right\}.
\end{equation}
Every choice of $d>0$ remains admissible; this optimization is a separate claim
about the fee.
\end{proposition}
\begin{proof}
The derivative of the log charge has the sign of
$d^2+d-t^2$. It changes sign exactly once on $(0,\infty)$. At its root,
$t^2=d^2+d$, giving the displayed charge.
\end{proof}

For fixed $y$, the optimized factor has the expansion
\begin{equation}\label{adm:eq:inslocal}
 \log F_{\lambda,d_*(s|\lambda|)}(y)
 =\lambda y-\frac12\lambda^2y^2+o(\lambda^2),\qquad \lambda\to0.
\end{equation}
Indeed, $d_*\sim s^2\lambda^2$, and both the fee and the insurance term at a
fixed $y$ vanish faster than every power of $|\lambda|$. Thus the expansion is
that of $\log(1+\lambda y)$. Define the factor at $\lambda=0$ to be one.
This is a pointwise expansion, not a uniform assertion over unbounded $y$.

\subsection{The complete admissible class under a moment bound}\label{adm:sec:moment}

The moment class admits a particularly concrete complete answer, valid even
when $1<p\le2$.

\begin{lemma}[Moment-boundary witnesses at every point]\label{adm:lem:momentwitness}
For every $y\in\R$ there is a law $P_y\in\Pm$ with $P_y\{y\}>0$ and
$\E_{P_y}|Y|^p=M$.
\end{lemma}
\begin{proof}
Choose sufficiently small $\varepsilon>0$ and define
\[
 r_\varepsilon^p=\frac{M-\varepsilon|y|^p}{1-\varepsilon},\qquad
 m_\varepsilon=-\frac{\varepsilon y}{1-\varepsilon}.
\]
As $\varepsilon\downarrow0$, $r_\varepsilon\to M^{1/p}>0$ and
$m_\varepsilon\to0$, so $r_\varepsilon>|m_\varepsilon|$. Place mass
$\varepsilon$ at $y$ and the remaining mass at $\pm r_\varepsilon$, with
conditional probabilities $(1\pm m_\varepsilon/r_\varepsilon)/2$. The mean
is zero and the absolute $p$th moment is exactly $M$.
\end{proof}

\begin{theorem}[Complete moment-class characterization]\label{adm:thm:momentcomplete}
For $p>1$ and $M>0$, the pointwise admissible e-values for $\Pm$ are exactly
the globally nonnegative functions
\begin{equation}\label{adm:eq:canonical}
 h_{a,b}(y)=1+ay+b(|y|^p-M),\qquad a\in\R,\quad b\ge0.
\end{equation}
Every e-value for $\Pm$ is pointwise dominated by one of these functions.
\end{theorem}
\begin{proof}
\emph{Validity and admissibility.}
If \eqref{adm:eq:canonical} is nonnegative, it is valid by the two constraints.
Every law in Lemma~\ref{adm:lem:momentwitness} gives $\E h_{a,b}=1$. The exact
witness case of Lemma~\ref{adm:lem:witness} proves admissibility.

\emph{A canonical dominator for every e-value.}
Let $e$ be valid. Lemma~\ref{adm:lem:momentwitness} implies $e(y)<\infty$ for all
$y$. Fix $r$ with $0<r<M^{1/p}$. For any finite set $S$ containing $0,\pm r$,
consider the finite linear program
\[
 \max\left\{\sum_{y\in S}w_y e(y):
 w_y\ge0,\ \sum w_y=1,\ \sum w_yy=0,\ \sum w_y|y|^p\le M\right\}.
\]
It is feasible and has value at most one. Its dual minimizes $\kappa$ subject
to
\[
 \kappa+ay+b(|y|^p-M)\ge e(y)\quad(y\in S),\qquad b\ge0.
\]
Finite linear-programming strong duality gives such coefficients with
$\kappa\le1$. Increasing $\kappa$ to one preserves domination. The
inequalities at $0,\pm r$, together with $e\ge0$, imply
\[
 0\le b\le1/M,\qquad
 |a|r\le1-b(M-r^p)\le1.
\]
Consequently, the coefficients can always be chosen in the same compact
rectangle $K=[-1/r,1/r]\times[0,1/M]$. For each real $y$, the condition
$1+ay+b(|y|^p-M)\ge e(y)$ defines a closed subset of $K$. Every finite
collection of these closed sets has a common point, by the finite program
above. Compactness gives a common point for all $y\in\R$. This gives a
globally nonnegative canonical dominator $h_{a,b}\ge e$. If $e$ is admissible,
it must equal its dominator everywhere.
\end{proof}

\begin{proposition}[An explicit parameter region]\label{adm:prop:region}
For $b>0$, global nonnegativity of \eqref{adm:eq:canonical} is equivalent to
\begin{equation}\label{adm:eq:region}
 1-bM\ge(p-1)p^{-p/(p-1)}
       |a|^{p/(p-1)}b^{-1/(p-1)}.
\end{equation}
For $b=0$, it is equivalent to $a=0$. When $p=2$, the condition is
$a^2\le4b(1-bM)$.
\end{proposition}
\begin{proof}
For $b>0$, minimize $ay+b|y|^p$ over $y$. Its minimum is the negative of
the right side of \eqref{adm:eq:region}. For $b=0$, a nonconstant affine function
cannot be nonnegative on all of $\R$.
\end{proof}

\subsection{Self-normalized exponentials are inadmissible}
\label{adm:sec:sn}

The common unnormalized factor is
\begin{equation}\label{adm:eq:selfnorm}
 f_\lambda(y)=\exp\{\lambda y-\lambda^2y^2/2\}.
\end{equation}
Its observed quadratic term is attractive for variance adaptation. It need not
be an e-value on its own. We first establish the elementary inequality used to
repair it, then construct admissible dominators of the repaired versions.

\begin{lemma}[A power remainder]\label{adm:lem:power}
For $2<q\le3$, set $D_q=3^{(1-q)/2}$. Then, for every real $u$,
\begin{equation}\label{adm:eq:power}
 \ee^{u-u^2/2}\le1+u+D_q(u_-)^q
                  \le1+u+D_q|u|^q.
\end{equation}
In particular, the cubic coefficient is $D_3=1/3$.
\end{lemma}
\begin{proof}
For $u\ge0$, $\log(1+u)\ge u-u^2/2$: the derivative of their difference is
$u^2/(1+u)\ge0$. For $u=-x<0$, let
$R(x)=\ee^{-x-x^2/2}-1+x$. We have $R(x)\le x$, since the exponential is at
most one. Also $R(x)\le x^3/3$. To see the latter, define
$g(x)=1-x+x^3/3-\ee^{-x-x^2/2}$. Then $g(0)=g'(0)=0$ and
\[
 g''(x)=x\{2-(2+x)\ee^{-x-x^2/2}\}\ge0.
\]
For the last inequality, $(2+x)\ee^{-x-x^2/2}\le2$, as follows by
differentiating its logarithm. Finally, splitting at $x=\sqrt3$ gives
\[
 \min\{x^3/3,x\}\le3^{(1-q)/2}x^q,
 \qquad 2<q\le3.
\]
For $x\le\sqrt3$, use $x^{3-q}/3\le3^{(1-q)/2}$; for larger $x$, use
$x^{1-q}\le3^{(1-q)/2}$.
\end{proof}

If a mean-zero law obeys $\E|Y|^q\le M_q$, this implies
$\E f_\lambda(Y)\le1+c_\lambda$ with
$c_\lambda=D_qM_q|\lambda|^q$. Two valid repairs are therefore
\begin{equation}\label{adm:eq:tworepairs}
 S_\lambda^{\exp}(y)=\ee^{\lambda y-\lambda^2y^2/2-c_\lambda},
 \qquad
 S_\lambda^{\log}(y)=\frac{f_\lambda(y)}{1+c_\lambda}.
\end{equation}
When $c_\lambda>0$, $1+c_\lambda<\ee^{c_\lambda}$, so
$S_\lambda^{\log}>S_\lambda^{\exp}$ everywhere. This already proves that
the exponential-fee version is inadmissible. The next results show that the
logarithmic-fee version is also inadmissible for the charges specified here.

\subsubsection{A completion for the actual \texorpdfstring{$p$}{p}-moment class}

\begin{theorem}[Admissible self-normalized dominator under $p>2$]
\label{adm:thm:momentsn}
Suppose $p>2$ and the known bound is $\E|Y|^p\le M$. Put
\[
 q=\min\{p,3\},\qquad r=M^{1/p},\qquad
 B_{p,q,r}(y)=\frac qp r^{q-p}|y|^p+
                 \left(1-\frac qp\right)r^q,
 \qquad c_\lambda=D_qr^q|\lambda|^q.
\]
Then
\begin{equation}\label{adm:eq:momentsn}
 G_\lambda^{\rm mom}(y)=
 \frac{1+\lambda y+D_q|\lambda|^q B_{p,q,r}(y)}{1+c_\lambda}
\end{equation}
is admissible for the full class $\Pm$. It dominates both factors in
\eqref{adm:eq:tworepairs} with this $c_\lambda$. For $\lambda\ne0$ it strictly
dominates the logarithmically repaired factor somewhere.
\end{theorem}
\begin{proof}
Concavity of $z\mapsto z^{q/p}$ and its tangent at $z=r^p$ give
\begin{equation}\label{adm:eq:tangent}
 |y|^q\le B_{p,q,r}(y).
\end{equation}
Consequently, \eqref{adm:eq:power} shows that the numerator in
\eqref{adm:eq:momentsn} dominates $f_\lambda(y)>0$. Its expectation is at most
$1+c_\lambda$, because $\E B_{p,q,r}(Y)\le r^q$. Hence it is valid.
More specifically, writing $D=1+c_\lambda$, it equals
\[
 1+\frac\lambda D y+
 \frac{D_q|\lambda|^q(q/p)r^{q-p}}{D}(|y|^p-M).
\]
Thus it is canonical for the \emph{actual $p$-moment class}, and
Theorem~\ref{adm:thm:momentcomplete} proves admissibility. Domination of both
repairs follows from the numerator majorization and $1+c_\lambda\le
\ee^{c_\lambda}$. The polynomial numerator grows in the tails, whereas
$f_\lambda(y)\to0$ as $|y|\to\infty$, proving a strict gain.
\end{proof}

For $2<p\le3$, the formula simplifies to
\[
 G_\lambda^{\rm mom}(y)=
 \frac{1+\lambda y+D_p|\lambda y|^p}{1+D_pM|\lambda|^p}.
\]
When $p>3$, merely applying a third-moment result would prove validity but
would not establish admissibility for the smaller $p$-moment class. The tangent
completion \eqref{adm:eq:tangent} is what supplies that stronger conclusion.

\subsubsection{A sub-Gaussian completion}\label{adm:sec:sgsn}

\begin{theorem}[Admissible self-normalized dominator under sub-Gaussianity]
\label{adm:thm:sgsn}
Let $\lambda\ne0$ and define
\begin{equation}\label{adm:eq:sgsnconstants}
 \gamma=\sqrt3\,\ee^{-3/2},\qquad
 \theta=\frac{\sqrt3}{s},\qquad
 \varepsilon=\sign(\lambda),\qquad
 c_\lambda=\gamma s^3|\lambda|^3.
\end{equation}
The factor
\begin{equation}\label{adm:eq:sgsn}
 G_\lambda^{\rm sg}(y)=
 \frac{1+\lambda y+c_\lambda\exp\{-\varepsilon\theta y-3/2\}}
      {1+c_\lambda}
\end{equation}
is admissible for $\Ps$ and dominates both
\begin{equation}\label{adm:eq:sgsnrepairs}
 \frac{f_\lambda(y)}{1+c_\lambda}
 \quad\text{and}\quad
 \exp\{\lambda y-\lambda^2y^2/2-c_\lambda\}.
\end{equation}
Both repaired factors are valid and inadmissible for $\lambda\ne0$.
\end{theorem}
\begin{proof}
For $x\ge0$ and $\theta>0$, maximizing $x^3\ee^{-\theta x}$ gives
\begin{equation}\label{adm:eq:cubicexp}
 \frac{x^3}{3}\le\frac9{\ee^3\theta^3}\ee^{\theta x}.
\end{equation}
Apply this to $x=-\varepsilon y$ when that quantity is nonnegative. For the
opposite sign the negative-part cube is zero. With \eqref{adm:eq:sgsnconstants},
\[
 \frac{((\lambda y)_-)^3}{3}
 \le\sqrt3\,s^3|\lambda|^3\ee^{-3}\ee^{-\varepsilon\theta y}
 =c_\lambda E_{-\varepsilon\theta}(y).
\]
The $q=3$ case of Lemma~\ref{adm:lem:power} now gives
\[
 f_\lambda(y)\le1+\lambda y+c_\lambda E_{-\varepsilon\theta}(y).
\]
The numerator of \eqref{adm:eq:sgsn} is therefore positive. Its factor is an
affine completion of the admissible MGF e-value
$E_{-\varepsilon\theta}$, so Lemma~\ref{adm:lem:affine} proves admissibility.
The same inequality proves validity and domination of the repaired factors.
As $y\to+\infty$ with $\lambda>0$, or $y\to-\infty$ with $\lambda<0$,
the numerator of the completion grows linearly while $f_\lambda$ tends to
zero. The logarithmically repaired factor is thus strictly dominated
somewhere; the exponential repair is strictly smaller still.
\end{proof}

Define $G_0^{\rm sg}=1$. No bound on a moment separate from sub-Gaussianity
was imposed: the MGF controls the insurance term directly.

\section{Using Gaffke's method}\label{app:gaffke}

The preceding constructions build supermartingales and are valid for adapted
data. A different route is available when the observations are independent:
convert the one-round e-values from Appendix~\ref{app:admissible} directly to
a fixed-time $p$-value using Gaffke's Dirichlet average. This appendix gives a
Gaussian-efficient extension to unbounded observations. Its endpoints require
numerical quantiles or root finding, which is why it is separated from the
analytically closed-form constructions in the main paper.

\subsection{The independent e-to-p merger}

For a nonnegative vector $\mathbf e=(e_1,\ldots,e_n)$, define
\begin{equation}\label{eq:gaffke-merger}
 K_n(\mathbf e)=\mathbb P_{\mathbf D}
 \left\{\sum_{i=1}^n e_iD_i\le1\right\},\qquad
 (D_0,D_1,\ldots,D_n)\sim\operatorname{Dirichlet}(1,\ldots,1).
\end{equation}
The probability is conditional on the observed vector. The coordinate $D_0$
has coefficient zero and must be included. Equivalently, for independent
unit exponentials $E_0,E_1,\ldots,E_n$,
\begin{equation}\label{eq:gaffke-exponential}
 K_n(\mathbf e)=\mathbb P\left\{
       \sum_{i=1}^n(e_i-1)E_i\le E_0\right\}.
\end{equation}
Thus $K_n$ is nonincreasing in every input. The theorem of
Vlassis and Thomas~\cite{VlassisThomas2026}, establishing a conjecture of
Gaffke~\cite{Gaffke2005}, states that
\begin{equation}\label{eq:gaffke-validity}
 \mathbb P\{K_n(\mathbf e)\le u\}\le u,\qquad 0\le u\le1,
\end{equation}
whenever the random inputs are independent, nonnegative, and satisfy
$\mathbb E e_i\le1$ individually. We use this independent merging theorem as
an external result; its proof is not an optional-stopping argument.
Conditional e-validity relative to the past, by itself, does not imply the
independence required in \eqref{eq:gaffke-validity}. A predictable betting
policy cannot be inserted into this merger without a separate dependence
theorem. Conditioning on an independent pilot is permissible when the scored
observations and their e-values remain conditionally independent.

\subsection{A monotone extension to unbounded means}

Choose a deterministic $R_n>0$ and a stop-loss charge $\delta_n\ge0$ satisfying
\[
 \mathbb E(\pm(X_i-\mu)-R_n)_+\le\delta_n
 \quad\text{for every }i.
\]
The following choices follow from the assumptions in the main paper:
\begin{align}
 \delta_n^{\mathrm{sg}}
 &=\int_{R_n}^\infty\exp\{-u^2/(2s^2)\}\,du,
       &&\text{under sub-Gaussianity},\label{eq:gaffke-sg-charge}\\
 \delta_n^{(p)}
 &=C_pM R_n^{1-p},\qquad
 C_p=\frac{(p-1)^{p-1}}{p^p},
       &&\text{if }\mathbb E|X_i-\mu|^p\le M.
       \label{eq:gaffke-moment-charge}
\end{align}
For completeness, the first follows by integrating the Chernoff tail bound.
The second follows from
$(x-R)_+\le C_pR^{1-p}|x|^p$, obtained by maximizing
$(x-R)/x^p$ over $x>R$.

For a candidate $m$, form
\begin{equation}\label{eq:gaffke-inputs}
 \begin{aligned}
 e_{i,n}^+(m)&=\frac{(R_n+X_i-m)_+}{R_n+\delta_n},\qquad
 e_{i,n}^-(m)=\frac{(R_n+m-X_i)_+}{R_n+\delta_n},\\
 p_n^\pm(m)&=K_n\bigl(e_{1,n}^\pm(m),\ldots,e_{n,n}^\pm(m)\bigr).
 \end{aligned}
\end{equation}
At $m=\mu$, each input has expectation at most one, since
$(R+y)_+=R+y+(-R-y)_+$. If the observations are independent, so are the
inputs. Both $p$-values are therefore valid at the true mean. Moreover,
$p_n^+(m)$ is nondecreasing in $m$ and $p_n^-(m)$ is nonincreasing in $m$.
For $q=\alpha/2$, let
\begin{equation}\label{eq:gaffke-endpoints}
 L_n=\inf\{m:p_n^+(m)\ge q\},\qquad
 U_n=\sup\{m:p_n^-(m)\ge q\}.
\end{equation}
These endpoints are finite. The closure of the joint nonrejection set is an
interval whenever it is nonempty. To ensure an interval-valued output for
every sample, including extreme samples for which the two bounds cross,
define
\begin{equation}\label{eq:gaffke-interval}
 I_n^{\mathrm G}=[\min\{L_n,U_n\},\max\{L_n,U_n\}].
\end{equation}
This convention enlarges the possibly empty joint nonrejection set and does
not affect the eventual endpoints in the efficiency theorem below.

\begin{proposition}[Fixed-time coverage]\label{prop:gaffke-coverage}
For independent observations with common mean $\mu$ and either of the stated
envelopes, $\mathbb P\{\mu\in I_n^{\mathrm G}\}\ge1-\alpha$.
\end{proposition}
\begin{proof}
If both $p_n^+(\mu)>q$ and $p_n^-(\mu)>q$, monotonicity gives
$L_n\le\mu\le U_n$. Therefore failure of coverage implies
$p_n^+(\mu)\le q$ or $p_n^-(\mu)\le q$. Apply
\eqref{eq:gaffke-validity} and a union bound.
\end{proof}

\subsection{Gaussian efficiency}

The mechanism is transparent in the exponential representation. Without
clipping near the relevant candidates, Gaffke's construction behaves like a
Dirichlet-weighted empirical mean with one extra lower or upper observation
of size $R_n$. This synthetic observation costs order $R_n/n$, which is
negligible if $R_n=o(\sqrt n)$. The stop-loss fee costs $\delta_n$ and must
also be negligible at the Gaussian scale.

\begin{theorem}[Unbounded Gaffke intervals attain the Gaussian endpoints]
\label{thm:gaffke-unbounded}
Suppose $X_i$ are iid with mean $\mu$ and variance $0<\sigma^2<\infty$.
Assume the deterministic tuning satisfies
\begin{equation}\label{eq:gaffke-conditions}
 \frac{R_n}{\sqrt n}\to0,\qquad
 \sqrt n\,\delta_n\to0,\qquad
 \limsup_{n\to\infty}\frac{\max_{i\le n}|X_i-\mu|}{R_n}<1
 \quad\text{almost surely},
\end{equation}
and $R_n\to\infty$. For every fixed $\alpha\in(0,1)$,
\[
 L_n=\bar X_n-\frac{\sigma z_{1-\alpha/2}}{\sqrt n}+o(n^{-1/2}),
 \qquad
 U_n=\bar X_n+\frac{\sigma z_{1-\alpha/2}}{\sqrt n}+o(n^{-1/2})
 \quad\text{almost surely}.
\]
The conditions hold with the charges above and the explicit choices
\begin{align}
 R_n&=2s\sqrt{\log(n+1)},
     &&\text{under sub-Gaussianity},\label{eq:gaffke-sg-R}\\
 R_n&=M^{1/p}(n+1)^{(p+2)/(4p)},
     &&\text{under a }p\text{th moment bound, }p>2.
       \label{eq:gaffke-moment-R}
\end{align}
\end{theorem}
\begin{proof}
Work on an almost-sure sample path where the sample mean and variance
converge and \eqref{eq:gaffke-conditions} holds. On every band
$|m-\bar X_n|\le C/\sqrt n$, no input in \eqref{eq:gaffke-inputs} is clipped
once $n$ is sufficiently large. Write $G_n=\sum_{i=1}^nE_i$ and
$W_n=\sum_{i=1}^nX_iE_i/G_n$. On that band, elementary rearrangement of
\eqref{eq:gaffke-exponential} gives
\begin{equation}\label{eq:gaffke-auxiliary}
 \begin{aligned}
 p_n^+(m)&=\mathbb P(A_n\le m),\qquad
 p_n^-(m)=\mathbb P(B_n\ge m),\\
 A_n&=W_n-\delta_n-(R_n+\delta_n)\frac{E_0}{G_n},\\
 B_n&=W_n+\delta_n+(R_n+\delta_n)\frac{E_0}{G_n}.
 \end{aligned}
\end{equation}
These probabilities are over exponentials only, with the data path fixed.

Let $Q_n=\sum_{i=1}^n(X_i-\bar X_n)^2$. The conditional weighted-mean CLT is
\begin{equation}\label{eq:gaffke-cond-clt}
 \frac{\sqrt n(W_n-\bar X_n)}{\sqrt{Q_n/n}}
 \ \Rightarrow\ N(0,1).
\end{equation}
Here is a direct verification. Since the centered data sum to zero,
\[
 W_n-\bar X_n
 =\frac{\sum_i(X_i-\bar X_n)(E_i-1)}{G_n}.
\]
Finite second moment implies
$\max_{i\le n}|X_i-\mu|/\sqrt n\to0$ almost surely: the individual
statement follows from Borel--Cantelli and
$\sum_i\mathbb P(|X_i-\mu|>\varepsilon\sqrt i)<\infty$, and then taking
maxima gives the claim. Consequently, the normalized coefficients
$b_{i,n}=(X_i-\bar X_n)/\sqrt{Q_n}$ satisfy
$\sum_i b_{i,n}^2=1$ and $\max_i|b_{i,n}|\to0$. Because unit exponentials
have finite third absolute centered moment,
$\sum_i|b_{i,n}|^3\le\max_i|b_{i,n}|\to0$, proving the numerator CLT by
the Lyapunov condition. Also $G_n/n\to1$ in probability. This proves
\eqref{eq:gaffke-cond-clt}.

The additive terms in \eqref{eq:gaffke-auxiliary}, multiplied by $\sqrt n$,
vanish in probability: $\sqrt n\delta_n\to0$ and
$\sqrt n(R_n+\delta_n)E_0/G_n=o_{\mathbb P}(1)$. Thus both auxiliary
variables have conditional Gaussian limits with center $\bar X_n$ and scale
$\sigma/\sqrt n$. Quantile convergence gives their fixed $q$ and $1-q$
quantiles. Choose $C>\sigma z_{1-q}$ and then a fixed positive margin.
The corresponding crossings lie in the local band for all sufficiently
large $n$. Monotonicity and \eqref{eq:gaffke-auxiliary} identify those
crossings with \eqref{eq:gaffke-endpoints}, proving the endpoint expansions.

For \eqref{eq:gaffke-sg-R}, Mills' elementary upper bound gives
$\delta_n\le(s^2/R_n)\exp\{-R_n^2/(2s^2)\}$, so
$\sqrt n\delta_n\to0$. The Chernoff bound and Borel--Cantelli give
$|X_i-\mu|\le s\sqrt{3\log i}$ eventually, whence the last limsup in
\eqref{eq:gaffke-conditions} is at most $\sqrt3/2<1$.
For the moment choice, put $a=(p+2)/(4p)$. Then
$1/p<a<1/2$ and $a>1/[2(p-1)]$. A finite $p$th moment gives
$\max_{i\le n}|X_i-\mu|=o(n^{1/p})$ almost surely by the same
Borel--Cantelli argument. Hence $\max_i|X_i-\mu|/R_n\to0$,
$R_n/\sqrt n\to0$, and
$\sqrt n\delta_n=O(n^{1/2-a(p-1)})\to0$.
\end{proof}

The same merger accepts the admissible completions from
Appendix~\ref{app:admissible}. Pointwise domination of each input decreases
$K_n$ at every fixed candidate and can therefore improve the accepted set.
However, these completed inputs generally lose monotonicity in $m$, so the
improved inversion need not remain an interval. One-round admissibility also
does not imply admissibility of Gaffke's merger; the two questions are distinct.
The positive-part version above deliberately preserves monotone inversion.

\subsection{Empirical results}
\label{sec:gaffke-experiments}

The fresh experiment in Table~\ref{tab:gaffke-unbounded} evaluates the
unbounded construction \eqref{eq:gaffke-inputs}--\eqref{eq:gaffke-sg-R} on
iid $N(0,\sigma^2)$ observations, with known proxy $s=1$. For each
$(\sigma,n)$, 120 independent data sets are reused across the three confidence
levels. Each entry is the mean width divided by the Chernoff width at that
same confidence level. The seed is 2026102807.

\begin{table}[t]
\centering
\caption{Mean-width ratio of the unbounded Gaffke interval to Chernoff.
Values below one favor Gaffke. These are numerical appendix methods, not
the closed-form intervals evaluated in the main paper.}
\label{tab:gaffke-unbounded}
\begin{tabular}{rrrrr}
\toprule
$\sigma$ & $n$ & $90\%$ & $95\%$ & $99\%$\\
\midrule
$0.5$ & $100$  & $0.649$ & $0.698$ & $0.804$\\
$0.5$ & $1000$ & $0.425$ & $0.447$ & $0.481$\\
$1.0$ & $100$  & $0.912$ & $0.962$ & $1.049$\\
$1.0$ & $1000$ & $0.752$ & $0.796$ & $0.859$\\
\bottomrule
\end{tabular}
\end{table}

We evaluated the uniform Dirichlet CDF deterministically by integrating its
normalized B-spline density and used 27 bisection steps per endpoint. For sorted
knots $a_0\le\cdots\le a_n$ with $a_n>a_0$, the density is
$nB(t)/(a_n-a_0)$, where $B$ is the degree-$(n-1)$ B-spline basis element on
those knots; repeated knots are handled directly. The knot vector here is
$(0,e_{1,n}^{\pm}(m),\ldots,e_{n,n}^{\pm}(m))$. If all knots are zero,
set $K_n=1$ directly. Thus no Monte Carlo
Dirichlet approximation is used. The outer experiment has only 120
replications per cell, so its coverage estimates are descriptive and do not
constitute high-precision validation at level $99\%$.
The Monte Carlo standard errors of these width ratios are below $0.004$.
At $\sigma=0.5$, the width reduction is already $20$--$35\%$ at $n=100$
and exceeds $50\%$ at $n=1000$. With the tight proxy $\sigma=s=1$, the
method is slightly wider than Chernoff at $n=100$ and confidence $99\%$,
but shorter in the other displayed cells. At $n=1000$, its width divided by
the Gaussian oracle width is $1.118$, $1.103$, and $1.085$ at the three
confidence levels. Observed coverages range from $0.908$ to $1.000$ at
$90\%$, from $0.967$ to $1.000$ at $95\%$, and from $0.992$ to $1.000$ at
$99\%$. The relatively high coverages in several small-sample cells reflect
the finite-sample price of the added zero coordinate and the chosen
stop-loss envelope. Typical run times in this experiment were $6$--$7$ ms
per interval at $n=100$ and $130$--$148$ ms at $n=1000$; these are ordinary
wall-clock measurements, not optimized performance benchmarks.

\section{Scalar confidence intervals: constructions and proofs}
\label{app:scalar-calibration}

This appendix supplies the details for Section~\ref{sec:subgaussian-ci}.
It begins with the local-MGF sample-scatter process of
Theorem~\ref{thm:subpsi-Q}, then gives optional tuning for the exponential
interval and the proof of Gaussian calibration. These are two different
tasks: the exponential process learns the leading variance, while the
Gaussian continuation supplies the exact terminal quantile. We prove the
latter under a moment bound so that Section~\ref{sec:moments} can reuse it.

\subsection{The sample-scatter process under a local MGF bound}
\label{app:subpsi-Q-proof}

Theorem~\ref{thm:subpsi-Q} uses the MGF only near zero. The main issue is
to replace a quadratic centered at the unknown mean by the observed sample
scatter. We first make the one-observation bound uniform over a forecast
chosen from the past, then use Welford's identity to assemble the scatter.
This explains both the cubic repair and the smaller harmonic MGF charge.

\begin{proof}[Proof of Theorem~\ref{thm:subpsi-Q}]
The proof follows four steps, all under the conditional-mean model.
\begin{itemize}
 \item \textbf{Step 1: extract moment bounds.} The two endpoint MGF values
 give the local quadratic and third-moment bounds needed below.
 \item \textbf{Step 2: allow any predictable forecast.} A scalar cubic
 inequality controls the exponential factor uniformly over its quadratic
 center.
 \item \textbf{Step 3: assemble the sample scatter.} Conditional H\"older
 interpolation and Welford's identity give the stated process.
 \item \textbf{Step 4: obtain the variants.} The same factor argument
 permits a pure cubic repair and the sharper MGF envelope.
\end{itemize}

\noindent\textbf{Step 1: extract moment bounds.}
We need bounds that depend only on the supplied MGF envelope, rather than
on unknown features of the conditional distribution.
Centering and Jensen's inequality imply $K(u)\ge0$ on the parameter
interval, so $A_R\ge0$. For $|u|\le R$, Taylor expansion gives the
pointwise bound
\begin{align*}
 \exp(uy)-1-uy
 &\le\sum_{k=2}^\infty\frac{|uy|^k}{k!}\\
 &\le\frac{u^2}{R^2}\{\exp(R|y|)-1-R|y|\}\\
 &\le\frac{u^2}{R^2}\{\exp(Ry)+\exp(-Ry)-2\}.
\end{align*}
Consequently,
\begin{equation}\label{eq:subpsi-Q-localMGF}
 \mathbb E_{i-1}\exp(u\xi_i)
 \le1+\frac{A_Ru^2}{R^2}
 \le\exp(S^2u^2/2),\qquad |u|\le R.
\end{equation}
Also, for every real $z$,
\[
 |z|^3\le 2\sqrt3\left(\frac{z^2}{2}+\frac{z^4}{24}\right)
           \le2\sqrt3\{\cosh z-1\}.
\]
The first inequality is equivalent, away from zero, to
$(|z|-\sqrt{12})^2\ge0$. Taking conditional expectations with
$z=R\xi_i$ yields
\begin{equation}\label{eq:subpsi-Q-moments}
 \mathbb E_{i-1}|\xi_i|^3\le M_3,\qquad
 \mathbb E_{i-1}\xi_i^2\le S^2.
\end{equation}
For the variance bound, one may either take the limit in
\eqref{eq:subpsi-Q-localMGF} or use
$R^2\xi_i^2\le\exp(R\xi_i)+\exp(-R\xi_i)-2$.
If $A_R=0$, this last inequality implies $\xi_i=0$ almost surely;
\eqref{eq:subpsi-Q-localMGF}--\eqref{eq:subpsi-Q-moments} then hold
with the stated arbitrary positive $S$ and $M_3=0$.

\noindent\textbf{Step 2: allow any predictable forecast.}
The following scalar inequality removes the quadratic-order error and
leaves only a one-sided cubic remainder. For every real $u$, with
$u_- =\max\{-u,0\}$,
\begin{equation}\label{eq:subpsi-Q-scalar}
 \exp(u-u^2/2)\le1+u+\frac{(u_-)^3}{3}.
\end{equation}
For $u\ge0$, this follows from
$\log(1+u)\ge u-u^2/2$, since the derivative of the difference is
$u^2/(1+u)\ge0$. For $u=-x$, $x\ge0$, define
$g(x)=1-x+x^3/3-\exp(-x-x^2/2)$. Then $g(0)=g'(0)=0$ and
\[
 g''(x)=x\{2-(2+x)\exp(-x-x^2/2)\}\ge0.
\]
Indeed, $(2+x)\exp(-x-x^2/2)$ starts at two and is decreasing,
as is seen by differentiating its logarithm. This proves
\eqref{eq:subpsi-Q-scalar}.

We now allow the quadratic to be centered at a forecast. This is the
part that will make the final variance clock independent of the candidate
mean. We claim that, for every finite $\mathcal F_{i-1}$-measurable real $a$,
\begin{equation}\label{eq:subpsi-Q-forecast}
 \mathbb E_{i-1}\exp\left\{
       \lambda\xi_i-\frac{(\lambda\xi_i+a)^2}{2}\right\}
 \le 1+\frac{4M_3}{9}|\lambda|^3+\frac{40S^4}{23}\lambda^4
 =\exp\{\rho(\lambda)\}.
\end{equation}
The statement is immediate when $\lambda=0$. Otherwise, suppress
conditioning, put $U=\lambda\xi_i$ and $t=S|\lambda|\le1/5$, and
apply \eqref{eq:subpsi-Q-moments}. We have
\[
 \mathbb E U=0,\quad \mathbb E U^2\le t^2,\quad
 \mathbb E U_-\le t/2,\quad
 \mathbb E U_-^3\le M_3|\lambda|^3.
\]
Here $\mathbb E U_-=\mathbb E|U|/2$ follows from centering.
The following inequalities hold for each value of the predictable shift
and therefore also conditionally.

There are three cases. A nonnegative shift only reduces the negative
tail in the cubic remainder; a sufficiently large negative shift is paid
for directly by the MGF; the remaining small negative shifts require a
short polynomial bound.

If $a\ge0$, \eqref{eq:subpsi-Q-scalar} gives
\begin{align*}
 \mathbb E\exp\{U-(U+a)^2/2\}
 &\le\exp(-a)\left\{1+a+\frac{\mathbb E(U+a)_-^3}{3}\right\}\\
 &\le1+\frac{M_3}{3}|\lambda|^3.
\end{align*}
Now let $a=-b<0$ and define
$b_0=t/(1-t)$ and $q=1/(1-t)=1+b_0$.
When $b\ge b_0$, completing the square gives, pointwise in $U$,
\[
 U-\frac{(U-b)^2}{2}
 \le qU-bb_0+\frac{b_0^2}{2}.
\]
The parameter restriction \eqref{eq:subpsi-Q-domain} ensures
\[
 |q\lambda|=\frac{|\lambda|}{1-S|\lambda|}\le R.
\]
Thus the local MGF bound \eqref{eq:subpsi-Q-localMGF}, with
$t^2q^2=b_0^2$, yields
\[
 \mathbb E\exp\{U-(U-b)^2/2\}
 \le\exp\{b_0(b_0-b)\}\le1.
\]
Only an MGF parameter inside $[-R,R]$ was used, even when the shift
$b$ is arbitrarily large.

It remains to consider $0<b<b_0\le1/4$. Using
$(U-b)_-\le U_-+b$ in \eqref{eq:subpsi-Q-scalar}, we obtain
\begin{align*}
 \mathbb E\exp\{U-(U-b)^2/2\}
 \le\exp(b)\left\{1-b+\frac{M_3}{3}|\lambda|^3
                         +bt^2+\frac{b^2t}{2}+\frac{b^3}{3}\right\}.
\end{align*}
For $b\le1/4$, $\exp(b)\le4/3$ and
$\exp(b)(1-b)\le1-b^2/2$. The latter follows by differentiating
both sides and comparing at zero. Since $t\le1/5$,
\begin{align*}
 \mathbb E\exp\{U-(U-b)^2/2\}
 &\le1+\frac{4M_3}{9}|\lambda|^3+\frac43bt^2
       +b^2\left\{-\frac12+\frac43\left(\frac t2+\frac b3\right)\right\}\\
 &\le1+\frac{4M_3}{9}|\lambda|^3+\frac43bt^2-\frac{23}{90}b^2\\
 &\le1+\frac{4M_3}{9}|\lambda|^3+\frac{40}{23}t^4.
\end{align*}
The last step maximizes the displayed quadratic in $b$ over all real
$b$. This proves \eqref{eq:subpsi-Q-forecast} in every case.

\noindent\textbf{Step 3: assemble the sample scatter.}
Welford's identity weights the squared forecast residual at time $i$ by
$(i-1)/i$. We reproduce that weight by interpolating between the repaired
factor and the ordinary MGF factor. At the true candidate $m=\mu$, let, for $i\ge2$,
\[
 r_i=\frac{i-1}{i},\qquad c_i=\bar X_{i-1},\qquad
 A_i=\exp\{\lambda\xi_i-\lambda^2(X_i-c_i)^2/2\},\qquad
 B_i=\exp(\lambda\xi_i).
\]
Applying \eqref{eq:subpsi-Q-forecast} with
$a=\lambda(\mu-c_i)$ bounds $\mathbb E_{i-1}A_i$ by
$\exp\{\rho(\lambda)\}$, uniformly over the predictable forecast.
Conditional H\"older's inequality therefore gives
\begin{align}
 \mathbb E_{i-1}\exp\{\lambda\xi_i-r_i\lambda^2(X_i-c_i)^2/2\}
 &=\mathbb E_{i-1}[A_i^{r_i}B_i^{1-r_i}]\notag\\
 &\le (\mathbb E_{i-1}A_i)^{r_i}
       (\mathbb E_{i-1}B_i)^{1-r_i}\notag\\
 &\le\exp\{r_i\rho(\lambda)+(1-r_i)K(\lambda)\}.
 \label{eq:subpsi-Q-holder}
\end{align}
For $i=1$, put $r_1=0$ and use only the natural MGF factor, bounded
by $\exp\{K(\lambda)\}$.
The update formula for the sample mean implies Welford's identity
\[
 Q_i-Q_{i-1}=\frac{i-1}{i}(X_i-\bar X_{i-1})^2\quad(i\ge2).
\]
Consequently,
\[
 Q_n=\sum_{i=2}^nr_i(X_i-\bar X_{i-1})^2,\qquad
 \sum_{i=1}^nr_i=n-H_n,\qquad
 \sum_{i=1}^n(1-r_i)=H_n.
\]
Multiplying the normalized factors in \eqref{eq:subpsi-Q-holder},
including the first-round MGF factor, gives exactly
\eqref{eq:subpsi-Q-process}. Each factor has conditional expectation
at most one, proving the nonnegative supermartingale assertion.

\noindent\textbf{Step 4: obtain the variants.}
Only the one-observation upper bounds change; the multiplication argument
from Step 3 stays the same.
Since $\log(1+x)\le x$ for $x\ge0$ and
$|\lambda|\le1/(5S)$,
\[
 \rho(\lambda)\le\frac{4M_3}{9}|\lambda|^3
                         +\frac{40S^4}{23}\lambda^4
 \le C|\lambda|^3.
\]
The forecast bound thus also holds with the larger repair
$C|\lambda|^3$, and the same conditional-factor proof applies.
Finally, both the assumed bound and \eqref{eq:subpsi-Q-localMGF}
control the natural MGF factor, so $K_R$ in
\eqref{eq:subpsi-Q-localenvelope} is also a valid envelope in
\eqref{eq:subpsi-Q-holder} and in the first round.
\end{proof}

\subsection{Optional finite-sample tuning}\label{app:scalar-tuning}

The simple bet in \eqref{eq:ci-simplebet} balances the evidence threshold
against the pilot estimate of variance. The refinement below also accounts
for the cubic cost when it matters at moderate sample sizes. Every quantity
is determined by the pilot and planned horizon, so the conditional coverage
argument is unchanged.

First cap the bet by the theorem's allowed range and by the balances with
the quadratic and cubic charges:
\begin{equation}\label{eq:ci-pilottuning}
 \lambda_0=\min\left\{\frac1{5s},
 \sqrt{\frac{2\log(2/\delta)}{N(\widehat v+s^2/k+s^2H/N)}},
 \left(\frac{\log(2/\delta)}{2s^3(N-H)}\right)^{1/3}\right\}.
\end{equation}
Next take one Newton step toward minimizing the pilot estimate of the
radius:
\begin{equation}\label{eq:ci-newton}
 \lambda=\min\left\{\frac1{5s},\quad
 \lambda_0\frac{6\log(2/\delta)
       -N(\widehat v+s^2/k+s^2H/N)\lambda_0^2}
       {4\log(2/\delta)+4s^3(N-H)\lambda_0^3}\right\}.
\end{equation}
The numerator is positive by the square-root constraint defining
$\lambda_0$. This is a fixed arithmetic refinement, with no search or
iteration. Omitting it changes neither coverage nor the limiting conclusions.

\subsection{Gaussian calibration under a moment bound}
\label{app:scalar-heat}

The Gaussian continuation is a smooth version of a terminal rejection
indicator. Updating its variance by observed squares cancels the quadratic
drift; a decreasing reserve pays for what remains. Three devices make that
argument valid for unbounded observations: a variance floor keeps the
potential smooth, freezing handles large increments, and a coarse interval
allows explicit inversion. The pilot only chooses the budget and does not
need a separate accuracy guarantee.

Use the pilot and scoring statistics from Section~\ref{subsec:ci-practical}.
We state the assumptions more generally to cover the later moment extension:
conditionally on the past,
\begin{equation}\label{eq:ci-heatassumptions}
 \E[X_i-\mu\mid\F_{i-1}]=0,\qquad
 \E[|X_i-\mu|^q\mid\F_{i-1}]\le M,\qquad
 \E[(X_i-\mu)^2\mid\F_{i-1}]\le V,
 \qquad 2<q\le3,
\end{equation}
where $M,V>0$ are known and the bounds hold almost surely. A zero
moment bound is the degenerate case already handled in
Section~\ref{sec:moments}.

For error level $\varepsilon\in(0,1)$, the following quantities have the
same roles as in the main text. The slack $\rho$ divides the error budget,
$R$ is the coarse radius, $\tau$ is the smoothing floor, and $A$ is the
initial variance budget:
\begin{equation}\label{eq:ci-heat-compact-general}
 \begin{aligned}
 \rho&=\frac1{8\log(n+e)},&
 R&=\sqrt{\frac{V}{2\varepsilon N\rho}},\\
 \tau&=\left\{\frac{6MN}{\varepsilon\rho(1-2\rho)}\right\}^{2/q},&
 A&=N(\widehat v+8V\rho)+NR^2+2\tau.
 \end{aligned}
\end{equation}
Use the radius \eqref{eq:ci-heat-simple-radius} when
$Q/N\le\widehat v+8V\rho$, and set $r_H(\varepsilon)=R$ otherwise.
The pilot variance has no coverage requirement. The observable check ensures
$A-Q\ge NR^2+2\tau>0$ whenever its square root is used.

Set the Gaussian terminal cutoff to
\[
 b=\sqrt A\,z_{1-\varepsilon(1-2\rho)^2/2},
\]
and use the continuation potential \eqref{eq:ci-heat-main-potential}.
\begin{theorem}[Closed-form Gaussian calibration]\label{thm:ci-heat}
Under \eqref{eq:ci-heatassumptions}, for every pilot-measurable
$\varepsilon\in(0,1)$, the interval
$I_H(\varepsilon)=[z-r_H(\varepsilon),z+r_H(\varepsilon)]$ satisfies
\begin{equation}\label{eq:ci-heatcoverage}
 \Pp\{\mu\in I_H(\varepsilon)\mid\F_k\}\ge1-\varepsilon.
\end{equation}
If the observations are iid with variance $\sigma^2>0$ and
$\varepsilon=\varepsilon_n\to\varepsilon_\infty\in(0,1)$, then almost surely
\begin{equation}\label{eq:ci-heatefficient}
 \sqrt n\,r_H(\varepsilon_n)\longrightarrow
          \sigma z_{1-\varepsilon_\infty/2},\qquad
 \sqrt n(z-\bar X_n)\longrightarrow0.
\end{equation}
\end{theorem}

We prove validity by first bounding one guarded update, then converting
the terminal rejection bound into interval coverage. Efficiency is a
separate pathwise argument: the pilot learns the variance and the reserve
becomes negligible. The first ingredient is the guarded drift bound below.
Small increments use the heat cancellation; large increments use only the
boundedness of the potential.

\begin{lemma}[Uniform guarded heat drift]\label{lem:ci-heatdrift}
Fix $\tau>0$. At an active state $(v,x)$ with $v\ge2\tau$, propose an increment
$u$. Accept the update $(v',x')=(v-u^2,x+u)$ if $v-u^2\ge\tau$, and otherwise
leave the state unchanged and freeze it permanently. Also freeze after an
accepted update leaving $v'<2\tau$. Then, for $2<q\le3$,
\begin{equation}\label{eq:ci-driftbound}
 |F(v',x')-F(v,x)-F_x(v,x)u|
 \le3\tau^{-q/2}|u|^q.
\end{equation}
The pointwise bound is asserted at active states only. At an already
frozen state the actual increment and its conditional drift are zero.
At the true mean, its conditional expected drift is at most
$3M\tau^{-q/2}$, including the rejected-update and frozen cases.
\end{lemma}
\begin{proof}
The proof separates the smooth, small-increment calculation from the tail
case, where a Taylor approximation need not be accurate.
\begin{itemize}
 \item \textbf{Step 1: control the derivatives.} The variance floor bounds
 the Gaussian derivatives uniformly in the current sum.
 \item \textbf{Step 2: bound the update.} Heat cancellation controls small
 increments, and boundedness controls large accepted or rejected increments.
 \item \textbf{Step 3: take conditional expectations.} Centering removes
 the linear term and the moment bound pays the remainder.
\end{itemize}

\noindent\textbf{Step 1: control the derivatives.}
Writing $y=(x-b)/\sqrt v$, direct differentiation gives
\[
 F_x=v^{-1/2}\phi(y),\quad F_{xxx}=v^{-3/2}(y^2-1)\phi(y),
 \quad F_{xv}=\tfrac12F_{xxx},\quad
 F_{vv}=\frac{(3y-y^3)\phi(y)}{4v^2}.
\]
In particular,
$|F_x|\le v^{-1/2}$, $|F_{xxx}|\le v^{-3/2}$,
$|F_{xv}|\le\frac12v^{-3/2}$, and $|F_{vv}|\le v^{-2}$.
These follow by maximizing the polynomials times the Gaussian density;
for example $(1+y^2)\phi(y)\le(1+2/\mathrm e)/\sqrt{2\pi}<1$.

\noindent\textbf{Step 2: bound the update.}
Write $\Delta F=F(v',x')-F(v,x)$. If $|u|\le\sqrt\tau$, the update is accepted, and all intermediate
variances are at least $\tau$. Expand first in space through second order,
then in variance through first order. The heat identity cancels the quadratic
terms, while the space remainder, the mixed derivative, and the variance
remainder give
\[
 |\Delta F-F_xu|
 \le\frac23\tau^{-3/2}|u|^3+\frac12\tau^{-2}u^4
 \le\frac76\tau^{-q/2}|u|^q.
\]
The last step uses $|u|/\sqrt\tau\le1$ and $q\le3$. If
$|u|>\sqrt\tau$, whether accepted or rejected,
$|\Delta F|\le1$ and hence
\[
 |\Delta F-F_xu|\le1+|u|/\sqrt\tau
 \le2\tau^{-q/2}|u|^q.
\]
\noindent\textbf{Step 3: take conditional expectations.}
At an active state, $F_x(v,x)$ is measurable from the past. Conditional
centering therefore removes $F_x(v,x)u$, and the moment bound gives
conditional upward drift at most $3M\tau^{-q/2}$. This includes rejected
proposals, since their contribution was bounded in Step 2. At a frozen
state the actual increment is zero, so its drift is zero directly. Thus
freezing is part of the drift calculation rather than an appeal to
optional stopping for an unguarded process.
\end{proof}

\begin{proof}[Proof of Theorem~\ref{thm:ci-heat-supermartingale}]
The remaining reserve decreases by a fixed amount at every observation.
It suffices to check that this decrease pays for the largest possible
conditional upward drift of the potential.
Condition on $\F_k$, so that $A$, $\tau$, and $b$ are fixed. At $m=\mu$,
either signed increment is conditionally centered, and
Lemma~\ref{lemma:cubic_inequality} gives
\[
 \E[|Z_i-\mu|^3\mid\F_{k+i-1}]\le6\gamma s^3.
\]
Apply Lemma~\ref{lem:ci-heatdrift} with $q=3$ and $M=6\gamma s^3$.
At active states the potential's conditional upward drift is at most
$18\gamma s^3/\tau^{3/2}$; at frozen states its drift is zero. Therefore
\begin{align*}
 \E[W_i(\mu)\mid\F_{k+i-1}]
 &\le F(v_{i-1},x_{i-1})+\frac{18\gamma s^3}{\tau^{3/2}}
       +\frac{18\gamma s^3(N-i)}{\tau^{3/2}}\\
 &=W_{i-1}(\mu).
\end{align*}
The guards keep $v_i\ge\tau$, so the Gaussian potential is well defined
and lies in $[0,1]$. The reserve is finite and nonnegative for
$0\le i\le N$, proving conditional integrability and nonnegativity.
Finally, $W_0>0$ is fixed given the pilot, so normalization preserves the
supermartingale property and gives initial value one.
\end{proof}

\begin{proof}[Proof of coverage in Theorem~\ref{thm:ci-heat}]
The terminal potential gives a test, but its states depend on the candidate
mean. The coarse interval lets us identify all relevant terminal states
and turn that test into the reported explicit radius.
\begin{itemize}
 \item \textbf{Step 1: calibrate the terminal test.} Add the drift reserve
 and apply conditional Markov's inequality to each sign.
 \item \textbf{Step 2: prevent freezing for localized candidates.} The
 variance check ensures that their squared residuals fit the initial budget.
 \item \textbf{Step 3: invert and account for error.} An excluded localized
 mean triggers one signed test; localization pays the remaining error.
\end{itemize}

\noindent\textbf{Step 1: calibrate the terminal test.}
Condition on $\F_k$. For a candidate $m$ and each sign, start the guarded
process of Lemma~\ref{lem:ci-heatdrift} at $(v_0,x_0)=(A,0)$ and propose
$u_i=\pm(Z_i-m)$. Keep the state fixed after freezing. At $m=\mu$,
\[
 W_i=F(v_i,x_i)+3M(N-i)\tau^{-q/2}
\]
is a nonnegative supermartingale: its remaining reserve decreases by
$3M\tau^{-q/2}$ per step, paying the conditional drift bound. The choice
of $\tau$ gives
\[
 3MN\tau^{-q/2}=\frac{\varepsilon\rho(1-2\rho)}2,\qquad
 W_0=\frac{\varepsilon(1-2\rho)(1-\rho)}2.
\]
At time $N$, $W_N=F(v_N,x_N)$, so conditional Markov implies
\begin{equation}\label{eq:ci-heattail}
 \Pp\{F(v_N,x_N)\ge1-\rho\mid\F_k\}
 \le\frac{\varepsilon(1-2\rho)}2
\end{equation}
for each sign.

\noindent\textbf{Step 2: prevent freezing for localized candidates.}
Martingale orthogonality gives
$\E[(z-\mu)^2\mid\F_k]\le V/N$. Therefore
$\Pp\{|z-\mu|>R\mid\F_k\}\le2\varepsilon\rho$.
On the observable variance check, every candidate satisfying $|z-m|\le R$
obeys
\[
 \sum_{i=1}^N(Z_i-m)^2
 =Q+N(z-m)^2\le Q+NR^2\le A-2\tau.
\]
All partial sums of these nonnegative squares satisfy the same bound, so no
update freezes for any such candidate. Its terminal potential is consequently
\[
 \Phi\left(\frac{\pm N(z-m)-b}
 {\sqrt{A-Q-N(z-m)^2}}\right).
\]
\noindent\textbf{Step 3: invert and account for error.}
Bounding the terminal denominator by $\sqrt{A-Q}$ removes its dependence
on the candidate mean. Specifically, if $|z-m|\le R$ but
\[
 |z-m|>\frac{b+z_{1-\rho}\sqrt{A-Q}}N,
\]
the appropriate signed potential exceeds $1-\rho$, because its denominator
is at most $\sqrt{A-Q}$. Thus exclusion of the true mean is contained in the
localization failure event or one of the two events in
\eqref{eq:ci-heattail}. Their total conditional probability is at most
$2\varepsilon\rho+\varepsilon(1-2\rho)=\varepsilon$.
If the variance check fails, the reported radius is $R$ and the coarse
guarantee applies directly.
\end{proof}

We record a moment-rate fact to make the almost-sure efficiency argument
self-contained, including its later $p$-moment application.
\begin{lemma}[A sufficient strong-law rate]\label{lem:ci-lawrate}
If $Y_i$ are iid with $\E|Y_1|^r<\infty$, $1<r\le2$, then for every
$0<d<(r-1)/2$,
$m^{-1}\sum_{i\le m}Y_i-\E Y_1=O_{\mathrm{a.s.}}(m^{-d})$.
Consequently, under the iid assumptions in Theorem~\ref{thm:ci-heat},
\begin{equation}\label{eq:ci-variancerates}
 \widehat v-\sigma^2=O_{\mathrm{a.s.}}(n^{-d_0}),\qquad
 Q/N-\sigma^2=O_{\mathrm{a.s.}}(n^{-d_0})
\end{equation}
for some $d_0>0$.
\end{lemma}
\begin{proof}
Truncation gives a bounded-variable maximal bound on each dyadic block.
The extra moment makes both truncation failures and centered fluctuations
summable, so Borel--Cantelli turns them into an almost-sure rate.
For dyadic $B=2^j$, truncate $Y_i$ to $Y_i^{(B)}=Y_i\mathbf1\{|Y_i|\le B\}$.
Writing $C=\E|Y_1|^r$, the probability that any of the first $2B$ terms is
changed is at most $2CB^{1-r}$. Their truncated variances are at most
$CB^{2-r}$. Kolmogorov's maximal inequality bounds the probability that a
centered truncated partial sum through $2B$ exceeds $B^{1-d}$ by
$2CB^{1-r+2d}$. Both bounds are summable over dyadic $B$. The truncation bias
per term is at most $CB^{1-r}$. Borel--Cantelli and division by $m\in[B,2B]$
give the asserted rate. The maximal inequality itself follows by partitioning
according to the first crossing; independence makes the future centered sum
orthogonal to the stopped partial sum.

Apply the rate to $(X_i-\mu)^2$ with $r=q/2>1$, and to $X_i-\mu$ with $r=2$.
The squared empirical-mean correction is negligible. Since
$k\asymp\sqrt n$, the pilot has a polynomial variance rate. Subtracting the
first $k$ partial sums from the first $n$ proves a polynomial rate for the
scored variance as well. This gives \eqref{eq:ci-variancerates}.
\end{proof}

\begin{proof}[Proof of Gaussian efficiency in Theorem~\ref{thm:ci-heat}]
We now assume iid sampling; this assumption was not used for coverage.
The radius has a Gaussian cutoff term and a leftover-variance term. The
former must converge to the desired quantile, while the latter must vanish
at the $n^{-1/2}$ scale.
\begin{itemize}
 \item \textbf{Step 1: learn the variance budget.} Strong-law rates make
 the observable check eventually pass and identify the remaining variance.
 \item \textbf{Step 2: evaluate the radius.} The leftover-variance term
 vanishes, leaving the true-variance Gaussian cutoff.
 \item \textbf{Step 3: restore the full-sample center.} The coarse radius
 and the omitted pilot have no first-order effect.
\end{itemize}

\noindent\textbf{Step 1: learn the variance budget.}
By Lemma~\ref{lem:ci-lawrate}, both variance errors are $o(\rho)$, so the
observable check eventually passes almost surely. Since
$\varepsilon_n\to\varepsilon_\infty>0$, the deterministic terms obey
\[
 \tau/N=O(N^{2/q-1}\rho^{-2/q})=o(\rho),\qquad
 R^2=O((N\rho)^{-1})=o(\rho).
\]
It follows that
\[
 A/N\to\sigma^2,\qquad
 (A-Q)/N=\widehat v-Q/N+8V\rho+R^2+2\tau/N\sim8V\rho.
\]
\noindent\textbf{Step 2: evaluate the radius.}
The remaining variance is of order $N\rho$, while the extra Gaussian
quantile grows only logarithmically in $1/\rho$.
Indeed, $\varepsilon_n(1-2\rho)^2/2\to\varepsilon_\infty/2$ and
$z_{1-\rho}=O(\sqrt{\log(1/\rho)})$. Consequently
\[
 \sqrt N\,\frac{b+z_{1-\rho}\sqrt{A-Q}}N
 =\sqrt{A/N}\,z_{1-\varepsilon_n(1-2\rho)^2/2}
   +\sqrt{(A-Q)/N}\,z_{1-\rho}
 \longrightarrow\sigma z_{1-\varepsilon_\infty/2}.
\]
\noindent\textbf{Step 3: restore the full-sample center.}
The coarse alternative has
$\sqrt N R=\sqrt{V/(2\varepsilon_n\rho)}\to\infty$ and hence does not alter
the limit. Finally,
\[
 \sqrt n(z-\bar X_n)
 =\frac{k}{\sqrt n}(z-\bar X_k)\longrightarrow0
\]
by the strong law and $k/\sqrt n\to1$. Since $N/n\to1$, the claimed limits
follow.
\end{proof}

The sub-Gaussian result is now a direct specialization, with no new
calibration argument. Taking $q=3$, $M=6\gamma s^3$, and $V=s^2$ specializes this construction
to \eqref{eq:ci-heat-simple-reserves}--\eqref{eq:ci-heat-simple-radius}
and proves Theorem~\ref{thm:ci-heat-main}. The third-moment bound follows by
applying Lemma~\ref{lemma:cubic_inequality} to both signs.

\section{Vector-valued constructions and proofs}\label{app:multi-details}

This appendix supports Section~\ref{sec:multi}. The constructions for a
vector mean and for adaptive regression follow the same progression:
first obtain a valid bound for a centered vector sum, then translate that
bound into an ellipsoid for the unknown parameter. The sample-scatter
mixture uses observed variation for practical width; the heat construction
adds the exact Gaussian cutoff. For regression, the centered vector sum
is a score, and the inverse design converts score uncertainty into
parameter uncertainty.

We first give the multivariate-mean calculations and then the regression
recipes and proofs in Appendix~\ref{app:linear-bandits}. In each proof,
finite-sample coverage uses the conditional model alone. Covariance
convergence and stable designs enter only when establishing efficiency.

\subsection{Optimizing the normal-mixture benchmark}\label{app:multi-normal}

The benchmark uses the supplied proxy in every direction. Its only
choice is how widely to spread the Gaussian bets. At a planned horizon,
we can make that choice before observing the scored data and optimize it
exactly. For a Gaussian betting prior with precision $\gamma Ns^2I_d$, the usual
proxy-based mixture gives the radius
\[
 \frac{s}{\sqrt N}
 \sqrt{(1+\gamma)\{2\log(1/b)+d\log(1+1/\gamma)\}},\qquad\gamma>0.
\]
The optimizing scale depends only on $d$ and $b$. Differentiating in $\gamma$ gives
$\gamma^{-1}-\log(1+\gamma^{-1})=2\log(1/b)/d$.
At the optimum the expression under the square root is
$d(1+\gamma^{-1})$, yielding the analytic formula
\begin{equation}\label{eq:multi-normal-closed}
 w(N,d,b,s)=s\sqrt{-\frac dN\,W_{-1}\!\left(-\exp\{-1-2\log(1/b)/d\}\right)}.
\end{equation}
Here $W_{-1}$ is the lower real branch of the Lambert $W$ function. This
evaluates the minimum in \eqref{eq:multi-normal-radius} directly. Using
this horizon-optimized benchmark avoids a logarithmic penalty from a
prior scale held fixed across all horizons.

\subsection{A capped Gaussian mixture of the sample-scatter processes}
\label{app:multi-mixture}

The key to explicit inversion is that the sample scatter is shared by
all candidate means. We average directional bets with a Gaussian prior,
while capping their lengths to respect Theorem~\ref{thm:qprocess}.
The cap complicates the Gaussian integral, but a lower bound on a small
ball of allowed bets will suffice; we never need to evaluate the mixture.

\paragraph{Pilot and scoring statistics.}
For a planned sample size $n\ge4$, use
$k=\max\{2,\lfloor n^{2/5}\rfloor\}<n$ pilot observations and score the
remaining $N=n-k$. Define the scoring mean and scatter by
\[
 z=\frac1N\sum_{j=1}^N X_{k+j},\qquad
 Q=\sum_{j=1}^N(X_{k+j}-z)(X_{k+j}-z)^\top.
\]
Compute the pilot covariance
$k^{-1}\sum_{i=1}^k(X_i-\bar X_k)(X_i-\bar X_k)^\top$ and clip its
eigenvalues to $[s^2/\log(n+e),s^2]$ to obtain $B\succ0$.
These same statistics are used by the heat construction below. No pilot
accuracy guarantee is needed for coverage. For $n<4$, the full-sample
normal-mixture ball is a valid fallback.

\paragraph{The mixture and its quadratic bound.}
Condition on the pilot. Let $R\succ0$ be a pilot-selected prior precision,
and let $0<a\le1/(5s)$ be a cap chosen before scoring. Write
$\pi_R=N_d(0,R^{-1})$. At scoring time $t$, let $Q_t$ be the sample
scatter of the first $t$ scored observations and let
$H_t=\sum_{i=1}^t i^{-1}$, with $Q_0=0$ and $H_0=0$.
The exact process used by the construction is
\begin{equation}\label{eq:multi-q-mixture}
 \begin{split}
 M_t(m)={}&\int_{\|\theta\|\le a}
 \exp\!\left\{\theta^\top\sum_{j=1}^t(X_{k+j}-m)
       -\frac12\theta^\top(Q_t+s^2H_tI_d)\theta
       -s^3(t-H_t)\|\theta\|^3\right\}\,\pi_R(d\theta)\\
 &+\pi_R\{\|\theta\|>a\},\qquad 0\le t\le N.
 \end{split}
\end{equation}
The prior mass outside the allowed bets is held at wealth one. Consequently
$M_0=1$, and Theorem~\ref{thm:qprocess}, applied to each projection
$\theta^\top X_i$, makes $M_t(\mu)$ a nonnegative test supermartingale.
The cap is essential: the theorem does not justify an unrestricted
Gaussian mixture. Neither the integral nor the excluded prior mass needs
to be evaluated by the reported procedure.

For a concrete default, use
\[
 a=\min\left\{\frac1{5s},\frac1{sN^{1/4}}\right\},\qquad
 R=\frac{N\operatorname{tr}(B)}{\sqrt d}\,I_d,
\]
with the clipped pilot covariance $B$ defined above.
To pass from the exact process to a computable ellipsoid, we first
replace the cubic charge by a quadratic upper bound on the allowed bets.
This folds the charge into the same matrix as the observed scatter.
At the planned horizon, write $Q=Q_N$ and define
\begin{equation}\label{eq:multi-practical-stats}
 \begin{aligned}
 V&=Q+s^2H_NI_d+2s^3(N-H_N)aI_d,\\
 A&=R+V,\qquad D=\log\det A-\log\det R.
 \end{aligned}
\end{equation}
The extra isotropic term in $V$ bounds the cubic repair uniformly over
the allowed bets. It is a deterministic charge given the pilot, and it
vanishes relative to $N$ for the default shrinking cap.

The following finite calculation turns the mixture into an ellipsoid.
After completing the square, the integral has covariance $A^{-1}$.
We will retain a ball of radius $a/2$ under that Gaussian law. The next
two quantities express a conservative lower bound on its probability:
$\kappa_a$ is the radius of a standard Gaussian ball known to fit, and
$p_a$ is its probability. Let $F_{\chi_d^2}$ denote the chi-square
distribution function, and put
\begin{equation}\label{eq:multi-cap-calibration}
 \kappa_a=\frac a2\sqrt{\lambda_{\min}(A)},\qquad
 p_a=F_{\chi_d^2}(\kappa_a^2).
\end{equation}
Here $0<p_a<1$. Define the explicit cutoff
\begin{equation}\label{eq:multi-H}
 q_b(A,R,a)=
 \begin{cases}
 \sqrt{2\log(1/b)+D-2\log p_a},
       &2\log(1/b)+D-2\log p_a\le\kappa_a^2,\\[4pt]
 \displaystyle\frac{\log(1/b)+D/2-\log p_a}{\kappa_a}
                  +\frac{\kappa_a}{2},&\text{otherwise}.
 \end{cases}
\end{equation}
The first branch is the usual Gaussian-mixture quadratic inversion.
The second accounts for the cap when the required bet would otherwise
be too large. Both branches use only matrix operations, logarithms,
square roots, and a chi-square distribution-function evaluation.
No candidate search, mixture integration, or numerical optimization is
needed. For integer $d$, the chi-square distribution function can also be
evaluated by finite recurrences from exponential and Gaussian functions.

\begin{proposition}[Sample-scatter mixture ellipsoid]
\label{prop:multi-practical}
Under \eqref{eq:multi-model}, fix $b\in(0,1)$ and choose $R\succ0$ and
$0<a\le1/(5s)$ using only the pilot. With the preceding definitions, the ellipsoid
\begin{equation}\label{eq:multi-practical-ellipse}
 \mathcal P_b=\{m:N^2(m-z)^\top A^{-1}(m-z)
                        \le q_b(A,R,a)^2\}
\end{equation}
has conditional coverage at least $1-b$ given the pilot. It is a global
outer bound on $\{m:M_N(m)\le1/b\}$ and needs no localization event.
\end{proposition}

\begin{proof}[Proof of Proposition~\ref{prop:multi-practical}]
The mixture is already a valid test process. The remaining task is
deterministic: show that a candidate far enough from the scoring mean
must make the mixture large, even with the bet cap.
\begin{itemize}
 \item \textbf{Step 1: keep a useful ball of bets.} Replace the cubic
 charge by its quadratic bound, and shift a small ball toward the
 candidate's most favorable Gaussian bet.
 \item \textbf{Step 2: lower-bound its contribution.} Gaussian symmetry
 removes the remaining linear fluctuation and leaves a scalar bound in
 the candidate's squared Mahalanobis distance.
 \item \textbf{Step 3: invert and calibrate.} Solve the scalar bound
 explicitly and apply conditional Markov's inequality.
\end{itemize}

\noindent\textbf{Step 1: keep a useful ball of bets.}
Conditional validity of \eqref{eq:multi-q-mixture} follows from conditional
mixing of the projected sample-scatter supermartingales. For the analytic
inversion, fix an arbitrary candidate $m$ and write
$S=N(z-m)$. Since $\|\theta\|^3\le a\|\theta\|^2$ on the cap,
\[
 M_N(m)\ge\int_{\|\theta\|\le a}
       e^{\theta^\top S-\theta^\top V\theta/2}\,\pi_R(d\theta).
\]
Combining the prior precision with the quadratic charge gives the
matrix $A=R+V$. Put $h=A^{-1}S$ and $r=S^\top A^{-1}S$. For $r>0$, choose
$\eta=\min\{1,\kappa_a/\sqrt r\}$; for $r=0$, take $\eta=1$.
Because $\|h\|\le\sqrt{r/\lambda_{\min}(A)}$, the shifted ball
\[
 \{\eta h+g:\|g\|\le a/2\}
\]
lies entirely inside the allowed bet ball. Restricting the integral to
this shifted ball and completing the square yields
\[
 M_N(m)\ge e^{-D/2+(\eta-\eta^2/2)r}
 \E\!\left[e^{(1-\eta)G^\top S}
             \ind\{\|G\|\le a/2\}\right],
 \qquad G\sim N_d(0,A^{-1}).
\]
The shifted ball is used only to bound the integral; the mixture prior and
bet cap remain fixed before scoring.
\smallskip
\noindent\textbf{Step 2: lower-bound its contribution.}
The centered Gaussian law and the ball are symmetric. Replacing the
exponential inside the expectation by its average at $G$ and $-G$
therefore gives a hyperbolic cosine, which is at least one. Moreover,
writing $G=A^{-1/2}Z$ with $Z\sim N_d(0,I_d)$ gives
\[
 \Pp(\|G\|\le a/2)
 \ge\Pp(\|Z\|\le\kappa_a)=p_a.
\]
Thus the following bound holds for every candidate, without localization:
\[
 \log M_N(m)\ge-D/2+\log p_a+
 \begin{cases}
 r/2,&\sqrt r\le\kappa_a,\\
 \kappa_a\sqrt r-\kappa_a^2/2,&\sqrt r>\kappa_a.
 \end{cases}
\]
\smallskip
\noindent\textbf{Step 3: invert and calibrate.}
The bound is quadratic near the center and linear farther out: once the
Gaussian optimum would leave the cap, we move only as far as the cap
allows. Inverting this increasing function of $\sqrt r$ at level $\log(1/b)$
gives exactly \eqref{eq:multi-H}. Hence every candidate satisfying
$M_N(m)\le1/b$ belongs to
\[
 \{m:N^2(z-m)^\top A^{-1}(z-m)\le q_b(A,R,a)^2\}.
\]
Conditional Markov's inequality at the true mean proves error at most
$b$. There is no separate localization event and no error split.
\end{proof}

\paragraph{What covariance adaptation means here.}
The ellipsoid uses a regularized sample covariance. To identify its
limiting geometry, fix an iid law with $\Sigma\succ0$. The pilot and
scoring sample covariances are strongly consistent, $N/n\to1$, and
$a\sim s^{-1}N^{-1/4}$. Therefore, almost surely,
\[
 V/N\longrightarrow\Sigma,\qquad
 R/N\longrightarrow\frac{\operatorname{tr}(\Sigma)}{\sqrt d}\,I_d,
 \qquad
 A/N\longrightarrow\Sigma+
             \frac{\operatorname{tr}(\Sigma)}{\sqrt d}\,I_d.
\]
Also $\kappa_a^2$ has order $N^{1/2}$, so $p_a\to1$ and the quadratic
branch of \eqref{eq:multi-H} is eventually used at every fixed error
level $b$. The regularized shape $A$ should not be confused with the true
covariance itself. This practical mixture does not assert the exact
chi-square constant; that calibration is supplied by the matrix heat
construction.

\subsection{Matrix heat drift and Gaussian calibration}\label{app:multi-heat}

The heat argument has two jobs. A declining reserve pays for the
continuation probability's possible upward drift, giving a
supermartingale. A coarse localization ball then ensures that, for the
candidates we need to protect, the observed covariance never exhausts
the initial budget. This makes the terminal state explicit.

Accordingly, $W_G$ below is a coarse radius in the original coordinates,
$M_B$ bounds a whitened observation's third moment, and $\tau$ is the
smoothing floor. The initial budget $A$ covers nominal covariance,
localization, and smoothing. Use the same slack
$\rho=1/[8\log(n+e)]$ and the pilot covariance $B$
defined in Appendix~\ref{app:multi-mixture}. At component error level $b$, set
\begin{equation}\label{eq:multi-heat-tuning}
 \begin{aligned}
 W_G&=s\sqrt{\frac d{2bN\rho}},&
 Q_B&=B^{-1/2}QB^{-1/2},\\
 M_B&=\frac{4d^{3/2}s^3}{\lambda_{\min}(B)^{3/2}},&
 \tau&=\left\{\frac{4M_BN}{b\rho(1-2\rho)}\right\}^{2/3},\\
 A&=N(1+8\rho)+\frac{NW_G^2}{\lambda_{\min}(B)}+2\tau.
 \end{aligned}
\end{equation}
The matrix $Q_B$ is observed after scoring; the initial parameters are
fixed given the pilot. Check
\begin{equation}\label{eq:multi-heat-audit}
 \lambda_{\max}(Q_B)\le N(1+8\rho).
\end{equation}
On this event use \eqref{eq:multi-heat-radius}; otherwise report the
radius-$W_G$ ball. These choices give
\[
 \frac{4M_BN}{\tau^{3/2}}=b\rho(1-2\rho),
\]
so the process \eqref{eq:multi-heat-process} is equivalently
$F(V_j,x_j)+4M_B(N-j)/\tau^{3/2}$. This is the drift-reserve form used
in the proof below.

\begin{proof}[Proof of Theorem~\ref{thm:multi-heat}]
The first three steps prove conditional coverage without requiring the
pilot covariance to be accurate. Only the last step assumes iid sampling.
\begin{itemize}
 \item \textbf{Step 1: control one heat update.} The matrix heat identity
 cancels the quadratic term; the guards keep the remaining derivatives
 bounded.
 \item \textbf{Step 2: pay the drift.} Use the third-moment bound to
 turn that remainder into a deterministic charge and obtain the
 supermartingale and its terminal tail bound.
 \item \textbf{Step 3: obtain the reported ellipsoid.} On the coarse
 localization event and observable audit, no update freezes. The
 Gaussian continuation then excludes every candidate outside the
 displayed radius.
 \item \textbf{Step 4: identify the iid limit.} The pilot learns the
 covariance, the audit eventually passes, and the leftover covariance
 contributes a vanishing part of the radius.
\end{itemize}

\noindent\textbf{Step 1: control one heat update.}
Condition on the pilot. The initial budget $A$, smoothing reserve $\tau$,
pilot covariance $B$, and terminal cutoff in $F$ are then fixed. Gaussian
convolution satisfies the rank-one heat identity
\[
 \frac{d}{dt}F(V+tuu^\top,x)\bigg|_{t=0}
 =\tfrac12D_x^2F(V,x)[u,u].
\]
For $V\succeq\tau I_d$, differentiating the Gaussian kernel and applying
Cauchy--Schwarz to the third and fourth Hermite polynomials gives
\begin{align*}
 |D_x^3F(V,x)[u^3]|&\le\sqrt6\,\tau^{-3/2}\|u\|^3,\\
 |D_x^4F(V,x)[u^4]|&\le\sqrt{24}\,\tau^{-2}\|u\|^4,
 \qquad \|\nabla F(V,x)\|\le\tau^{-1/2}.
\end{align*}
At an active state, $V\succeq2\tau I_d$. If $\|u\|\le\sqrt\tau$, the
proposed update is accepted, and all intermediate covariance matrices are
at least $\tau I_d$. A spatial Taylor expansion followed by a covariance
Taylor expansion cancels the quadratic term and gives
\[
 |\Delta F-\nabla F(V,x)^\top u|
 \le\left(\frac{2\sqrt6}{3}+\frac{\sqrt{24}}8\right)
          \tau^{-3/2}\|u\|^3
 <3\tau^{-3/2}\|u\|^3.
\]
For larger accepted increments, boundedness of $F$ gives
$|\Delta F-\nabla F(V,x)^\top u|
\le1+\|u\|/\sqrt\tau\le2\tau^{-3/2}\|u\|^3$.
The same bound covers a rejected increment, for which $\Delta F=0$.
At an already frozen state the actual potential increment is zero.

\smallskip
\noindent\textbf{Step 2: pay the drift.}
At an active state, conditional centering removes the linear part of
the update. Thus only the absolute third-moment remainder needs to be
paid. Under the null, the proposed increment
$u=B^{-1/2}(X_{k+j}-\mu)$ is conditionally centered, and
\[
 \E[\|u\|^3\mid\mathcal F_{k+j-1}]\le M_B.
\]
Indeed, each coordinate is
$s\lambda_{\min}(B)^{-1/2}$-sub-Gaussian, its absolute third moment is
at most four times the cube of that proxy, and
$\|u\|^3\le\sqrt d\sum_l|u_l|^3$.
The conditional potential drift is therefore at most
$4M_B\tau^{-3/2}$, including the rejected and frozen cases. The reserve in
\eqref{eq:multi-heat-process} decreases by this amount at each step, so
$(W_j(\mu))_{j=0}^N$ is a nonnegative supermartingale. Its conditional
integrability follows from $0\le F\le1$ and the finite initial reserve.

The choices of cutoff and smoothing reserve give
\[
 F(AI_d,0)=b(1-2\rho)^2,\qquad
 \frac{4M_BN}{\tau^{3/2}}=b\rho(1-2\rho),\qquad
 W_0=b(1-2\rho)(1-\rho).
\]
Since the reserve vanishes at time $N$, conditional Markov implies
\[
 \Pp\{F(V_N,x_N)\ge1-\rho\mid\mathcal F_k\}
 \le b(1-2\rho).
\]
\smallskip
\noindent\textbf{Step 3: obtain the reported ellipsoid.}
The terminal tail bound is useful once the candidate's terminal sum and
covariance can be written in terms of the observed mean and scatter.
Localization and the audit provide exactly that guarantee.
Martingale orthogonality and \eqref{eq:multi-model} give
$\E[\|z-\mu\|^2\mid\mathcal F_k]\le ds^2/N$.
Consequently, the radius-$W_G$ ball misses $\mu$ with conditional
probability at most $2b\rho$.

On the observable audit, every candidate in that ball satisfies
\[
 \sum_{j=1}^N u_ju_j^\top
 =Q_B+NB^{-1/2}(z-m)(z-m)^\top B^{-1/2}
 \preceq(A-2\tau)I_d.
\]
Every partial sum of these positive semidefinite matrices satisfies the
same bound. Hence no update freezes for such a candidate, and
\[
 x_N=NB^{-1/2}(z-m),\qquad
 V_N\preceq\{A-\lambda_{\min}(Q_B)\}I_d.
\]
If $\|B^{-1/2}(z-m)\|>r_G$, then on the Gaussian event
$\|Z\|\le c_d(\rho)$ the triangle inequality gives
\[
 \|x_N+V_N^{1/2}Z\|
 >\sqrt A\,c_d\bigl(b(1-2\rho)^2\bigr).
\]
Thus $F(V_N,x_N)\ge1-\rho$. Exclusion of the true mean is contained in
the localization failure event or the terminal heat event, with total
conditional probability at most $2b\rho+b(1-2\rho)=b$.
If the audit fails, the reported ball has error at most $2b\rho\le b$;
its failure is already included in the same localization event.

\smallskip
\noindent\textbf{Step 4: identify the iid limit.}
For efficiency, iid sub-Gaussian moment bounds give polynomial almost-sure
rates for the pilot and scoring covariances. The pilot has
$k\asymp n^{2/5}$ observations, and clipping its eigenvalues preserves
these rates because $\Sigma\succ0$ and $\Sigma\preceq s^2I_d$. Therefore
\[
 B-\Sigma=o(\rho),\qquad Q_B/N-I_d=o(\rho)
 \quad\text{almost surely}.
\]
The audit eventually passes. When $b\to\alpha\in(0,1)$,
$W_G^2=O((N\rho)^{-1})$, $M_B=O(1)$ eventually, and
$\tau/N=O(N^{-1/3}\rho^{-2/3})=o(\rho)$. It follows that
\[
 A/N\to1,\qquad
 \{A-\lambda_{\min}(Q_B)\}/N=8\rho+o(\rho).
\]
Since $c_d(\rho)=O(\sqrt{\log(1/\rho)})$ for fixed $d$, the remaining
covariance term in $\sqrt N r_G$ tends to zero. Its cutoff term tends to
$c_d(\alpha)$ because $b(1-2\rho)^2\to\alpha$.
This proves $B\to\Sigma$ and $\sqrt N r_G\to c_d(\alpha)$ almost surely.
Also $N/n\to1$ and
\[
 \sqrt n(z-\bar X_n)=\frac{k}{\sqrt n}(z-\bar X_k)\longrightarrow0
 \quad\text{almost surely},
\]
since $k\asymp n^{2/5}$ and both sample means are strongly consistent.
Thus recentering at the full-sample mean and rescaling by $\sqrt n$
gives the Hausdorff limit in \eqref{eq:multi-goal}.
\end{proof}

\subsection{Adaptive regression and linear-bandit scores}\label{app:linear-bandits}

We now apply the vector constructions to the regression scores of
Section~\ref{subsec:linear-bandits}. The action can depend on previous
observations, but the score at the true parameter remains conditionally
centered. This is the property needed by both the exponential and heat
arguments. The design matrix then maps a score bound back into parameter
space, producing the usual sandwich covariance in the limit.

The forecast clock and score are defined in
\eqref{eq:lb-observed-scores}. Throughout this part, write
$G_n=D_n$ and $v_i=\E[\varepsilon_i^2\mid\F_{i-1}]$.

\subsubsection{Stability conditions and the normal-mixture benchmark}
\label{app:lb-benchmark}

Coverage is valid for arbitrary predictable bounded actions. Efficiency
requires more: both the information supplied by the actions and the
noise covariance of the scores must settle to full-rank limits. The
planned-time heat calibration also needs their errors to disappear
faster than its logarithmic slack. The following polynomial condition
ensures that separation: for
deterministic $D_\infty,K_\infty\succ0$ and some $\kappa>0$,
\begin{equation}\label{eq:lb-stability}
 \|D_n/n-D_\infty\|+\|K_n/n-K_\infty\|=O(n^{-\kappa})
 \quad\text{almost surely}.
\end{equation}
This rules out the design degeneracies under which an ordinary regression
Gaussian limit need not hold~\cite{Deshpande2018,Whitehouse2025}.
The score covariance $K_n$ permits
heteroskedastic noise; under conditional homoskedasticity,
$K_\infty=\sigma^2D_\infty$.

For orientation, in model~\eqref{eq:lb-model} put $V_n=\rho I_d+D_n$ and
$\widehat\theta_n^{\rho}=V_n^{-1}\sum_{i\le n}x_iY_i$. If
$\|\theta_*\|\le S_0$ is known, the Abbasi-Yadkori, P\'al, and
Szepesv\'ari confidence sequence~\cite{abbasi2011} has the familiar form
\[
 \|\theta-\widehat\theta_n^{\rho}\|_{V_n}
 \le s\sqrt{\log\frac{\det V_n}{\rho^d}+2\log(1/\alpha)}
                +\sqrt\rho S_0.
\]
For a fixed $\rho>0$ and $D_n/n\to D_\infty\succ0$, its parameter radius
has order $s\sqrt{\log n/n}$. We target instead the true-variance
chi-square ellipsoid at a planned time and the true-variance
$\sqrt{2\log\log n/n}$ shape over time. Variance-adaptive linear-bandit
confidence sets are already studied, for example by Jun and
Kim~\cite{JunKim2024}; the present goals are these exact limiting shapes
and explicit reported ellipsoids. No improved regret theorem is asserted.

\subsubsection{The planned-time recipe and its proof}\label{app:lb-fixed}

The planned-time construction has a simple order. First use a standard
score mixture to localize the parameter. Next check that the design is
invertible and that all localized candidates fit inside a common
covariance budget. Only then use the heat ellipsoid; otherwise keep the
score-mixture fallback. The pilot chooses the coordinates and initial
budgets before scoring, so adaptive actions do not invalidate the
calibration.

Use the first $k=\max\{2,\lfloor n^{2/5}\rfloor\}<n$ observations as a pilot
and let $N=n-k$. All sums in the
following display run over the scored observations $i=k+1,\ldots,n$:
\begin{equation}\label{eq:lb-score-stats}
 G=\sum x_ix_i^\top,\quad h=\sum x_iY_i,\quad
 S(\theta)=h-G\theta.
\end{equation}
From the pilot compute the minimum-norm least-squares estimate
$\widehat\theta_P=D_k^\dagger\sum_{i\le k}x_iY_i$ and the residual score
covariance $k^{-1}\sum_{i\le k}x_ix_i^\top(Y_i-x_i^\top\widehat\theta_P)^2$.
Project its eigenvalues onto $[s^2L^2/\log(n+e),s^2L^2]$ to obtain $B\succ0$.
Its finite-sample accuracy is not assumed.

The next calibration allocates error $\beta$ to localization and the
remaining error $a$ to the heat bound. The threshold $g$ checks the
smallest design eigenvalue; it turns a score bound into the radius $R$
in parameter space. The quantity $T$ pays for replacing a candidate's
residual covariance by the observable least-squares residual covariance.
Finally, $M_B$, $\tau$, and $A$ play the same moment, smoothing, and
initial-budget roles as in the vector-mean heat construction.
For an error budget $b\in(0,1)$, set
\begin{align}
 \beta&=b/(4\log(n+e)),\quad a=b-\beta,\quad
 \eta=1/(8\log(n+e)),\quad a_0=a(1-2\eta),\nonumber\\
 C&=d\log(1+NL^2/d)+2\log(1/\beta),\quad g=N/\log(n+e),\nonumber\\
 R^2&=s^2C(1+1/g)/g,\quad T=NL^4R^2/\lambda_{\min}(B),
 \label{eq:lb-local-tuning}\\
 M_B&=4s^3L^3\lambda_{\min}(B)^{-3/2},\quad
 \tau=\{4M_BN/(a\eta)\}^{2/3},\nonumber\\
 A&=N(1+1/\log(n+e))^2+(1+\log(n+e))T+2\tau.\label{eq:lb-heat-A}
\end{align}
\paragraph{Check the design and observed covariance.}
The first audit is $\lambda_{\min}(G)\ge g$. If it passes, compute
\[
 \widehat\theta=G^{-1}h,\quad
 \widehat Q=\sum x_ix_i^\top(Y_i-x_i^\top\widehat\theta)^2,\quad
 \widehat Q_B=B^{-1/2}\widehat Q B^{-1/2}.
\]
The second audit and the heat-selection gate are
\begin{equation}\label{eq:lb-audits}
 \lambda_{\max}(\widehat Q_B)\le N(1+1/\log(n+e)),\qquad 2\tau/N\le1/\log(n+e).
\end{equation}
\paragraph{Report the heat ellipsoid or the fallback.}
When the checks pass, $q_-$ below is a common lower bound on the
covariance spent by every localized candidate. Subtracting it from $A$
therefore bounds the remaining Gaussian continuation covariance.
If both audits and the gate pass, set
\begin{align}
 q_-&=(1-1/\log(n+e))\lambda_{\min}(\widehat Q_B)-(\log(n+e)-1)T,\nonumber\\
 r&=\sqrt A\,c_d(a_0)+\sqrt{A-q_-}\,c_d(\eta),\qquad
 c_d(u)=\sqrt{\chi^2_{d,1-u}},\label{eq:lb-heat-cutoff}
\end{align}
and report the closed-form ellipsoid
\begin{equation}\label{eq:lb-heat-set}
 \mathcal L_{n,b}=
 \{\theta:(\theta-\widehat\theta)^\top G B^{-1}G
                (\theta-\widehat\theta)\le r^2\}.
\end{equation}
If either audit or the gate fails, report instead the score-mixture set
\begin{equation}\label{eq:lb-fallback}
 \mathcal L_{n,b}=
 \{\theta:S(\theta)^\top(I_d+G)^{-1}S(\theta)\le s^2C\}.
\end{equation}
This is an ellipsoid whenever $G$ is invertible and otherwise a cylinder.
For sample sizes too small to form the pilot, use
\eqref{eq:lb-fallback} on all observations with
$C=d\log(1+nL^2/d)+2\log(1/b)$.

The fixed-time recipe costs $O(nd^2+d^3)$ arithmetic operations and uses
two passes over the scored observations, or stores them for the residual
covariance calculation.

\begin{proof}[Proof of Theorem~\ref{thm:lb-fixed}]
The difficulty is that a candidate changes the residual covariance, and
the design is adaptive. We handle both issues with a score-based
localizer fixed in radius before scoring, followed by deterministic
covariance comparisons. The proof proceeds as follows.
\begin{itemize}
 \item \textbf{Step 1: localize with a score mixture.} Obtain a valid
 parameter-radius bound whenever the design audit passes, without
 assuming a bound on the parameter itself.
 \item \textbf{Step 2: compare residual covariances.} Bound the
 candidate-dependent covariance above and below using least-squares
 residuals and the localization radius.
 \item \textbf{Step 3: apply the heat calibration.} On the audited,
 localized event, invert the terminal heat bound and combine the two
 error budgets; retain the mixture fallback on failure of a check.
 \item \textbf{Step 4: recover the sandwich limit.} Under the stated
 stability condition, the checks eventually pass and the additional
 reserves vanish at the Gaussian scale.
\end{itemize}

\noindent\textbf{Step 1: localize with a score mixture.}
Conditional on the pilot, the natural score-mixture e-value is
\[
 \det(I_d+G)^{-1/2}
 \exp\{S(\theta_*)^\top(I_d+G)^{-1}S(\theta_*)/(2s^2)\}.
\]
It follows by mixing the directional exponential supermartingales, exactly
as in the Abbasi-Yadkori construction. Since
$\det(I_d+G)\le(1+NL^2/d)^d$, \eqref{eq:lb-fallback} covers with
conditional error at most $\beta$. On this event and the first audit,
\[
 \|\widehat\theta-\theta_*\|^2
 \le \frac{1+1/g}{g}
 S(\theta_*)^\top(I_d+G)^{-1}S(\theta_*)\le R^2.
\]
This supplies a deterministic-radius localizer without assuming a parameter
norm bound or conditioning on the realized adaptive design.

\smallskip
\noindent\textbf{Step 2: compare residual covariances.}
A candidate near least squares changes each residual by a bounded
amount. The following two-sided square bounds keep the resulting
covariance comparison explicit; their logarithmic slack tends to zero.
For any candidate $\theta$ in that localizer put
$\delta=\widehat\theta-\theta$ and
$Q(\theta)=\sum x_ix_i^\top(Y_i-x_i^\top\theta)^2$.
The scalar inequalities
\[
 (u+v)^2\le(1+1/\log(n+e))u^2+(1+\log(n+e))v^2,\qquad
 (u+v)^2\ge(1-1/\log(n+e))u^2-(\log(n+e)-1)v^2
\]
and $\sum x_ix_i^\top(x_i^\top\delta)^2\preceq NL^4R^2I_d$ give
\begin{align*}
 B^{-1/2}Q(\theta)B^{-1/2}
 &\preceq(1+1/\log(n+e))\widehat Q_B+(1+\log(n+e))TI_d,\\
 B^{-1/2}Q(\theta)B^{-1/2}
 &\succeq(1-1/\log(n+e))\widehat Q_B-(\log(n+e)-1)TI_d.
\end{align*}
On both audits these bounds are, respectively, at most $(A-2\tau)I_d$
and at least $q_-I_d$.

\smallskip
\noindent\textbf{Step 3: apply the heat calibration.}
The upper covariance bound prevents freezing, while the lower bound
controls how much covariance is left at the terminal time. Apply the
matrix heat argument of Theorem~\ref{thm:multi-heat} to the
score increments $u_i=B^{-1/2}x_i(Y_i-x_i^\top\theta)$, starting from
$(AI_d,0)$ with cutoff $\sqrt A c_d(a_0)$ and covariance floor $\tau I_d$.
At the true parameter these increments are conditionally centered and
\[
 \E[\|u_i\|^3\mid\F_{i-1}]
 \le4s^3L^3\lambda_{\min}(B)^{-3/2}=M_B.
\]
The same Taylor cancellation therefore gives the nonnegative
supermartingale $F(V_j,z_j)+(N-j)4M_B\tau^{-3/2}$.
Its initial value is $a(1-\eta)$, so the probability that its terminal heat
value is at least $1-\eta$ is at most $a$. On the localized, audited event
there is no freezing; its terminal score is
$z_N=B^{-1/2}G(\widehat\theta-\theta)$ and
$V_N\preceq(A-q_-)I_d$. The triangle inequality for the Gaussian
continuation then shows that every $\theta$ outside
\eqref{eq:lb-heat-set} has heat value at least $1-\eta$.
Combining the heat error $a$ and localizer error $\beta$ proves coverage,
including the audit- or gate-failure fallback. Crucially, $A,B,R$ were all fixed
before the scored observations; neither the realized design nor the residual
clock was used to choose the initial heat covariance.

\smallskip
\noindent\textbf{Step 4: recover the sandwich limit.}
Only this step uses the stability assumptions. Bounded actions and
conditional sub-Gaussian errors give polynomial almost-sure rates for
the martingale sums
$\sum x_i\varepsilon_i$ and
$\sum x_ix_i^\top(\varepsilon_i^2-v_i)$. For example, each matrix or
vector entry has uniformly bounded conditional second moment; the
martingale maximal inequality on dyadic blocks and Borel--Cantelli give
$O_{\mathrm{a.s.}}(n^{1/2+\epsilon})$ for every fixed $\epsilon>0$.
Together with
\eqref{eq:lb-stability}, least-squares residual expansion gives
\[
 B-K_\infty=o(1/\log(n+e)),\qquad
 G/N-D_\infty=o(1/\log(n+e)),\qquad
 \widehat Q_B/N-I_d=o(1/\log(n+e))
 \quad\text{almost surely}.
\]
The eigenvalue clipping is asymptotically inactive except possibly at the
upper endpoint, which cannot change these rates since
$K_\infty\preceq s^2L^2I_d$. Thus both audits eventually pass.
Furthermore $R^2=O((\log n)^2/n)$, $T=O((\log n)^2)$,
$\tau/N=O(N^{-1/3}\{\log(n+e)\}^{2/3})$, and
\[
 A/N\to1,\qquad (A-q_-)/N=O(1/\log(n+e)).
\]
The rate for $\tau/N$ ensures that the gate eventually passes as well.
Because $c_d(\eta)=O(\sqrt{\log\log(n+e)})$, we obtain
$r/\sqrt N\to c_d(b)$ almost surely. Consequently the shape matrix
$nr^2G^{-1}BG^{-1}$ tends to
$\chi^2_{d,1-b}D_\infty^{-1}K_\infty D_\infty^{-1}$.
Finally, deleting $k\asymp n^{2/5}$ pilot observations changes ordinary
least squares by $o(n^{-1/2})$ almost surely, using the same martingale
strong rates. This proves \eqref{eq:lb-gaussian-limit}. The standard
martingale central limit theorem applies under
\eqref{eq:lb-stability}: conditional variances converge to $K_\infty$ and
the sub-Gaussian bound supplies conditional Lindeberg.
\end{proof}

The homoskedastic specialization has
$\mathcal L_{n,b}=\widehat\theta_n+
\{z:z^\top D_nz\le\sigma^2\chi^2_{d,1-b}(1+o(1))\}$ in the corresponding
normalized Hausdorff sense. Thus it replaces the supplied proxy by the true
noise variance and the exponential-mixture cutoff by the Gaussian
chi-square cutoff. This fixed-time conclusion is not a time-uniform
guarantee, and it does not hold for every unstable or insufficiently excited
linear-bandit design.

\subsubsection{The confidence-sequence recipe and its proof}\label{app:lb-cs}

For continuous monitoring, a single heat calibration at a planned
horizon is replaced by simultaneous exponential bets. There are two
choices to protect: the whitening matrix determines the shape, and the
bet size determines the radius. We allocate error across countable
families of both, then select one of each using the observed clock.
This permits adaptation to the current data without substituting an
uncalibrated random bet into a fixed-bet process.

Use the observable clock $R_n$ in \eqref{eq:lb-observed-scores}, with any
finite predictable forecast. The forecast-uniform cubic inequality applied
conditionally with
$U=(u^\top x_i)\varepsilon_i$ and
$a=(u^\top x_i)x_i^\top(\theta_*-\widetilde\theta_{i-1})$
proves the test-supermartingale claim in
\eqref{eq:lb-forecast-process}.

The convergence assumptions in Theorem~\ref{thm:lb-cs} are
\begin{equation}\label{eq:lb-cs-stability}
 G_n/n\longrightarrow D_\infty\succ0,\qquad
 K_n/n\longrightarrow K_\infty\succ0,\qquad
 \widetilde\theta_n\longrightarrow\theta_*.
\end{equation}
The limits are deterministic and positive definite, and convergence is
almost sure. No polynomial rate is required.

\paragraph{Choose a covariance frame.}
Here is a fully explicit ellipsoid obtained from these processes. A
slowly refining matrix grid makes the chosen frame converge to the score
covariance. Its counting penalty is lower order than $\log\log n$.
The matrices and direction nets below serve only to establish a common
coverage event: no grid or net is enumerated by the algorithm. Put
$D=d(d+1)/2$ and, for integers $k\ge2$, let
\begin{align*}
 \mathcal H_k&=\{H=k^{-1}A:A=A^\top,\ A_{ij}\in\mathbb Z,
       \ |A_{ij}|\le k^2,\ \lambda_{\min}(H)\ge k^{-1}\},\\
 \mathsf K_k&=(2k^2+1)^D,\qquad \mathsf J_k=(1+2k)^d,
 \qquad \zeta_k=6/(\pi^2k^2).
\end{align*}
There are at most $\mathsf K_k$ matrices in $\mathcal H_k$, which contains $I_d$;
the unit sphere has a deterministic $1/k$-net of size at most $\mathsf J_k$.
At time $n$, choose
\[
 k_n=\max\{2,\lfloor\log\log(e^e+n)\rfloor\},\qquad
 H_n^\circ=\{\operatorname{round}(k_nR_n/n)+2dI_d\}/k_n.
\]
Round the symmetric entries symmetrically. Set $H_n=H_n^\circ$ if its
finite entries and eigenvalues satisfy the displayed definition of
$\mathcal H_{k_n}$, and otherwise set $H_n=I_d$.
\paragraph{Choose a bet within that frame.}
The matrix $H$ whitens the score. Its smallest eigenvalue determines a
safe cap on the scalar bet, while the largest eigenvalue of the
whitened clock bounds variation in all directions. For a
positive-definite matrix $H$ define
\[
 b_H=\frac{\sqrt{\lambda_{\min}(H)}}{5sL},\qquad
 c_H=\frac{s^3L^3}{\lambda_{\min}(H)^{3/2}},\qquad
 v_n(H)=\lambda_{\max}(H^{-1/2}R_nH^{-1/2}).
\]
Use the dense scalar grid of Section~3 with cap $b_H$:
\[
 m_j=\lceil\log(j+e)\rceil,\quad
 w_j=\log2\{1/\log(j+2)-1/\log(j+3)\},\quad
 \lambda_{H,j,h}=b_H2^{-j-h/m_j},\quad \pi_{j,h}=w_j/m_j.
\]
Here $j\ge0$, $0\le h<m_j$, and $\sum_{j,h}\pi_{j,h}=1$.
For a selected point $(j,h)$ put
\begin{equation}\label{eq:lb-cs-bound}
 B_n(H,\lambda)=\frac1\lambda\log\frac{\mathsf K_k\mathsf J_k}{\alpha\zeta_k\pi_{j,h}}
                      +\frac{\lambda v_n(H)}2+nc_H\lambda^2.
\end{equation}
The first term pays for the simultaneous choices; the next two are the
observed-variance and cubic charges. As in the scalar construction, we
balance the threshold and variance terms and cap the result to protect
against large early bets. An explicit selection with the desired limit
is obtained by setting
\begin{align*}
 g_n&=\log\log(e^e+n)+\log(1/\alpha)
                         +(2D+d+4)\log(k_n+1)+1,\\
 \lambda_*&=\min\left\{b_{H_n},
       \sqrt{\frac{2g_n}{v_n(H_n)\vee n^{-1}}},
       \left(\frac{g_n}{2nc_{H_n}}\right)^{1/3}\right\}.
\end{align*}
With $t=\log_2(b_{H_n}/\lambda_*)$, take
$j=\lfloor t\rfloor$, $h=\lfloor m_j(t-j)\rfloor$ and
$\lambda_n=\lambda_{H_n,j,h}$. Thus
$1\le\lambda_n/\lambda_*\le2^{1/m_j}$.
\paragraph{Map the score bound into parameter space.}
When $G_n$ is nonsingular, write
$\widehat\theta_n=G_n^{-1}\sum_{i\le n}x_iY_i$ and define
\begin{equation}\label{eq:lb-cs-ellipsoid}
 \mathcal C_n^{\rm LIL}=\left\{\theta:
 (\theta-\widehat\theta_n)^\top G_nH_n^{-1}G_n
 (\theta-\widehat\theta_n)
 \le\left\{\frac{B_n(H_n,\lambda_n)}{1-1/k_n}\right\}^2\right\}.
\end{equation}
Set $\mathcal C_n^{\rm LIL}=\R^d$ when $G_n$ is singular.

\begin{proof}[Proof of Theorem~\ref{thm:lb-cs}]
First prove one event that covers every frame, direction, bet, and time.
The reported choices then inherit coverage. For efficiency, show that
the selected frame learns the score covariance and that its extra
selection penalties vanish relative to the LIL threshold.
\begin{itemize}
 \item \textbf{Step 1: obtain simultaneous directional bounds.}
 Apply Ville's inequality with the allocated weights to the forecast
 score processes.
 \item \textbf{Step 2: turn directions into an ellipsoid.} Use the
 sphere net to control the whole whitened score, and map it through
 the design matrix.
 \item \textbf{Step 3: learn the covariance frame.} Under the efficiency
 assumptions, the forecast clock and its rounded frame converge to
 the true score covariance.
 \item \textbf{Step 4: identify the LIL radius.} The fine bet grid and
 lower-order matrix penalties yield the sharp scalar threshold in
 that frame, and hence the claimed sandwich shape.
\end{itemize}

\noindent\textbf{Step 1: obtain simultaneous directional bounds.}
Allocate error $\alpha\zeta_k/(\mathsf K_k\mathsf J_k)$ to each matrix and net-direction
pair, and split this over the dense bets with masses $\pi_{j,h}$.
Apply Ville's inequality to \eqref{eq:lb-forecast-process} with
$u=\lambda H^{-1/2}v$, where $v$ is a unit net direction. The cap is
respected since $sL\lambda/\sqrt{\lambda_{\min}(H)}\le1/5$.
Simultaneously for every time and every such tuple, except on an event of
probability at most $\alpha\sum_{k\ge2}\zeta_k\le\alpha$,
\begin{align*}
 v^\top H^{-1/2}S_n(\theta_*)
 &\le\lambda^{-1}\log\frac{\mathsf K_k\mathsf J_k}{\alpha\zeta_k\pi_{j,h}}
       +\lambda v^\top H^{-1/2}R_nH^{-1/2}v/2
       +s^3\lambda^2\sum_{i\le n}|v^\top H^{-1/2}x_i|^3\\
 &\le B_n(H,\lambda).
\end{align*}
\smallskip
\noindent\textbf{Step 2: turn directions into an ellipsoid.}
The sphere-net inequality
$\|z\|\le(1-1/k)^{-1}\max_v v^\top z$ proves
\eqref{eq:lb-cs-ellipsoid}. Selecting a matrix or bet using the data is
permitted because this event covers every selection.

\smallskip
\noindent\textbf{Step 3: learn the covariance frame.}
Now impose the efficiency assumptions. Bounded actions and the
conditional sub-Gaussian fourth
moment bound imply the martingale strong law
$n^{-1}\sum x_ix_i^\top(\varepsilon_i^2-\E_{i-1}\varepsilon_i^2)\to0$.
Consistency and Cauchy--Schwarz remove the forecast-error terms, so
$R_n/n\to K_\infty$. The rounding error of $H_n^\circ$ is at most
$5d/(2k_n)$ in operator norm. Its membership test eventually passes;
thus $H_n\to K_\infty$, $v_n(H_n)/n\to1$, and $c_{H_n}$ converges to a
finite positive constant.

\smallskip
\noindent\textbf{Step 4: identify the LIL radius.}
The covariance frame now has asymptotically unit variation, so the
remaining calculation is the scalar fine-grid argument. The two caps
are eventually inactive and
$\lambda_n\sim\sqrt{2\log\log n/n}$. Its dense index is of order
$\log n$, whence
$-\log\pi_{j,h}=\log j+3\log\log j+O(1)
\sim\log\log n$.
The additional matrix and net penalties are
$O_d(\log k_n)=O_d(\log\log\log n)=o(\log\log n)$.
Consequently $\log\frac{\mathsf K_{k_n}\mathsf J_{k_n}}{\alpha\zeta_{k_n}\pi_{j,h}}\sim\log\log n$ and
$B_n(H_n,\lambda_n)\sim\sqrt{2n\log\log n}$.
The support matrix of the ellipsoid is
$B_n^2G_n^{-1}H_nG_n^{-1}/(1-1/k_n)^2$, proving the assertion.
\end{proof}

In particular, for homoskedastic noise of variance $\sigma^2$,
$K_\infty=\sigma^2D_\infty$, and the limiting confidence-sequence shape is
\[
 (\theta-\widehat\theta_n)^\top G_n(\theta-\widehat\theta_n)
 \le(1+o(1))2\sigma^2\log\log n.
\]
For example, the predictable ridge least-squares forecast
$(\rho I_d+G_n)^{-1}\sum_{i\le n}x_iY_i$, for any fixed $\rho>0$,
is strongly consistent under the stable full-rank design assumption:
the martingale score divided by $n$ tends to zero. No bound on
$\|\theta_*\|$ is required for forecasting or for coverage.
The stability conditions concern efficiency alone; the simultaneous
coverage statement permits arbitrary predictable designs.

\subsubsection{Practical experts, hybrids, and computation}\label{app:lb-practical}

The proof grids are not enumerated, but their logarithmic penalties can be
large at moderate sample sizes. We therefore supplement the sharp sequence
with two inexpensive experts. The standard normal-mixture score sequence is
\begin{equation}\label{eq:lb-ay-score}
 \mathcal C_n^{\rm AY}(\delta)=\left\{\theta:
 S_n(\theta)^\top(I_d+G_n)^{-1}S_n(\theta)
 \le s^2\{\log\det(I_d+G_n)+2\log(1/\delta)\}\right\}.
\end{equation}
It is the natural Abbasi-Yadkori score mixture, written around ordinary least
squares rather than a ridge center, and requires no parameter-norm bound.
At nonsingular $G_n$ its shape matrix is
$s^2\{\log\det(I_d+G_n)+2\log(1/\delta)\}G_n^{-1}(I_d+G_n)G_n^{-1}$.

\paragraph{A cheap expert for moderate sample sizes.}
The sharp sequence's slowly refining grids can be supplemented by a fixed
frame with a small number of bets. A deterministic $1/4$-net of the unit
sphere has at most $9^d$ points. For
$\Lambda=\{(5sL)^{-1}2^{-h}:h=0,1,2,3\}$, define
\begin{equation}\label{eq:lb-coarse-expert}
 B_n^{\rm coarse}(\delta)=\frac43\min_{\lambda\in\Lambda}
 \left\{\frac{\log(4\cdot9^d/\delta)}\lambda
       +\frac\lambda2\lambda_{\max}(R_n)
       +s^3L\lambda_{\max}(G_n)\lambda^2\right\}.
\end{equation}
Then $\|S_n(\theta_*)\|\le B_n^{\rm coarse}(\delta)$ for every $n$
with probability at least $1-\delta$. Indeed, for each unit net direction
$v$,
\[
 \sum_{i\le n}|v^\top x_i|^3
 \le L\sum_{i\le n}(v^\top x_i)^2
 \le L\lambda_{\max}(G_n).
\]
Ville's inequality for the four bets in
\eqref{eq:lb-forecast-process}, a union bound over the net, and its
$4/3$ interpolation factor prove the claim. No net is enumerated. At a
nonsingular $G_n$ this expert is the ellipsoid with center
$\widehat\theta_n$ and shape
$(B_n^{\rm coarse})^2G_n^{-2}$.
\paragraph{Combine the sequential experts.}
The error split below protects every expert simultaneously. On this
common coverage event, choosing the smallest-volume ellipsoid preserves
coverage.
For a single practical confidence sequence, give the sharp sequence
error $\alpha/4$, this cheap expert error $\alpha/4$, and the
normal-mixture sequence \eqref{eq:lb-ay-score} error $\alpha/2$; return
the centered ellipsoid of smallest volume. The union bound proves its
simultaneous coverage. Under \eqref{eq:lb-cs-stability} the cheap expert
has a nonvanishing support radius, while the normal-mixture radius has
order $\sqrt{\log n/n}$. The sharp sequence is therefore eventually
selected, preserving the true-variance LIL limit.
Until $G_n$ has full rank, report $\R^d$ for this hybrid sequence.

\paragraph{A practical fixed-time hybrid retaining exact Gaussian efficiency.}
The coarse forecast expert can also improve the onset of the planned-time
construction. Let $\mathcal B_n(\delta)$ denote the score-ball ellipsoid
in \eqref{eq:lb-coarse-expert}, using the full-sample ordinary least-squares
center. At a planned sample size $n$, put
\begin{equation}\label{eq:lb-planned-hybrid}
 \delta_n=\alpha/\log(n+e),\qquad b_n=\alpha-\delta_n,
 \qquad
 \mathcal L_n^{\rm hyb}
 =\text{the smaller-volume member of }
       \{\mathcal L_{n,b_n},\mathcal B_n(\delta_n)\}.
\end{equation}
For these volume comparisons, assign any unbounded set infinite volume and
break ties by a fixed convention.
Use $\mathcal L_{n,b_n}$ if the full-sample design is singular.
Again, selection causes no additional error: on the intersection of
the two component coverage events, either reported set contains the
parameter. The union bound gives fixed-time coverage at least $1-\alpha$, irrespective
of which component is selected; their centers need not coincide.
Under \eqref{eq:lb-stability}, $b_n\to\alpha$ and the proof of
Theorem~\ref{thm:lb-fixed} applies unchanged to the heat component.
That component is eventually selected: its volume is of order $n^{-d/2}$,
whereas the coarse expert's four fixed bets have positive cubic charges
proportional to $\lambda_{\max}(D_n)$ and hence nonvanishing parameter
radii under the stable full-rank design. Thus
\eqref{eq:lb-planned-hybrid} has the same exact Gaussian limit
\eqref{eq:lb-gaussian-limit}, with $b=\alpha$.
This sample-size-dependent error split defines a fixed-time confidence set;
the preceding confidence sequence is the appropriate procedure for
continuous monitoring.

At full rank, the reported regression sets are explicit ellipsoids; matrix arithmetic costs
$O(d^3)$ per report. The forecast-clock statistics and ridge forecasts update
in $O(d^2)$ per observation. These constructions establish confidence
certificates and their limiting shapes, not improved bandit regret constants.

\section{Moment-bound constructions and proofs}\label{app:moments}

This appendix supplies the constructions and proofs for
Section~\ref{sec:moments}. The main change from the sub-Gaussian case is
in the inversion, rather than the limiting width. The exponential process
still charges the observed squared residuals, but these residuals depend
on the candidate mean. We first restrict that candidate to a coarse valid
interval and then solve an explicit quadratic inequality. The additional
moment makes the repair negligible at both the Gaussian and LIL scales.

\subsection{Moment inputs and the power-repaired process}\label{app:mp-process}

Assume $M_p>0$; the case $M_p=0$ was settled in the main text. We use at
most a third moment to bound the exponential remainder. Conditional
H\"older's inequality converts the supplied $p$th moment into both this
remainder bound and a variance bound. Explicitly, define
\begin{equation}\label{eq:mp-inputs}
 q=\min(p,3),\quad r=M_p^{1/p},\quad M_q=r^q,\quad
 V=r^2,\quad D_q=3^{(1-q)/2},\quad \kappa=D_qM_q.
\end{equation}
Then $\E_{i-1}|X_i-\mu|^q\le M_q$ and
$\E_{i-1}(X_i-\mu)^2\le V$. Here $r$ is the envelope scale, $V$ pays
for coarse localization, and $\kappa$ pays for the exponential remainder.
None is an estimate of the true variance, which may be arbitrarily smaller
than $V$.

The scalar inequality behind the construction is
\begin{equation}\label{eq:mp-powerineq}
 e^{u-u^2/2}\le1+u+D_q(u_-)^q\le1+u+D_q|u|^q.
\end{equation}
For $u\ge0$, use $\log(1+u)\ge u-u^2/2$. For $u=-x<0$, the cubic
argument from Section~2 gives
$e^{-x-x^2/2}-1+x\le\min\{x^3/3,x\}$; splitting at $x=\sqrt3$
bounds this minimum by $D_qx^q$. Substituting $u=\lambda(X_i-\mu)$
and taking conditional expectations removes the linear term. Thus the
unrepaired factor has expectation at most
$1+\kappa|\lambda|^q\le e^{\kappa|\lambda|^q}$.
Multiplying the repaired factors proves the test-supermartingale claim
in \eqref{eq:mp-sm}. The same argument permits predictable bets and the
sharper repair $\log(1+\kappa|\lambda|^q)$; the formulas below use the
power upper bound. To test a candidate $m$, evaluate the expression in
\eqref{eq:mp-sm} with $\mu$ replaced by $m$.

\subsection{The planned-time recipe and its proof}\label{app:mp-ci}

The reported interval combines two branches with the same center. The
exponential branch can be useful at moderate sample sizes; the Gaussian
continuation branch supplies the exact limiting quantile. We describe
their calibrations separately before proving that their intersection has
both properties.

For $n\ge4$, use the pilot and scoring split of Section~2:
$k=\lfloor\sqrt n\rfloor$, $N=n-k$, pilot variance $\widehat v$, scoring
mean $z$, and scoring scatter $Q$. The heat argument uses only the
inputs $(q,M_q,V)$ in \eqref{eq:mp-inputs}. At error level
$\varepsilon$, first reserve a small part of the error for coarse
localization, then divide the remaining error between the two tails.
The following choices also specify the smoothing margin and variance
budget:
\begin{align}
 \ell&=\log(n+e),& \beta_H&=\varepsilon/(4\ell),&
 a&=(\varepsilon-\beta_H)/2,&\eta&=1/(8\ell),\notag\\
 R_H&=\sqrt{V/(N\beta_H)},&
 \tau&=\{3M_qN/(a\eta)\}^{2/q},&
 A_H&=N(\widehat v+V/\ell)+NR_H^2+2\tau.
 \label{eq:mp-heat-tuning}
\end{align}
Here $R_H$ is the coarse radius, $\tau$ keeps the continuation smooth,
and $A_H$ pays for the pilot variance, a vanishing slack, localization,
and smoothing. If the observed scoring variance fits the pilot budget,
$Q/N\le\widehat v+V/\ell$, set
\begin{equation}\label{eq:mp-heat-radius}
 r_H(\varepsilon)=\min\left\{R_H,
 \frac{\sqrt{A_H}z_{1-a(1-2\eta)}
       +z_{1-\eta}\sqrt{A_H-Q}}{N}\right\};
\end{equation}
otherwise set $r_H(\varepsilon)=R_H$. The two terms in the numerator
are the terminal Gaussian cutoff and the allowance for remaining variance.
By Lemma~\ref{lem:ci-heatdrift}, the heat identity cancels the
quadratic drift and a reserve of $3M_q\tau^{-q/2}$ per scoring step
pays the remainder. No MGF is needed.

For the practical branch, localization prevents distant candidate means
from re-entering the exponential acceptance set. Indeed, replacing the
true mean by a candidate $m$ gives
$\sum_i(Z_i-m)^2=Q+N(z-m)^2$: the only extra term is the squared
distance from the reported center. On a coarse localization event this
term has a known bound, and the smaller quadratic root can sometimes
sharpen that bound.

Allocate a vanishing error budget to this practical branch so that the
Gaussian branch asymptotically receives the full level $\alpha$. Set
\begin{equation}\label{eq:mp-ci-levels}
 \delta_n=\frac{10\alpha}{10+\log(n+e)},\qquad
 \alpha_H=\alpha-\delta_n,\qquad
 \beta=\delta_n/4,\qquad L=\log\frac{2}{\delta_n-\beta},
 \qquad R=\sqrt{\frac{V}{N\beta}}.
\end{equation}
The part $\beta$ pays for localization; the remainder pays for the two
signed exponential tests. Choose the pilot-measurable bet
\begin{equation}\label{eq:mp-ci-bet}
 \lambda=\min\left\{\frac1r,
 \sqrt{\frac{2L}{N(\widehat v+V/k)}},
 \left(\frac{L}{\kappa N}\right)^{1/q}\right\},
 \qquad
 B=\frac{L}{N\lambda}+\frac{\lambda Q}{2N}
                  +\kappa\lambda^{q-1}.
\end{equation}
The square-root choice balances the evidence threshold against the pilot
variance, while the power cap prevents the repair from dominating when
the pilot variance is small. The basic practical radius is
$r_P=\min\{R,B+\lambda R^2/2\}$. If $\lambda R\le1$ and
$2\lambda B<1$, also take the smaller-root bound:
\begin{equation}\label{eq:mp-ci-root}
 r_P=\min\left\{r_P,
       \frac{2B}{1+\sqrt{1-2\lambda B}}\right\}.
\end{equation}
This is an explicit quadratic formula, without root finding. For $n\ge4$,
report
\begin{equation}\label{eq:mp-ci-final}
 I_n^p=z\ \pm\ \min\{r_H(\alpha_H),r_P\}.
\end{equation}
For $n<4$, use $\bar X_n\pm\sqrt{V/(n\alpha)}$.

\begin{proof}[Proof of Theorem~\ref{thm:mp-ci}]
The finite-sample argument uses only the conditional moment model. The
last two steps specialize to one fixed iid law with variance
$\sigma^2>0$.
\begin{itemize}
 \item \textbf{Step 1: localized exponential coverage.} Use a coarse
 radius to select the valid part of the quadratic inversion.
 \item \textbf{Step 2: Gaussian-branch coverage.} The heat reserve and
 the coarse localization event account for the two sources of error.
 \item \textbf{Step 3: the Gaussian limit.} Consistent variance estimates
 make the remaining-variance term negligible.
 \item \textbf{Step 4: the reported endpoints.} The practical branch is
 eventually wider, so the intersection inherits the Gaussian limit.
\end{itemize}

\noindent\textbf{Step 1: localized exponential coverage.}
Condition on the pilot. Conditional Chebyshev gives $|z-\mu|\le R$
outside an event of probability $\beta$. Markov's inequality for the
two signs of \eqref{eq:mp-sm}, scored with the pilot-chosen bet, gives
outside an event of probability $\delta_n-\beta$,
\[
 h-\frac\lambda2h^2\le B,\qquad h=|z-\mu|.
\]
Here $\sum_i(Z_i-\mu)^2=Q+Nh^2$. On $h\le R$, this implies
$h\le B+\lambda R^2/2$. Under the refinement conditions, the function
$h\mapsto h-\lambda h^2/2$ increases on $[0,R]$, and the larger
quadratic root exceeds $1/\lambda\ge R$. Therefore the inequality
forces $h$ below the smaller root in \eqref{eq:mp-ci-root}. The
practical branch has conditional error at most $\delta_n$.

\medskip\noindent\textbf{Step 2: Gaussian-branch coverage.}
For a generic level $\varepsilon$, apply the guarded drift bound in
Lemma~\ref{lem:ci-heatdrift} with the calibration above. The resulting
heat process has initial potential
$a(1-2\eta)$ and initial drift reserve
$3M_qN\tau^{-q/2}=a\eta$. Its initial wealth is thus $a(1-\eta)$.
At the terminal time the reserve is zero. For each sign, conditional
Markov bounds the probability that its potential exceeds $1-\eta$ by
$a$. Conditional Chebyshev pays $\beta_H$ for the coarse radius
$R_H$. On the variance check, every candidate in that coarse interval
has total squared residuals at most $A_H-2\tau$, so its process never
freezes. The same deterministic inversion as in
Appendix~\ref{app:scalar-heat} gives \eqref{eq:mp-heat-radius}.
If the check fails, the coarse radius applies directly. Total
miscoverage is at most $2a+\beta_H=\varepsilon$.
Taking $\varepsilon=\alpha_H$ and intersecting with Step 1 gives error
at most $\alpha_H+\delta_n=\alpha$. The $n<4$ claim follows directly
from Chebyshev's inequality.

\medskip\noindent\textbf{Step 3: the Gaussian limit.}
Now assume iid sampling. The moment of order $q>2$ gives polynomial
almost-sure rates for the pilot and scoring variances by
Lemma~\ref{lem:ci-lawrate}. Hence $Q/N-\widehat v=o(1/\ell)$, and the
heat check eventually passes. Moreover
\[
 \tau/N=O(N^{2/q-1}\ell^{2/q})=o(1/\ell),\qquad
 R_H^2=o(1/\ell),\qquad A_H/N\to\sigma^2,\qquad
 (A_H-Q)/N\sim V/\ell.
\]
The last relation identifies the remaining variance as vanishing slack,
rather than a leading variance charge. Since $\alpha_H\to\alpha$,
$a(1-2\eta)\to\alpha/2$, and
$z_{1-\eta}=O(\sqrt{\log\ell})$, the Gaussian cutoff in
\eqref{eq:mp-heat-radius} tends to the desired constant while the
remaining-variance term vanishes. Also
$\sqrt N R_H=\sqrt{V/\beta_H}\to\infty$, so clipping by $R_H$ does
not affect the limit. Thus
$\sqrt N r_H(\alpha_H)\to\sigma z_{1-\alpha/2}$ almost surely.

\medskip\noindent\textbf{Step 4: the reported endpoints.}
The practical branch cannot change this limit. Its threshold has
$L\to\infty$, while arithmetic--geometric mean gives
\[
 B\ge\frac{L}{N\lambda}+\frac{\lambda Q}{2N}
   \ge\frac{\sqrt{2LQ}}{N}.
\]
Both the base outer bound and the smaller root, when used, are at least
$B$. Since $Q/N\to\sigma^2>0$ and
$\sqrt N R=\sqrt{V/\beta}\to\infty$, it follows that
$\sqrt N r_P\to\infty$. The minimum therefore eventually selects the
heat branch. Finally $N/n\to1$ and
$\sqrt n(z-\bar X_n)=(k/\sqrt n)(z-\bar X_k)\to0$ almost surely,
which transfers the radius limit to the two stated endpoints.
\end{proof}

\subsection{The confidence-sequence recipe and its proof}\label{app:mp-cs}

Continuous monitoring requires one localization event and one exponential
event that hold at every time. The first can be coarse: it only needs to
shrink to zero. The second covers a calibrated grid of bets, allowing
current observations to select a useful bet without invalidating the
guarantee. We then apply the same localized quadratic inversion as above.

For the coarse event, divide time into dyadic blocks and allocate summable
error probabilities across them. Let $j(n)=\lfloor\log_2n\rfloor$ and
\begin{equation}\label{eq:mp-localizer}
 \beta_j=\frac\alpha2\frac{6}{\pi^2(j+1)^2},\qquad
 R_n=\frac{\sqrt{V2^{j(n)+1}/\beta_{j(n)}}}{n}.
\end{equation}
The process $S_t=\sum_{i\le t}(X_i-\mu)$ is a square-integrable
martingale, with $\E S_T^2\le VT$. Doob's maximal inequality gives
\[
 \Pp\!\left(\max_{t\le2^{j+1}}|S_t|
      >\sqrt{V2^{j+1}/\beta_j}\right)\le\beta_j.
\]
Summing over $j$ proves
$\Pp(\forall n:\ |\bar X_n-\mu|\le R_n)\ge1-\alpha/2$.
In particular, $R_n^2=O((\log n)^2/n)\to0$, although this coarse
radius is not LIL-optimal.

For the exponential event, use the dense grid of Section~3 with cap
$b=1/(5r)$:
\begin{align}
 m_j&=\lceil\log(j+e)\rceil,&
 w_j&=(\log2)\{1/\log(j+2)-1/\log(j+3)\},\notag\\
 \lambda_{j,k}&=b2^{-j-k/m_j},&
 \pi_{j,k}&=w_j/m_j,\qquad 0\le k<m_j.
 \label{eq:mp-grid}
\end{align}
As the bet size decreases, the grid becomes fine enough to avoid a fixed
inflation in the LIL constant. Four additional fixed bets,
$b,b/2,b/4,b/8$, improve finite-sample behavior. Give each fixed bet
mass $1/5$ and the dense grid the remaining mass $1/5$; the dense masses
$\pi_{j,k}$ sum to one within that grid. Allocate total error $\alpha/2$
to this entire two-sided betting family. The thresholds are
\begin{equation}\label{eq:mp-grid-levels}
 L_{j,k}=\log\frac{20}{\alpha\pi_{j,k}},\qquad
 L_h=\log\frac{20}\alpha\quad(h=0,1,2,3).
\end{equation}
Ville's inequality and a union bound give one event covering every bet
at every time. This simultaneous event justifies the data-dependent
selection below.

The selector balances the evidence threshold with observed variation,
using the coarse radius to bound the candidate-dependent part of the
quadratic. A power cap also controls the repair when the observed
variation is small. Explicitly, put
\begin{align}
 V_n^\circ&=Q_n+nR_n^2,&
 g_n&=\log\log(e^e+n)+\log(20/\alpha),\notag\\
 \lambda_n^*&=\min\left\{b,\sqrt{2g_n/V_n^\circ},
                (g_n/(2\kappa n))^{1/q}\right\},&
 u_n&=\log_2(b/\lambda_n^*),\notag\\
 j_n&=\lfloor u_n\rfloor,&
 k_n&=\lfloor m_{j_n}(u_n-j_n)\rfloor.
 \label{eq:mp-grid-select}
\end{align}
Write $\lambda_n=\lambda_{j_n,k_n}$. For each of the five selected
pairs $(\lambda,L)$---the dense choice and the four fixed bets---define
\[
 B_n(\lambda,L)=\frac{L}{n\lambda}
             +\frac{\lambda Q_n}{2n}+\kappa\lambda^{q-1},
 \qquad
 U_n(\lambda,L)=\min\{R_n,B_n(\lambda,L)+\lambda R_n^2/2\}.
\]
When $\lambda R_n\le1$ and $2\lambda B_n(\lambda,L)<1$, also take the
minimum with
\[
 \frac{2B_n(\lambda,L)}{1+\sqrt{1-2\lambda B_n(\lambda,L)}}.
\]
Let $r_n^p$ be the minimum of the resulting five values of $U_n$.
Each term is explicit; no search over candidate means or grid points
is performed.

\begin{proof}[Proof of Theorem~\ref{thm:mp-cs}]
Coverage follows from simultaneous events under the conditional moment
model. To identify the leading width, we then specialize to iid data
and show that both localization and the power repair are negligible.
\begin{itemize}
 \item \textbf{Step 1: simultaneous coverage.} Intersect the common
 localization and betting events, then invert on their intersection.
 \item \textbf{Step 2: the selected bet.} The dense grid approximates
 the true-variance LIL bet with vanishing relative error.
 \item \textbf{Step 3: the reported width.} The repair and localization
 terms disappear at first order, and the fixed bets cannot change the
 limit.
\end{itemize}

\noindent\textbf{Step 1: simultaneous coverage.}
The localizer event and the two-sided grid event for \eqref{eq:mp-sm}
have total error at most $\alpha$. On their intersection, the practical
CI argument gives
$|\bar X_n-\mu|\le U_n(\lambda,L)$ simultaneously for every calibrated
bet and every $n$. This includes all five selected pairs. Taking their
minimum proves coverage, even though the dense pair depends on current
observations.

\medskip\noindent\textbf{Step 2: the selected bet.}
At a fixed iid law, the strong law gives $Q_n/n\to\sigma^2$, and
$R_n\to0$ gives $V_n^\circ/n\to\sigma^2$. The square-root term in
\eqref{eq:mp-grid-select} eventually lies below the fixed cap. Its ratio
to the power cap is
$O((\log\log n/n)^{1/2-1/q})\to0$, so that cap is eventually inactive
as well. Rounding to the dense grid changes the bet by a factor smaller
than $2^{1/m_{j_n}}$, which tends to one. Consequently
\[
 \lambda_n\sim\sqrt{\frac{2\log\log n}{n\sigma^2}},\qquad
 L_{j_n,k_n}=\log\log n+3\log\log\log n+O(1).
\]
For the threshold relation, use
$j_n\sim\log n/(2\log2)$,
$w_j\sim(\log2)/[j(\log j)^2]$, and $m_j\sim\log j$.
Thus the cost of choosing a bet is exactly of LIL order, with only
lower-order overhead from the grid refinement.

\medskip\noindent\textbf{Step 3: the reported width.}
The first two terms in $B_n(\lambda_n,L_{j_n,k_n})$ are each asymptotic
to half the target scale. The repair is smaller because $q>2$:
\[
 \kappa\lambda_n^{q-1}
 =O\!\left((\log\log n/n)^{(q-1)/2}\right)
 =o\!\left(\sqrt{\log\log n/n}\right).
\]
The additional localization term $\lambda_nR_n^2/2$, relative to the
target scale, is $O(R_n^2)\to0$. The quadratic refinement, when used,
is asymptotic to $B_n$ because $2\lambda_nB_n\to0$. Meanwhile the
coarse radius $R_n$ is larger than the target by a diverging factor.
Therefore the dense branch has the limit in \eqref{eq:mp-lil}.

Finally, each fixed positive bet has
$B_n\to\lambda\sigma^2/2+\kappa\lambda^{q-1}>0$. Since $R_n\to0$,
its clipped radius eventually equals $R_n$; the quadratic refinement
also tends to a positive value whenever it is used. The four fixed bets
therefore cannot make the eventual minimum smaller than the dense-branch
scale. This proves \eqref{eq:mp-lil} for the reported radius.
\end{proof}

\subsection{A counterexample to the unrestricted sample-scatter process}
\label{app:mp-obstruction}

A moment bound permits a rare observation to be arbitrarily large.
If every observation in a sample takes that same large value, the sample
scatter is zero although the sum is large. Its exponential contribution
can overwhelm the polynomially small probability of this event. This
explains why the preceding constructions retain candidate-centered
residuals and localize before inversion.

For a formal counterexample, fix $n$, $\lambda>0$, and a finite
deterministic repair $C$. For large $A$, let
$\Pp(Y=A)=\varepsilon_A=M_p/(2A^p)$ and put the remaining mass at
$-\varepsilon_AA/(1-\varepsilon_A)$. These laws are centered. Their
positive atom contributes $M_p/2$ to the absolute $p$th moment, and the
other contribution tends to zero, so the total is at most $M_p$ for
sufficiently large $A$. Under iid sampling, the event
$Y_1=\cdots=Y_n=A$ has $Q_n=0$, giving
\[
 \E\exp\{\lambda\textstyle\sum_{i=1}^nY_i-\lambda^2Q_n/2-C\}
 \ge e^{-C}(\varepsilon_Ae^{\lambda A})^n\longrightarrow\infty.
\]
Reflection handles $\lambda<0$. Thus no finite deterministic repair can
make this unrestricted pure-$Q_n$ exponential valid over the moment
class. This obstruction concerns that particular exponential form; the
localized inversion above still gives explicit intervals with both
sharp limiting widths.

\subsection{Student moment envelopes}\label{app:mp-experimental-inputs}

For completeness, a Student variable with $\nu>2$ degrees of freedom,
rescaled to have variance one, has absolute central moment
\begin{equation}\label{eq:mp-student-envelope}
 M_p=(\nu-2)^{p/2}
 \frac{\Gamma((p+1)/2)\Gamma((\nu-p)/2)}
 {\sqrt\pi\,\Gamma(\nu/2)},\qquad p<\nu.
\end{equation}
The restriction $p<\nu$ is exactly the condition for this moment to be
finite.

\end{document}